\documentclass{article}

\usepackage{PRIMEarxiv}

\usepackage[utf8]{inputenc} 
\usepackage[T1]{fontenc}    
\usepackage{hyperref}       
\usepackage{url}            
\usepackage{booktabs}       
\usepackage{amsfonts}       
\usepackage{nicefrac}       
\usepackage{microtype}      
\usepackage{lipsum}
\usepackage{fancyhdr}       
\usepackage{graphicx}       
\graphicspath{{media/}}     
\usepackage{amssymb,mathrsfs,dsfont, mathtools,amsthm}
\usepackage{bbm}
\usepackage{bm,url,wasysym}
\usepackage{amsmath}
\usepackage{mleftright}
\usepackage[round]{natbib}
\usepackage{enumitem}
\usepackage{xcolor}
\usepackage{longtable}

\theoremstyle{plain}
\newtheorem{theorem}{Theorem}[section]
\newtheorem{lemma}[theorem]{Lemma}
\newtheorem{corollary}[theorem]{Corollary}

\theoremstyle{remark}
\newtheorem{Def}{Definition}[section]
\newtheorem{assumption}{Assumption}[section]
\newtheorem{remark}{Remark}[section]

\DeclareMathOperator*{\argmax}{arg\,max}
\DeclareMathOperator*{\argmin}{arg\,min}

\newcommand{\Bcal}{\mathcal{B}}
\newcommand{\Ccal}{\mathcal{C}}
\newcommand{\Fcal}{\mathcal{F}}
\newcommand{\Gcal}{\mathcal{G}}
\newcommand{\Hcal}{\mathcal{H}}

\newcommand{\Mcal}{\mathcal{M}}
\newcommand{\Ncal}{\mathcal{N}}

\newcommand{\Scal}{\mathcal{S}}
\newcommand{\Ucal}{\mathcal{U}}
\newcommand{\Wcal}{\mathcal{W}}
\newcommand{\Xcal}{\mathcal{X}}

\newcommand{\Zcal}{\mathcal{Z}}

\newcommand{\Ebb}{\mathbb{E}}

\newcommand{\Nbb}{\mathbb{N}}
\newcommand{\Pbb}{\mathbb{P}}
\newcommand{\Rbb}{\mathbb{R}}

\newcommand{\ind}{\mathds{1}}

\title{Nonparametric Intensity Estimation for Covariate-Driven Poisson Processes Using Deep Learning
}

\author{
  Atsutomo Yara \\
  Graduate School of Engineering Science, \\
  the University of Osaka \\
  \texttt{a-yara@sigmath.es.osaka-u.ac.jp} \\
   \And
  Yoshikazu Terada \\
  Graduate School of Engineering Science, \\
  the University of Osaka \\
  RIKEN Center for Advanced Intelligence Project (AIP) \\
  \texttt{terada.yoshikazu.es@osaka-u.ac.jp} \\
}

\begin{document}
\maketitle

\begin{abstract}
    We consider nonparametric maximum likelihood estimation of the intensity function of a Poisson process, 
    allowing the intensity to depend on observed covariates. 
    Such models describe how event rates vary with external conditions, 
    including conditions under which the intensity becomes arbitrarily small or vanishes. 
    However, risk analysis based on the Kullback–Leibler (KL) divergence typically requires the true intensity 
    to be bounded away from zero, a restrictive assumption in applications. 
    We develop a unified approach that avoids this assumption and derive an oracle inequality for the covariate-averaged Hellinger distance between the true intensity and the nonparametric maximum likelihood estimator (NPMLE).
    As an application of the oracle inequality, we establish convergence rates for NPMLEs based on deep neural networks (DNNs). Our results show that, as in nonparametric regression, 
    the ability of DNNs to capture underlying structure leads to improved convergence rates. 
    Specifically, DNN-based NPMLEs mitigate the curse of dimensionality when the true intensity has a compositional structure and automatically adapt to a low-dimensional Riemannian manifold structure of the event space. 
    For intensity functions with a compositional structure, 
    DNN-based NPMLEs achieve nearly minimax-optimal convergence rates 
    without explicitly incorporating this structure into the network architecture.
    Simulations demonstrate that DNN-based estimators outperform kernel estimators in most settings, 
    while manifold kernels remain competitive in some designs. 
    Finally, we analyze stimulus-evoked neural spiking using data from the Allen Brain Observatory Visual Coding Neuropixels project. 
    DNN-based NPMLEs provide better out-of-sample performance than a conditional kernel estimator, 
    as evaluated by held-out log-likelihood.
\end{abstract}


\section{Introduction} \label{sec:intro}
Poisson (point) processes are fundamental probabilistic models for random events occurring in time and space, 
and they have been widely applied in fields such as epidemiology, ecology, geology, seismology, and urban informatics.
Poisson processes (and their generalizations) have been used, for example, to model the occurrence of earthquakes \cite{Ogata98}, the spatial distribution of disease cases as a function of environmental exposure \cite{Diggle97}, 
and the frequency of taxi demand in urban areas \cite{Moreira-MatiasEtAl13}.
A Poisson process is characterized by its intensity function, 
which represents the mean occurrence rate of events at any point in time or space.
Estimating this intensity function from event data is an important problem, 
as it provides insights into the underlying event-generation mechanism and enables the prediction of future events.

Methods for estimating the intensity function can be broadly categorized according to the type of data available.
In some applications, only a single realization of a Poisson process is observed; 
that is, we can observe only one sequence of event times generated by the underlying process.
Examples include the recorded occurrence times of earthquakes in a region over a century 
or data on the failure time of a system observed only once.
In such cases, because only a single sample path of the underlying process is available, 
estimation of the intensity function is challenging without some structural assumptions.
In this situation, the intensity function is often modeled parametrically, 
and its parameters are estimated by maximum likelihood estimation.
Since such parametric approaches are interpretable, we can gain insights into the underlying structure.
Moreover, it is well known that the estimators have desirable theoretical properties, such as consistency and asymptotic normality, which can be established under suitable regularity conditions (\cite{Kutoyants12}).
However, when modeling complex phenomena in real-world applications, 
the true intensity function is often too complicated to be captured by simple parametric models.
In such situations, nonparametric estimation of the intensity function becomes important.
Classic examples of nonparametric estimators include kernel smoothing methods \cite{Diggle85}.
Recently, more sophisticated estimators have been proposed, such as those based on reproducing kernel Hilbert spaces \cite{FlaxmanEtAl17} and estimators for Poisson processes supported on Riemannian manifolds \cite{WardEtAl23}.

In other settings, 
one may observe multiple independent realizations of a Poisson process obtained 
under comparable experimental or environmental conditions.
Typical examples include spike sequences collected from repeated experimental trials or temporal event data 
from multiple patients.
Because these settings provide substantially more observations than the single-realization case, 
nonparametric intensity estimation becomes considerably more powerful.
A wide range of estimators has been proposed for this data type, 
including those based on thresholding rules \cite{Reynaud-BouretRivoirard10}, 
kernel smoothing techniques \cite{ChiangEtAl05}, 
and hierarchical Bayesian formulations \cite{KirichenkoZanten15}, 
accompanied by analyses of their convergence rates as the number of observed processes $n \to \infty$.

Covariate-driven point process models provide a flexible framework for describing heterogeneity across time, 
space, or experimental conditions.
In such models, the intensity function is allowed to depend on observed covariates, and this dependence is often interpreted through the framework of Cox processes, or doubly stochastic Poisson processes.
A number of estimation methods have been developed for covariate-dependent intensity functions, including parametric estimators (\citealp{Lawless87,O'Sullivan93}), kernel-smoothing-based methods \cite{KlutchnikoffMassiot16}, and more recent deep-learning-based approaches (\citealp{OkawaEtAl19, WangEtAl20}).
From a theoretical perspective, \cite{GiordanoEtAl2026} recently studied nonparametric Bayesian estimation for covariate-driven inhomogeneous Poisson point processes under growing-domain asymptotics.
They considered the setting in which a single spatial point pattern is observed together with a spatial covariate field $Z$, and the intensity is modeled as $\lambda(x)=\rho(Z(x))$.
Under suitable stationarity and ergodicity conditions on the covariate field, they established posterior contraction rates for several classes of priors, including Gaussian process priors, Gaussian mixture priors, and Pólya-tree-type priors.
In another related direction, \cite{GyotokuEtAl25} investigated deep learning for point processes whose intensities depend on mixing covariate processes.
Their framework covers Cox-type models, marked point processes, and multivariate point processes, and they derived an oracle inequality and convergence rate for the risk function.
In this work, we study the problem of estimating the true conditional intensity function $S_0(\cdot \mid \star)$ from i.i.d. observations of the form $(Z_1, \{x_{1, 1}, \dots, x_{1, N_1}\}), \dots, (Z_n, \{x_{n, 1}, \dots, x_{n, N_n}\})$, where each Poisson process is driven by the conditional intensity $S_0(\cdot \mid Z_i)$ given covariates $Z_i$.

Most existing theoretical studies of nonparametric intensity estimation focus on estimators with explicit forms, such as kernel, spline, or Bayesian estimators.
However, theoretical analysis remains limited for cases where the estimator lacks an explicit representation, such as deep-learning-based intensity estimation, despite its strong predictive performance in practice.
In particular, in nonparametric maximum likelihood estimation, standard generalization error analysis typically investigates the convergence of the Kullback-Leibler (KL) divergence between the true model and the estimated model.
However, such an approach could fail in our setting because the KL divergence diverges when the support of the estimated intensity function does not cover that of the true intensity function.
As a result, standard KL-based analysis becomes difficult unless additional assumptions are imposed on the true intensity function.
For example, it is often assumed that the true intensity function is bounded away from zero (e.g., \cite{KirichenkoZanten15}).
Under this assumption, the KL divergence reduces to the squared error, yielding convergence rates of the same form as in regression problems.
On the other hand, from the viewpoint of real-world applications, convergence in the KL divergence is not necessarily required.
Instead, it is often sufficient to establish consistency with respect to an appropriate metric.
Motivated by this, we apply the techniques of \cite{vandeGeer00} and derive an oracle inequality to evaluate the Hellinger distance rather than the KL divergence.
Using this oracle inequality, we can establish the consistency of the nonparametric maximum likelihood estimator (NPMLE) for the intensity function under realistic assumptions.
Moreover, as an analog of the $\alpha$-SVB condition for conditional probability estimation, which is introduced in \cite{BosSchmidt-Hieber22}, we show that when the distribution of the covariates and the class of the true intensity functions satisfy suitable conditions, the approximation error can be improved even if the true intensity function is not bounded away from zero.
As an important application of the oracle inequality, we derive convergence rates for the maximum likelihood estimator of the intensity function based on deep learning in the following two cases: (i) when the true intensity function has a compositional structure, and (ii) when the true intensity function is supported on a Riemannian manifold.
These results reveal that deep learning can avoid the curse of dimensionality in case (i), and can automatically adapt to the Riemannian manifold structure in case (ii), even without being given explicit information such as the Riemannian metric.
Furthermore, in case (i), we show that deep learning achieves a convergence rate that is nearly minimax optimal.

The main contributions of this study can be summarized as follows:
\begin{itemize}
\item Theoretical analysis of the NPMLE is difficult within the standard framework of generalization error analysis, due to the unboundedness of the KL divergence, unless one imposes assumptions on the true intensity function.
By adapting the density estimation argument of \cite{vandeGeer00}, we derive an oracle inequality that evaluates the performance of the NPMLE for the intensity function of a Poisson process under mild assumptions.
\item As an analog of the $\alpha$-SVB condition in conditional probability estimation, we introduce a corresponding condition for conditional intensity estimation and show that, under this condition, the approximation error can be improved.
\item When the true intensity function has a compositional structure, we derive convergence rates for the NPMLE based on deep learning. 
This shows that, in this setting, deep learning can avoid the curse of dimensionality.
We further show that the resulting convergence rates are nearly minimax optimal.
\item When the underlying Poisson process is supported on a low-dimensional Riemannian manifold, we show that the deep-learning-based NPMLE can automatically adapt to the low-dimensional structure.
\end{itemize}

\section{Preliminaries}
\subsection{Notation}
In this section, we introduce the notation used throughout the paper.
First, let $\Nbb$ denote the set of natural numbers, and let $\Rbb$ denote the set of real numbers.
For a metric space $(\Xcal,d)$, we write $\Bcal(\Xcal)$ for the Borel $\sigma$-algebra on $\Xcal$.
For two measures $\lambda$ and $\mu$ defined on the same measurable space, the notation $\lambda \ll \mu$ means that $\lambda$ is absolutely continuous with respect to $\mu$, that is, for every measurable set $B$, $\mu(B)=0$ implies $\lambda(B)=0$.
In this case, we write $d\lambda/d\mu$ for the density of $\lambda$ with respect to $\mu$, namely, the Radon--Nikodym derivative.
Following \cite{leskelä24}, for two $\sigma$-finite measures $\lambda$ and $\mu$ on the same measurable space, we define the Hellinger distance between $\lambda$ and $\mu$ by
\begin{equation*}
    H(\lambda, \mu) \coloneqq \left\{\frac{1}{2} \int (\sqrt{\smash[b]{f}} - \sqrt{\smash[b]{g}})^2 \, d\nu\right\}^{1/2},
\end{equation*}
where $f \coloneqq d\lambda/d\nu$ and $g \coloneqq d\mu/d\nu$ are the densities of $\lambda$ and $\mu$, respectively, with respect to a $\sigma$-finite measure $\nu$.
For two probability measures $P$ and $Q$ on the same probability space, we define the Kullback--Leibler (KL) divergence by
\begin{equation*}
    \mathrm{KL}(P \parallel Q) \coloneqq
    \begin{dcases*}
        \int \log\left(\frac{dP}{dQ}\right)\, dP & if $P \ll Q$, \\
        \infty & otherwise.
    \end{dcases*}
\end{equation*}
We also adopt the conventions $\log(0)\coloneqq -\infty$, $0\log(0)\coloneqq 0$, $1/0\coloneqq \infty$, and $0/0\coloneqq 0$.
For real numbers $a,b>0$, we write $a\lor b$ for the maximum of $a$ and $b$.
For a vector $v=(v_1,\dots,v_d)^\top\in\Rbb^d$ and a univariate function $f:\Rbb\to\Rbb$, we write
\[
f(v)\coloneqq (f(v_1),\dots,f(v_d))^\top.
\]
Throughout the paper, operations and inequalities involving vectors and vector-valued functions are interpreted componentwise.
For a vector $v\in\Rbb^d$, we define its $\ell_p$ norm by
\[
\lVert v\rVert_p \coloneqq \left(\sum_{i=1}^d \lvert v_i\rvert^p\right)^{1/p}.
\]
Moreover, for a measure $\mu$ on $\Rbb^d$ and a real-valued function $f:D\subset\Rbb^d\to\Rbb$, we define the $L^p$ norm by
\[
\lVert f\rVert_{L^p(D,\mu)} \coloneqq \left(\int_D \lvert f\rvert^p\,d\mu\right)^{1/p}.
\]
When the domain $D$ or the measure $\mu$ is clear from the context, we abbreviate this notation as $\lVert f\rVert_{L^p(D)}$, $\lVert f\rVert_{L^p(\mu)}$, or simply $\lVert f\rVert_{L^p}$.
We also define $\lVert f\rVert_\infty \coloneqq \sup_{x\in D} |f(x)|$.
Finally, for a real number $x\in\Rbb$, we write $\lceil x\rceil$ for the smallest integer greater than or equal to $x$, and $\lfloor x\rfloor$ for the largest integer less than or equal to $x$.

\subsection{Setting}
In this section, following \cite{LastPenrose17}, we introduce the general definitions of Poisson processes and Cox processes.
Let $(\Omega, \Fcal, \Pbb)$ be a probability space, and let $(\Xcal, \Bcal)$ be a measurable space.
We denote by $\Mcal_{< \infty} = \Mcal_{< \infty}(\Xcal)$ the space of all nonnegative integer-valued measures on $\Xcal$ such that $\mu(B) < \infty$ for all $B \in \Bcal$.
We further denote by $\Mcal = \Mcal(\Xcal)$ the space of all measures that can be written as countable sums of measures in $\Mcal_{< \infty}$, and by $\Ncal = \Ncal(\Xcal)$ the smallest $\sigma$-algebra on $\Mcal$ such that the mapping $\mu \mapsto \mu(B)$ is measurable for every $B \in \Bcal$.
We then define a point process on $\Xcal$ as follows.

\begin{Def}[Point process]
    A point process on $\Xcal$ is a random element of $(\Mcal, \Ncal)$.
    That is, it is a measurable mapping $X : \Omega \to \Mcal$.
\end{Def}

Intuitively, a point process describes a random collection of points on the space $\Xcal$.
Next, we define a Poisson process.

\begin{Def}[Poisson (point) process]
    Let $\lambda$ be a $\sigma$-finite measure on $\Xcal$.
    A Poisson (point) process on $\Xcal$ with intensity measure $\lambda$ is a point process $X$ on $\Xcal$ satisfying the following two properties:
    \begin{enumerate}[label=(\roman*)]
        \item For every $B \in \Bcal$, the random variable $X(B)$ follows the Poisson distribution with intensity parameter $\lambda(B)$.
        \item For every $N \in \Nbb$ and every collection of pairwise disjoint sets $B_1, \dots, B_N \in \Bcal$, the random variables $X(B_1), \dots, X(B_N)$ are independent.
    \end{enumerate}
\end{Def}

Hereafter, we write $P_\lambda$ for the distribution of a Poisson process $X$ with intensity measure $\lambda$.
That is, for any $N \in \Ncal$, we define $P_\lambda(N) \coloneqq \Pbb(X^{-1}(N))$.
For a Poisson process $X$ on $\Xcal$ with intensity measure $\lambda$, and for any $\sigma$-finite measure $\nu$ satisfying $\lambda \ll \nu$, the Radon--Nikodym derivative
\[
f(x) \coloneqq \frac{d\lambda}{d\nu}(x)
\]
is called the intensity function of $X$ with respect to $\nu$.
Intuitively, the intensity function describes how likely points are to occur at each location in $\Xcal$.
Next, we define a Cox process.

\begin{Def}[Cox process]
    Let $\xi$ be a random measure on $\Xcal$.
    A point process $\eta$ is called a Cox process directed by $\xi$ if
    \begin{equation*}
        \Pbb(\eta \in N \mid \xi) = P_\xi(N),\quad \Pbb\text{-a.s.},\; N \in \Ncal.
    \end{equation*}
    In this case, $\xi$ is called the directing random measure of $\eta$.
\end{Def}

In this study, given a covariate $Z$ taking values in $[0,1]^{d_z}$, we assume that $X(\cdot \mid Z)$ follows a Poisson process on a compact set $A \in \Bcal(\Rbb^{d_x})$ in the $d_x$-dimensional Euclidean space $\Rbb^{d_x}$, and we consider the problem of estimating the intensity function of $X(\cdot \mid Z)$.
That is, $X(\cdot \mid \cdot)$ is a Cox process whose directing random measure can be written as $\xi(\omega) = \xi(Z(\omega))$ for $\omega \in \Omega$, so that the randomness affects the process only through the covariate.

Hereafter, we set $d \coloneqq d_z + d_x$ and write $\lambda_0(\cdot \mid Z = z) \coloneqq \xi(z)$ for the directing random measure given the covariate value $Z=z$.
Using i.i.d.\ observations $(Z_1, X_1), \dots, (Z_n, X_n)$ from the Poisson process $X = X(\cdot \mid Z)$ with covariate $Z$ and intensity measure $\lambda_0(\cdot \mid Z)$, we study the problem of estimating the conditional intensity function of $X$ given $Z$, defined by
\[
S_0(x \mid z) \coloneqq \frac{d\lambda_0(\cdot \mid z)}{d\nu}(x).
\]
Here, $\nu$ is a measure on $A$, and we assume that $\lambda_0 \ll \nu$.
Throughout this paper, $\nu$ is basically taken to be the $d_x$-dimensional Lebesgue measure, while if $A$ is a low-dimensional manifold embedded in $\Rbb^{d_x}$, then $\nu$ is taken to be the volume element on that manifold.
We assume that $0\le\lambda_0(A\mid z)<\infty$.
If $\lambda_0(A\mid z)=0$, the corresponding Poisson process is
the empty point pattern almost surely.
Throughout, $dx$ denotes the $d_x$-dimensional Lebesgue measure.
Moreover, the observation corresponding to each $X_i$ given $z_i,\; i = 1, \dots, n$ is a random set of points $\{x_{i1}, \dots, x_{im_i}\} \subset A$, where the number of points $m_i = X_i(A \mid z_i)$ is a realization of a Poisson random variable with intensity parameter $\lambda_0(A \mid z_i)$.
We consider, as the true intensity measure $\lambda_0(\cdot \mid z)$ given the covariate $z$, a finite measure, and we assume that the true intensity function
is bounded.
Furthermore, we assume that the covariate $Z$ has a positive density on $[0, 1]^{d_z}$.

\begin{remark}
    The term ``conditional intensity function'' usually refers to the intensity function given the information available up to the previous time in the context of self-exciting stochastic processes.
    In this paper, however, we use the term ``conditional intensity function'' to mean the intensity function given the covariate.
\end{remark}

We briefly explain the significance of modeling a point process given covariates.
For example, suppose that we model the number of taxi pickups in a certain region by a Poisson process.
Then it is natural to think that the corresponding intensity function depends on covariates such as temperature and precipitation.
In this case, one could consider modeling the covariates such as temperature and precipitation themselves as part of a Poisson process, that is, modeling $X$ as a Poisson process on the Cartesian product of the two-dimensional spatial domain of taxi pickup locations and the covariate space.
However, this would imply that the data-generating mechanism produces not only taxi pickup locations but also temperature and precipitation as count data, which is clearly unnatural.
Therefore, in such situations, it is natural and practically useful to model the count data as a Poisson process given the covariates.

\subsection{Nonparametric Maximum Likelihood Estimation}
\subsubsection{Formulation of Nonparametric Maximum Likelihood Estimation}
In this section, we formulate the nonparametric maximum likelihood estimation of the intensity function.
We begin by reviewing the expression for the likelihood ratio between two Poisson processes.

\begin{lemma}[Theorem 1.3 of \citealp{Kutoyants12}] \label{lem: likelihood ratio}
    Let $\lambda$ and $\nu$ be finite intensity measures satisfying $\lambda \ll \nu$.
    Then the corresponding Poisson process distributions satisfy $P_\lambda \ll P_\nu$.
    Moreover, the likelihood ratio is given by
    \begin{equation*}
        \frac{dP_\lambda}{dP_\nu}(\eta)
        =
        \exp\left\{
            \int \log(\phi)\, d\eta
            -
            \int (\phi - 1)\, d\nu
        \right\},
    \end{equation*}
    where $\phi = d\lambda / d\nu$.
\end{lemma}

Now consider a Poisson process on a compact set $A \subset \Rbb^{d_x}$, and apply Lemma~\ref{lem: likelihood ratio}.
Then, writing the conditional intensity function given the covariate $z$ as
\[
S(x \mid z) = \frac{d\lambda(\cdot \mid z)}{d\nu}(x),
\]
we obtain
\begin{equation} \label{eq: likelihood ratio}
    \frac{dP_\lambda}{dP_\nu}(\eta)
    =
    \exp\left\{
        \int_A \log(S(x \mid z))\, d\eta
        -
        \int_A (S(x \mid z) - 1)\, d\nu(x)
    \right\}.
\end{equation}

Based on the likelihood ratio \eqref{eq: likelihood ratio}, for a given conditional intensity function $S$, we define the negative log-likelihood function for an observation $(X_i, z_i)$ by
\begin{equation*}
    \ell(X_i, z_i \mid S)
    \coloneqq
    - \int_A \log(S(x \mid z_i))\, dX_i(x \mid z_i)
    + \int_A S(x \mid z_i)\, d\nu(x)
    =
    - \sum_{j=1}^{N_i} \log(S(x_{ij} \mid z_i))
    + \int_A S(x \mid z_i)\, d\nu(x).
\end{equation*}
Here, $\{x_{i1}, \dots, x_{iN_i}\} \subset A$ denotes the observed points corresponding to $X_i$ given the covariate $z_i$, and $N_i = X_i(A \mid z_i)$.
Therefore, for i.i.d.\ observations $(X_1, z_1), \dots, (X_n, z_n)$ from the covariate $Z$ and the Poisson process $X(\cdot \mid Z)$ given $Z$, the negative log-likelihood function normalized by $1/n$ is given by
\begin{align}
    \ell_n(X_1, \dots, X_n, z_1, \dots, z_n \mid S)
    &\coloneqq
    \frac{1}{n} \sum_{i=1}^n \ell(X_i, z_i \mid S) \notag \\
    &=
    \frac{1}{n} \sum_{i=1}^n
    \left\{
        - \sum_{j=1}^{N_i} \log(S(x_{ij} \mid z_i))
        + \int_A S(x \mid z_i)\, d\nu(x)
    \right\}. \label{eq: negative log likelihood}
\end{align}
\[
\ell_n(S) \coloneqq \ell_n(X_1, \dots, X_n, z_1, \dots, z_n \mid S).
\]
In nonparametric maximum likelihood estimation, the estimator is defined as a function that minimizes the negative log-likelihood \eqref{eq: negative log likelihood} over a function space (hypothesis space) $\Fcal_n$, which is not necessarily parameterized by a finite-dimensional parameter.
The resulting estimator is called the nonparametric maximum likelihood estimator (NPMLE).
That is, the NPMLE $\hat{S}_n(x \mid z)$ of the conditional intensity function $S_0(x \mid z)$ of the Poisson process $X(\cdot \mid z)$ given the covariate $z$ is defined by
\begin{equation}
    \hat{S}_n \in \argmin_{S \in \Fcal_n} \ell_n(S). \label{eq: NPMLE}
\end{equation}

\subsubsection{Difficulty of Evaluating the Excess Risk} \label{sec: difficulty}
In the context of standard generalization error analysis, the performance of the NPMLE is typically evaluated by the excess risk
\begin{equation*}
    R(\hat{S}_n) \coloneqq \mathbb{E}[\ell(\hat{S}_n) \mid (Z_1, X_1), \dots, (Z_n, X_n)] - \inf_{S \in \Fcal_0}\mathbb{E}[\ell(S)],
\end{equation*}
where $\Fcal_0$ denotes the set of all nonnegative bounded measurable functions $f : A \times [0, 1]^{d_z} \to [0,\infty)$.
The excess risk depends on the estimator $\hat{S}_n$.
Thus, it is a random variable.
When no fonfusion is likely to arise from the context, we abbreviate $\mathbb{E}[\ell(\hat{S}_n) \mid (Z_1, X_1), \dots, (Z_n, X_n)]$ as $\mathbb{E}[\ell(\hat{S}_n)]$.
The following lemma gives an explicit expression for the excess risk.
Here, we write $\hat{\lambda}_n(\cdot \mid z)$ for the intensity measure corresponding to the estimator $\hat{S}_n(\cdot \mid z)$.

\begin{lemma} \label{lem: the excess risk}
    For the excess risk $R(\hat{S}_n)$, we have
    \begin{align}
        R(\hat{S}_n)
        &= \mathbb{E}_Z\!\left[\mathrm{KL}\!\left(P_{\lambda_0(\cdot \mid Z)} \parallel P_{\hat{\lambda}_n(\cdot \mid Z)}\right)\right] \nonumber \\
        &= \mathbb{E}_Z \left[
            \int_A S_0(x \mid Z)\log\left(\frac{S_0(x \mid Z)}{\hat{S}_n(x \mid Z)}\right)\, d\nu(x)
            + \int_A \bigl(\hat{S}_n(x \mid Z) - S_0(x \mid Z)\bigr)\, d\nu(x)
        \right], \label{eq: the excess risk}
    \end{align}
    where $\mathbb{E}_Z$ denotes expectation with respect to $Z$.
\end{lemma}

Equation \eqref{eq: the excess risk} diverges unless $\lambda_0(\cdot \mid z) \ll \hat{\lambda}_n(\cdot \mid z)$ on $A$ for every $z$.
Therefore, in order to evaluate the excess risk directly, one would need to assume, for example, that there exists a constant $c>0$ such that $\hat{S}_n > c$ for all estimators $\hat{S}_n$.
However, this is a rather strong and unrealistic assumption when discussing the convergence of the estimator to the true intensity function.
Moreover, from the viewpoint of real-world applications, convergence of the excess risk is not necessarily required; it is often sufficient to establish consistency with respect to an appropriate metric.
Therefore, by adapting the techniques of \cite{vandeGeer00}, we evaluate a lower bound of the excess risk $R(\hat{S}_n)$ instead of the excess risk itself, and establish consistency of that lower bound under realistic assumptions.
The following lemma shows that the Hellinger distance between intensity measures appears as a lower bound of the excess risk.
For notational simplicity, we also write the Hellinger distance between intensity measures in terms of the corresponding intensity functions, for example,
\[
H(S_1,S_2) \coloneqq H(\lambda_1,\lambda_2).
\]
We also define the expectation, with respect to the covariate, of the squared Hellinger distance between conditional intensity functions by
\begin{equation*}
    \widetilde{H}^2(S_1,S_2)
    \coloneqq
    \mathbb{E}_Z\bigl[H^2(S_1(\cdot \mid Z), S_2(\cdot \mid Z))\bigr].
\end{equation*}

\begin{lemma} \label{lem: lower bound of the excess risk}
    For the excess risk $R(\hat{S}_n)$, we have
    \begin{equation*}
        R(\hat{S}_n) \ge 2 \widetilde{H}^2(\hat{S}_n,S_0).
    \end{equation*}
\end{lemma}

\begin{remark}
    Since $x+1 \le e^x$ for any $x \in \mathbb{R}$, we have
    \begin{equation*}
        \widetilde{H}^2(\hat{S}_n,S_0)
        \ge
        \mathbb{E}_Z \bigl[1 - e^{-H^2(\hat{S}_n,S_0)}\bigr].
    \end{equation*}
    The quantity inside the expectation on the right-hand side is the Hellinger distance between the Poisson processes
    $P_{\hat{\lambda}_n}$ and $P_{\lambda_0}$ (see \citealp[Corollary 5.3]{leskelä24}).
    Thus, convergence in the Hellinger distance between intensity measures is stronger than convergence in the Hellinger distance between the corresponding Poisson processes.
\end{remark}

\section{Oracle Inequality}
In this section, by adapting the argument of \cite{vandeGeer00}, we establish an oracle inequality for evaluating the Hellinger distance between intensity measures.
We impose the following assumption on the true intensity function and the estimation model $\Fcal_n$.

\begin{assumption} \label{asm: weak absolute continuity}
    There exist a sequence of elements $S_n^* \in \Fcal_n$ and a constant $c_0 \ge 0$ such that
    \begin{equation*}
        \forall n \in \Nbb,\; \forall x \in A,\; \forall z \in [0, 1]^{d_z};\;
        \frac{S_0(x \mid z)}{S_n^*(x \mid z)} \leq c_0^2.
    \end{equation*}
\end{assumption}

\begin{remark}
    Assumption \ref{asm: weak absolute continuity} means that, if we suitably choose the hypothesis space $\Fcal_n$, we do not need to impose additional conditions on the true conditional intensity function $S_0$.
    We later show that the estimation model $\Fcal_n$ is taken to be a suitable class of deep neural networks; this assumption is automatically satisfied. See the proof of Theorem~\ref{thm: convergence rates for composition structured functions} for details.
\end{remark}

For a measure $Q$ and a measurable function $f$, we use the notation $Qf \coloneqq \int f\, dQ.$
We also write the population distribution as $P \coloneqq P_{\lambda_0 \mid Z} \otimes P_Z,$ and denote by $P_n$ the empirical distribution corresponding to the observations $(X_1, Z_1), \dots, (X_n, Z_n)$.
We first show, using the property of the maximum likelihood estimator, that the Hellinger distance between intensity measures can be bounded by a certain empirical process.
The following lemma is a key to deriving the oracle inequality.

\begin{lemma}[Basic Inequality] \label{lem: basic inequality}
    Suppose that Assumption \ref{asm: weak absolute continuity} holds.
    Then, for any $n \in \Nbb$ and any NPMLE $\hat{S}_n$, we have
    \begin{align*}
        \widetilde{H}^2\left(
            \frac{\hat{S}_n + S_n^*}{2},
            S_n^*
        \right)
        &\leq
        (P_n - P)
        \left\{
            \int_{A_n^*(Z)}
            \frac{1}{2}
            \log \left(
                \frac{\hat{S}_n(x \mid Z) + S_n^*(x \mid Z)}
                {2 S_n^*(x \mid Z)}
            \right)
            dX(x \mid Z)
        \right. \\
        &\quad\left.
            +
            \frac{1}{2}
            \int_{A}
            \left(
                S_n^*(x \mid Z)
                -
                \frac{\hat{S}_n(x \mid Z) + S_n^*(x \mid Z)}{2}
            \right)
            d\nu(x)
        \right. \\
        &\quad
        +
        2(c_0 + 1)
        \widetilde{H}\left(
            \frac{\hat{S}_n + S_n^*}{2},
            S_n^*
        \right)
        \widetilde{H}\left(
            S_n^*,
            S_0
        \right), 
    \end{align*}
    where $A_n^*(z) \coloneqq A \cap \{x \in A : S_n^*(x \mid z) > 0\}$ for $z \in [0, 1]^{d_z}$.
\end{lemma}

We denote the index sets of the empirical processes by
\begin{align*}
    \Gcal_n
    &\coloneqq
    \left\{
        g_S
        =
        \int_{A_n^*(Z)} \frac{1}{2}
        \log \left(
            \frac{S(x \mid Z) + S_n^*(x \mid Z)}{2S_n^*(x \mid Z)}
        \right)
        dX(x \mid Z)
        : S \in \Fcal_n
    \right\}, \\
    \Gcal_n^\prime
    &\coloneqq
    \left\{
        g_S^\prime
        =
        \int_{A} \frac{1}{2}
        \left(
            S_n^*(x \mid Z)
            -
            \frac{S(x \mid Z) + S_n^*(x \mid Z)}{2}
        \right)
        d\nu(x)
        : S \in \Fcal_n
    \right\}.
\end{align*}
The oracle inequality is derived by bounding the empirical process on the right-hand side of Lemma~\ref{lem: basic inequality} using a Bernstein-type concentration inequality (Theorem~\ref{thm: uniform bernstein inequality}).
When applying Theorem~\ref{thm: uniform bernstein inequality}, a key difficulty, compared with density estimation in \cite{vandeGeer00}, is that the function class $\Gcal_n$ indexing the empirical process is not uniformly bounded from below.
In density estimation, the corresponding term is like
\[
\frac{1}{2}\log\left(\frac{S(x \mid z)+S_n^*(x \mid z)}{2S_n^*(x \mid z)}\right),
\]
which is uniformly bounded from below.
In that setting, this fact, together with Lemma~\ref{lem: bound for bernstein semi-norm}, is used to show that the Bernstein condition \eqref{eq: bernstein condition} holds.
In contrast, in the estimation of the intensity function of a Poisson process, the lower bound of $\Gcal_n$ involves a random variable following a Poisson distribution.
Since such a quantity is generally not bounded from below, Lemma~\ref{lem: bound for bernstein semi-norm} cannot be applied directly.
To overcome this difficulty, we exploit the fact that, given the number of points of the Poisson process $X$ in $A_n^*$, namely $X(A_n^*) = N$, the individual points are independently distributed according to the distribution with density given by the normalized intensity function; see Proposition 3.8 of \cite{LastPenrose17}.
Using this property, we show that the Bernstein condition \eqref{eq: bernstein condition} is satisfied.

Before presenting the oracle inequality, we introduce several definitions and pieces of notation.
For a given hypothesis space $\mathcal{F}_n$ and $\delta > 0$, define
\begin{equation*}
    \bar{\Fcal}_n^{1/2}(\delta)
    \coloneqq
    \left\{
        \sqrt{\frac{S + S_n^*}{2}}
        :
        \widetilde{H}^2\left(\frac{S + S_n^*}{2}, S_n^*\right) \leq \delta^2,\;
        S \in \Fcal_n
    \right\}.
\end{equation*}

Next, for \(p \geq 1\) and a measure \(Q\), we define the bracketing number of a function class $\mathcal{F}$ with respect to the $L^p(Q)$ metric as follows.

\begin{Def}[Bracketing Number for the $L^p(Q)$ Metric]
    Let $N_{p, B}(\delta, \mathcal{F}, Q)$ denote the smallest integer $N$ such that there exist pairs of functions $\{(f_j^L, f_j^U)\}_{j=1}^N$ satisfying $\|f_j^U - f_j^L\|_{L^p(Q)} \leq \delta$ for all $j \in \{1, \dots, N\}$, and such that for every $f \in \mathcal{F}$, there exists an index $j = j(f) \in \{1, \dots, N\}$ for which
    \begin{equation*}
        f_j^L \leq f \leq f_j^U.
    \end{equation*}
    The quantity $N_{p, B}(\delta, \mathcal{F}, Q)$ is called the $\delta$-bracketing number of $\mathcal{F}$.
\end{Def}

Using the two definitions above, we define the bracketing entropy integral associated with the hypothesis space \(\mathcal{F}_n\) and \(S_n^*\) in Assumption~\ref{asm: weak absolute continuity} by
\begin{equation}
    J_B(\delta, \bar{\Fcal}_n^{1/2}(\delta))
    =
    \int_{\delta^2/(2^{16} C_\rho)}^\delta
    \sqrt{
        \log N_{2, B}(u, \bar{\Fcal}_n^{1/2}(\delta), \nu \otimes P_Z)
    }\; du
    \lor \delta, \label{eq: bracketing entropy integral}
\end{equation}
where an integral with lower limit greater than its upper limit is set to zero, and the definition of $C_\rho$ can be found in the statement of Lemma~\ref{lem: hellinger upper bound for bernstein semi-norm}.
The quantity \(J_B(\delta, \bar{\mathcal{F}}_n^{1/2}(\delta))\) measures the local complexity of the hypothesis space in a neighborhood of \(S_n^*\).

We are now ready to state the oracle inequality.

\begin{theorem}[Oracle Inequality] \label{thm: oracle inequality}
    Consider the problem of estimating the conditional intensity function $S_0(\cdot \mid Z)$ of a Poisson process $X(\cdot \mid Z)$.
    Let $\hat{S}_n$ be the NPMLE defined in \eqref{eq: NPMLE}, and suppose that Assumption \ref{asm: weak absolute continuity} holds.
    Assume further that $\sup_{f \in \Fcal_n} \lVert f \rVert_\infty \leq F$ for some $F > 0$.
    For $\delta > 0$, let $\Psi(\delta)$ satisfy $\Psi(\delta) \geq J_B(\delta, \bar{\Fcal}_n^{1/2}(\delta))$ and assume that $\Psi(\delta)/\delta^2$ is a nonincreasing function of $\delta$.
    Then there exists a constant $c \ge 1$ depending only on $c_0, F, \nu(A)$ such that, for any $\delta_n$ satisfying
    \begin{equation}
        \sqrt{n}\,\delta_n^2 \geq c \Psi(\delta_n),\quad n \in \Nbb, \label{eq: critical inequality}
    \end{equation}
    we have for $\delta \ge \delta_n$,
    \begin{equation}
        \Pbb \left(
            \widetilde{H}^2(\hat{S}_n, S_0)
            >
            1154(c_0+1)^2 \bigl( \delta^2 + \widetilde{H}^2(S_n^*, S_0) \bigr)
        \right)
        \leq
        c \exp\left(-\frac{n\delta^2}{c^2}\right). \label{eq: oracle inequality}
    \end{equation}
\end{theorem}

Furthermore, by integrating both sides of \eqref{eq: oracle inequality}, we obtain the following corollary concerning the expectation with respect to the observed data.

\begin{corollary} \label{cor: oracle inequality}
    Under the same conditions as in Theorem~\ref{thm: oracle inequality}, for every $\delta \ge \delta_n$, the following inequality holds:
    \begin{equation*}
        \Ebb \left[\widetilde{H}^2(\hat{S}_n, S_0)\right]
        \leq
        1154 (c_0 + 1)^2
        \left(
            \delta^2
            +
            \widetilde{H}^2(S_n^*, S_0)
            +
            \frac{c^3}{n}
        \right).
    \end{equation*}
    Here, the expectation is taken with respect to the observed data
    $(X_1,z_1),\dots,(X_n,z_n)$.
\end{corollary}

\begin{remark}
The oracle inequalities in Theorem~\ref{thm: oracle inequality} 
and Corollary~\ref{cor: oracle inequality} hold for any 
reference intensity $S_n^* \in \mathcal F_n$ satisfying 
Assumption~\ref{asm: weak absolute continuity}. 
This condition should be viewed as a uniform lower-domination requirement: 
$S_n^*$ must be positive wherever the true conditional intensity is positive 
and must dominate $S_0$ up to a constant factor. 
For instance, if the model class contains a conditional intensity function that is uniformly bounded away from zero, then the assumption is satisfied because $S_0$ is bounded.
Nevertheless, the choice of $S_n^*$ is not merely technical. 
The resulting oracle bound also contains the approximation error 
$\widetilde H^2(S_n^*,S_0)$, and a poorly chosen reference intensity may therefore lead to a loose bound. 
To derive sharp convergence rates, it is necessary to choose $S_n^*$ so 
that the domination condition and the approximation requirement are both 
satisfied. 
In the DNN analysis below, we construct such a sequence explicitly by 
combining a neural network approximation of $S_0$ with an appropriate 
lower clipping. 
This construction ensures Assumption~\ref{asm: weak absolute continuity} and yields a 
sufficiently small approximation error, which can then be substituted into 
Corollary~\ref{cor: oracle inequality}.
\end{remark}

Theorem~\ref{thm: oracle inequality} and Corollary~\ref{cor: oracle inequality} imply that, once we can control the bracketing number representing the complexity of the hypothesis space $\Fcal_n$ (or, under suitable conditions, a covering number that provides an upper bound on the bracketing number) together with the approximation error $\widetilde{H}^2(S_n^*, S_0)$, we can derive convergence rates for the NPMLE with respect to $\widetilde{H}^2$.

\section{Maximum Likelihood Estimation of Intensity Functions Using Deep Learning}

In this section, we consider maximum likelihood estimation for conditional intensity functions using deep learning.
That is, we take the hypothesis space $\Fcal_n$ to be a class of deep neural networks (DNNs).
By applying Corollary~\ref{cor: oracle inequality}, we derive convergence rates for the DNN-based NPMLE.
We first give a mathematical definition of the deep neural networks in Section~\ref{subsec: DNN}, and then present the main results on convergence rates in Section~\ref{subsec: convergence rates}.

\subsection{Deep Neural Networks} \label{subsec: DNN}

In this section, we give a mathematical definition of DNNs.
We use the ReLU (rectified linear unit) activation function
$\sigma:\Rbb \to \Rbb,\quad x \mapsto x \lor 0.$
For a bias vector $v=(v_1,\dots,v_r)^\top \in \Rbb^r$, we define the activation function with bias by
\begin{equation*}
    \sigma_{v}
    \begin{pmatrix}
    y_1 \\
    \vdots \\
    y_r
    \end{pmatrix}
    =
    \begin{pmatrix}
    \sigma(y_1 - v_1) \\
    \vdots \\
    \sigma(y_r - v_r)
    \end{pmatrix}.
\end{equation*}

Given a final-layer activation function $\psi$, we define a neural network for modeling an intensity function with depth $L \in \Nbb$ and width $r \in \Nbb$ at each hidden layer as follows:
\begin{equation}
    f:\Rbb^d \to \Rbb,\quad
    x \mapsto f(x)
    =
    \psi W_L\sigma_{v_L}W_{L-1}\sigma_{v_{L-1}}
    \cdots W_1\sigma_{v_1}W_0x.
    \label{eq: DNN function}
\end{equation}
Here, $W_0 \in \Rbb^{r \times d}$, $W_L \in \Rbb^{1 \times r}$, and $W_i \in \Rbb^{r \times r}$, $i=1,\dots,L-1$, are weight matrices, and $v_i \in \Rbb^r$, $i=1,\dots,L$, are bias vectors.
Although the widths may vary across layers in general, we use a common width $r$ for all hidden layers for simplicity.
We assume that all network parameters are bounded in absolute value by $B>0$.
Specifically, we consider the following class of DNNs:
\begin{equation*}
    \Fcal_\psi(L,r,B,F)
    \coloneqq
    \left\{
        f \text{ of the form } \eqref{eq: DNN function}
        :
        \max_{j=0,\dots,L}
        \left(
            \lVert W_j\rVert_\infty
            \lor
            \lVert v_j\rVert_\infty
        \right)
        \leq B,
        \|f\|_\infty \le F
    \right\},
\end{equation*}
where we set $v_0=0$.
Since intensity functions are nonnegative, we use the ReLU activation function in the final layer, that is, $\psi=\sigma$.
We also define the class of DNNs with a restricted number of nonzero parameters by
\begin{equation*}
    \Fcal_\psi(L,r,B,F,s)
    \coloneqq
    \left\{
        f \in \Fcal_\psi(L,r,B,F)
        :
        \sum_{j=0}^L
        \left(
            \lVert W_j\rVert_0
            +
            \lVert v_j\rVert_0
        \right)
        \leq s
    \right\}.
\end{equation*}
Here, $\lVert W_j\rVert_0$ and $\lVert v_j\rVert_0$ denote the number of nonzero components of $W_j$ and $v_j$, respectively.

\subsection{Convergence Rates} \label{subsec: convergence rates}

In this section, we take the hypothesis space $\Fcal_n$ to be the class of DNNs defined in the previous subsection, and apply Corollary~\ref{cor: oracle inequality} to derive convergence rates for the NPMLE of the conditional intensity function based on DNNs.
To discuss convergence rates, we need to impose certain restrictions on the function class to which the true conditional intensity function $S_0$ belongs.
We assume several structural conditions on the true conditional intensity function.

We first introduce an assumption common to all settings.
An interesting phenomenon occurs when the true intensity function takes values close to zero.
To quantify this behavior, we introduce below a quantity that measures the size of the region where the intensity function is small.
This condition is a modification, for conditional intensity functions, of the $\alpha$-SVB condition introduced for conditional probabilities in Definition 3.1 of \cite{BosSchmidt-Hieber22}.

\begin{Def}[Small value bound for conditional intensity functions] \label{def: alpha-SVB condition}
    Let $\alpha \geq 0$, and let $\Hcal$ be a class of conditional intensity functions.
    We say that $\Hcal$ is $\alpha$-small value bounded, or $\alpha$-SVB, with respect to the distribution of the covariate $Z$, if there exists a constant $C>0$ such that, for every $S \in \Hcal$,
    \begin{equation*}
        \int_{[0, 1]^{d_z}} \int_A \ind \left(S(x \mid z) \leq t\right)\, d\nu(x)dP_Z(z)
        \leq
        Ct^\alpha,
        \quad \text{for all } t \in (0,\infty).
    \end{equation*}
\end{Def}

By Tonelli's theorem, the above condition can be rewritten as
\begin{equation*}
    \int_A \Pbb_Z\left(S(x \mid Z) \leq t\right)\, d\nu(x)
    \leq
    Ct^\alpha,
    \quad \text{for all } t \in (0,\infty).
\end{equation*}
Therefore, this condition always holds with $\alpha=0$ and $C=\nu(A)$.
Thus, the condition is automatically satisfied with \(\alpha=0\).
For \(\alpha>0\), however, it imposes a restriction on the size of regions where the intensity is close to zero.
On the other hand, if the set of points $x \in A$ satisfying $\Pbb_Z(S(x \mid Z)=0)>0$ has positive $\nu$-measure, then this condition cannot hold for any $\alpha>0$.
If every function in $\Hcal$ is bounded below by a constant $B_0>0$ ($\Pbb_Z-a.s.$), then the condition holds for any $\alpha \geq 0$ by taking $C=\nu(A)B_0^{-\alpha}$.
Thus, the $\alpha$-SVB condition can be interpreted as quantifying, through $\alpha$, how sharply the true conditional intensity function decays near zero after averaging over the covariates.
We illustrate this interpretation with a simple example.

Let $d_z = d_x = 1$, and let $\nu$ be the Lebesgue measure on $A = [0, 1]$.
Let $Z$ be any random variable supported on $[0, 1]$.
For $\alpha>0$, define
\[
    S(x\mid z)=z+x^{1/\alpha}.
\]
Then, for any $t\in(0,2)$,
    $\{x\in A:S(x\mid z)\le t\}
    =
    \{x\in[0,1]:z+x^{1/\alpha}\le t\}.$
If $z>t$, this set is empty, while if $0\le z\le t$, it is equal to
    $[0,(t-z)^\alpha].$
Hence,
\[
    \nu\{x\in A:S(x\mid z)\le t\}
    =
    \{(t-z) \lor 0\}^\alpha,
\]
Therefore,
\[
    \int_{[0, 1]^{d_z}} \int_A \ind\{S(x\mid z)\le t\}\,d\nu(x)dP_Z(z)
    =
    \mathbb{E}_Z[((t-Z) \lor 0)^\alpha]
    \le
    t^\alpha.
\]
Thus, the $\alpha$-SVB condition holds with $C=1$.

The following Theorem~\ref{thm: convergence rates with SVB} shows that a larger value of $\alpha$ improves the approximation error.
The proof is almost identical to that of Theorem 3.2 in \cite{BosSchmidt-Hieber22}, but we provide the details in Appendix~\ref{app: proof of thm svb} for completeness.

\begin{theorem} \label{thm: convergence rates with SVB}
    Suppose that the function class to which the true conditional intensity function $S_0$ belongs is $\alpha$-SVB.
    Then, for any approximating function $S$ satisfying
    \[
    \lVert S-S_0\rVert_\infty \leq \frac{C_1}{M},
    \qquad
    \inf_{z\in [0, 1]^{d_z}}\inf_{x\in A} S(x \mid z) \geq \frac{1}{M}
    \]
    for some constant $C_1\ge0$ and some $M\ge C_1+1$, we have
    \begin{equation*}
        \widetilde{H}^2(S,S_0)
        \leq
        (C+1)
        \frac{(C_1+1)^{2+(\alpha\land 1)}}{M^{1+(\alpha\land 1)}}
        \left(
            \nu(A)
            +
            \frac{\ind(\alpha<1)}{1-\alpha}
            +
            \log(M)
        \right).
    \end{equation*}
    Here, $C$ is the constant in the SVB condition.
\end{theorem}

\subsubsection{Convergence Rates for Composition Structured Functions}

In this section, we discuss convergence rates when the true conditional intensity function $S_0$ belongs to a class of functions called composition-structured functions.
When the true function has such a compositional structure, various studies show that deep learning can avoid the curse of dimensionality in regression and conditional probability estimation (e.g., \citealp{Schmidt-Hieber20, BosSchmidt-Hieber22, YaraTerada2026}).

We now define composition-structured functions.
Suppose that the true conditional intensity function $S_0$ can be represented as the composition of $q+1$ functions:
\begin{equation*}
    S_0 = g_q \circ g_{q-1} \circ \cdots \circ g_1 \circ g_0.
\end{equation*}
Here, for $i=0,\dots,q$, let
\[
    g_i : [a_i,b_i]^{d^{(i)}} \to [a_{i+1},b_{i+1}]^{d^{(i+1)}},
    \qquad
    a_i,b_i \in \Rbb,\quad d^{(i)} \in \Nbb,
\]
where $d^{(0)}=d$ and $d^{(q+1)}=1$.
We write each component of $g_i$ as
\[
    g_i = (g_{i1},\dots,g_{i,d^{(i+1)}})^\top,
\]
and denote by $t_i$ the maximum number of variables on which each component function $g_{ij}$, $j=1,\dots,d^{(i+1)}$, depends.
For example, consider
\[
    f(x_1,x_2,x_3)
    =
    g_{11}\bigl(g_{01}(x_1,x_2),g_{02}(x_2,x_3)\bigr)
    =
    \exp\{-(x_1^2+x_2^2)\}
    +
    \exp\{-(x_2^2+x_3^2)\}.
\]
This function has a compositional representation with
    $q=1,\,
    g_{01}(x_1,x_2)=\exp\{-(x_1^2+x_2^2)\},\,
    g_{02}(x_2,x_3)=\exp\{-(x_2^2+x_3^2)\}$,
and
    $g_{11}(u_1,u_2)=u_1+u_2$.
In this case, we have
    $d^{(0)}=3,\, t_0=2,\, d^{(1)}=t_1=2,\, d^{(2)}=1$.

We further assume that each component function $g_{ij}$ belongs to a Hölder class.
For a domain $D$, a radius $K>0$, and a smoothness parameter $\beta>0$, the Hölder class is defined by
\begin{equation*}
    \Ccal^\beta(D,K)
    =
    \left\{
    f:D\to\Rbb:
    \sum_{\gamma:\lVert\gamma\rVert_1<\beta}
    \lVert \partial^{\gamma}f\rVert_\infty
    +
    \sum_{\gamma:\lVert\gamma\rVert_1=\lfloor\beta\rfloor}
    \sup_{x,y\in D,\;x\neq y}
    \frac{
        \lvert \partial^{\gamma}f(x)-\partial^{\gamma}f(y)\rvert
    }{
        \lVert x-y\rVert_\infty^{\beta-\lfloor\beta\rfloor}
    }
    \leq K
    \right\}.
\end{equation*}
Here, we use multi-index notation. That is, for
    $\gamma=(\gamma_1,\dots,\gamma_r)\in\Nbb^r$,
we write
    $\partial^{\alpha}
    =
    \partial^{\alpha_1}\cdots\partial^{\alpha_r}$.

Thus, the function class to which the true intensity function belongs is defined by
\begin{align*}
    \Gcal(q,t,\beta,K)
    \coloneqq
    \big\{
    S = g_q \circ \cdots \circ g_0 :\,
    &S \geq 0,\,
    g_i=(g_{ij})_j:
    [a_i,b_i]^{d^{(i)}}
    \to
    [a_{i+1},b_{i+1}]^{d^{(i+1)}}, \\
    &g_{ij}\in \Ccal^{\beta_i}([a_i,b_i]^{t_i},K),
    \text{ for some }
    \lvert a_i\rvert,\lvert b_i\rvert\leq K
    \big\},
\end{align*}
where $t = (t_0,\cdots,t_q)$, and $\beta=(\beta_0,\dots,\beta_q)$.

As an illustrative example, consider the conditional intensity function of a covariate-driven Poisson process describing taxi pickup requests in an urban area.
The covariates may include daily weather conditions such as precipitation and temperature.
It is natural to expect the intensity to increase under weather conditions that make walking less comfortable, such as heavy rain, extreme heat, or extreme cold.
A typical spatial pattern of taxi demand may consist of several peaks centered around major facilities, such as department stores, railway stations, or taxi stands.
These peaks can be modeled using Gaussian-density components whose locations and shapes reflect the city's spatial structure.
The effect of covariates can then be incorporated by allowing the height, scale, or anisotropy of each Gaussian-type component to vary with them.
For instance, under rainy conditions, the peaks around major facilities may become higher and more concentrated, whereas under mild weather conditions, they may become lower and more diffuse.
A conditional intensity function of this type can be represented schematically as
\begin{equation*}
    S_0(x \mid z)
    =
    \sum_{l=1}^{L}
    a_l(z)
    \exp\left(
        -
        \left\lVert
            B_l(z)(x-\mu_l)
        \right\rVert_2^2
    \right),
\end{equation*}
where $x$ denotes the spatial location, $z$ denotes the covariates, $\mu_l$ represents the location of the $l$-th major facility, and $a_l(z)$ and $B_l(z)$ control the covariate-dependent height and shape of the corresponding peak.
This is a Gaussian-mixture-type conditional intensity function whose mixture components are modulated by covariates.
Such a model naturally has a compositional structure.
Indeed, the inner functions compute low-dimensional quantities such as covariate-dependent weights, transformed distances from major facilities, and Gaussian-type local peaks, while the outer function aggregates these components by summation.
Since the summation map is infinitely differentiable, the effective smoothness of this outer component is high.
Consequently, this representation can lead to faster convergence rates than those suggested by the ambient dimension alone, reflecting the advantage of deep learning for composition-structured functions.

For fixed $P_Z$, $\alpha\ge0$, and $C>0$, let
$\Gcal_\alpha(q,t,\beta,K,C)$ denote the subclass of
$\Gcal(q,t,\beta,K)$ consisting of all functions $S$ satisfying
\[
    \int_{[0,1]^{d_z}}\int_A
    \ind\{S(x\mid z)\le u\}\,d\nu(x)\,dP_Z(z)
    \le Cu^\alpha
    \qquad\text{for every }u>0.
\]
We assume that $S_0$ belongs to this subclass and derive convergence
rates with respect to $\widetilde H^2$ for the DNN-based NPMLE.

\begin{theorem}[Convergence Rates for Composition Structured Functions Using Deep Learning] 
\label{thm: convergence rates for composition structured functions}
    Consider the problem of estimating the conditional intensity function $S_0$ of the Poisson process $X \mid Z$, where the covariate satisfies $Z \sim P_Z$.
    Suppose that
        $S_0 \in \Gcal_\alpha(q,t,\beta,K,C)$.
    Let $\hat{S}_n$ be the NPMLE defined in \eqref{eq: NPMLE}, and let the hypothesis space $\Fcal_n$ be the class of DNNs $\Fcal_{\sigma}(L,r,B,F)$ with $F \ge K$ satisfying
    \begin{equation*}
        \mathrm{(i)}\; L \asymp \log(n),\qquad
        \mathrm{(ii)}\; r \asymp \sqrt{n\phi_n},\qquad
        \mathrm{(iii)}\; B \asymp \max_{i=0,\dots,q}
        n^{\frac{2(\beta_i+1)}{(1+(\alpha\land 1))\beta_i^*+t_i}}.
    \end{equation*}
    Then there exists a constant $C^\prime>0$, depending only on
    $A,\nu,q,t,\beta,K,F,\alpha,C, c_0$, and $d^{(0)}, \dots, d^{(q)}$, such that
    \begin{equation*}
        \Ebb\left[\widetilde{H}^2(\hat{S}_n,S_0)\right]
        \leq
        C^\prime \phi_n \log(n)^3.
    \end{equation*}
    Here,
    \[
        \phi_n
        \coloneqq
        \max_{i=0,\dots,q}
        n^{
        -
        \frac{(1+(\alpha\land 1))\beta_i^*}
        {(1+(\alpha\land 1))\beta_i^*+t_i}
        },
        \qquad
        \beta_i^*
        \coloneqq
        \beta_i
        \prod_{l=i+1}^q
        (\beta_l\land 1).
    \]
\end{theorem}

\begin{remark}
    In this theorem (also in Theorem~\ref{thm: convergence rates for manifold}), we assume that the NPMLE defined in \eqref{eq: NPMLE} exists and is computable.
    More precisely, this estimator can be interpreted as a regularized estimator because the magnitudes of the parameters are bounded by $B$.
    Since we can easily extend this theorem to the case where we cannot compute the exact NPMLE (see \cite{Schmidt-Hieber20}, for example), we state the ideal case in this theorem for simplicity.
\end{remark}

\begin{remark}
    The network budgets in Theorems~\ref{thm: convergence rates for composition structured functions} and \ref{thm: convergence rates for manifold} are calibrated using unknown smoothness and structural parameters of the true intensity, as well as the SVB exponent. The resulting bounds characterize the performance of these calibrated network classes; they do not, by themselves, establish adaptation to the unknown parameters. In the manifold setting, the dependence on intrinsic dimension enters the choice of network budgets, while the networks receive ambient coordinates and do not require intrinsic coordinates or chart maps as inputs.
    A similar issue is discussed in \cite{KohlerLanger21}, which suggests selecting network depth and width via sample splitting. Such a selection procedure is a natural extension of our approach. Establishing that a data-driven choice retains the stated Hellinger rates, however, requires an additional model-selection argument and is not asserted here.
\end{remark}

The convergence rate $\phi_n$ in Theorem~\ref{thm: convergence rates for composition structured functions} does not depend on the input dimension $d$ directly.
Instead, it depends on the intrinsic dimensions $t_i$ of the component functions $g_i$, $i=0,\dots,q$, appearing in the compositional representation.
Therefore, if
\[
    i^* \in \argmax_{i=0,\dots,q}
    \left\{
        -\frac{\beta_i^*}{\beta_i^*+t_i}
    \right\}
\]
satisfies $t_{i^*} \ll d$, then deep learning can mitigate the curse of dimensionality.

\subsubsection{Minimax Optimality for Composition Structured Functions}

In this section, we show that the convergence rate for deep learning obtained in the previous section is minimax optimal up to logarithmic factors when $\alpha = 0$.
The following theorem gives this result.

\begin{theorem}[Minimax Optimality] \label{thm: minimax optimality}
    Let $A = [0, 1]^{d_x}$ and let $\nu$ be $d_x$-dimensional Lebesgue measure.
    Consider the unrestricted small-value class ($\alpha=0$) of conditional intensity functions.
    Suppose that the covariate $Z$ has a density with respect to the Lebesgue measure, and that this density is bounded above and below by positive constants.
    Then, for any nonnegative integer $q$, any smoothness vector $\beta$, any dimensions $d^{(0)}, \dots, d^{(q)}$ and $t$ satisfying
    \[
        t_i \leq \min\{d^{(0)},\dots,d^{(i-1)}\},
        \qquad i=1,\dots,q,
    \]
    and any sufficiently large $K>0$, there exists a positive constant $c$ such that, for all sufficient large $n$,
    \begin{equation*}
        \inf_{\hat{S}_n}
        \sup_{S_0 \in \Gcal(q,t,\beta,K)}
        \Ebb\left[
            \widetilde{H}^2(\hat{S}_n,S_0)
        \right]
        \geq
        c\max_{i=0, \dots, q} n^{-\frac{\beta_i^*}{\beta_i^* + t_i}}.
    \end{equation*}
    Here, the infimum is taken over all estimators, and the expectation is taken with respect to the distribution of the data determined by $S_0$ and $P_Z$.
\end{theorem}

The lower bound differs from the classical \(n^{-2\beta/(2\beta+d)}\) rate for densities or intensities bounded away from zero. This is because the present class allows intensities to approach or equal zero.
This leads to the slower rate.

\subsubsection{Convergence Rates for Poisson Processes on Riemannian Manifolds}

In this section, we consider the case where the true Poisson process is supported on a Riemannian manifold.
That is, the region in which events occur is a low-dimensional Riemannian manifold embedded in the ambient space where the events are observed.
We show that, in this setting, the convergence rate of the NPMLE of the intensity function based on deep learning does not depend on the dimension $d_x$ of the original event-observation space, but on the intrinsic dimension $d_\Mcal$ of the low-dimensional Riemannian manifold embedded in $\Rbb^{d_x}$.

As a motivating example, point processes supported on Riemannian manifolds arise in cellular biology and microbiology, where super-resolution microscopy records the spatial arrangement of proteins and other molecules on cellular membranes.
Since such membranes are curved surfaces embedded in an ambient Euclidean space, the support of the point process is more naturally modeled as a low-dimensional Riemannian manifold than as the full ambient Euclidean space \citep{WardEtAl23}.
Related theoretical work on point processes in manifolds also emphasizes that, when data are embedded in a high-dimensional Euclidean space but lie on a lower-dimensional manifold, intrinsic dimension and manifold learning become central objects of analysis \citep{PenroseEtAl2013}.
Spherical point processes provide another concrete example, where the support is the sphere, a canonical Riemannian manifold, and the goal is to estimate the intensity function of a non-homogeneous Poisson process on that manifold \citep{NgZammit-Mangion2022}.

Before stating the convergence rates, we introduce several definitions.
Let $\Mcal\subset\Rbb^{d_x}$ be a fixed compact $C^\infty$ embedded
manifold without boundary, of dimension $d_\Mcal$, and let $\nu$ be its
Riemannian volume measure, with $0<\nu(\Mcal)<\infty$.
The likelihood uses this known reference measure.

We fix a finite atlas $\{(V_i,P_i)\}_{i=1}^J$ of $\Mcal$ in which each
$P_i$ is the restriction of an affine ambient map and
$P_i(V_i)=\Omega_i$ is an open ball in $\Rbb^{d_\Mcal}$.  Such charts
are obtained by restricting projections onto tangent spaces.  Fix also a
smooth partition of unity $\{\rho_i\}_{i=1}^J$ with
$\operatorname{supp}(\rho_i)\subset V_i$.  All geometric choices remain
fixed as the intensity and the sample size vary.

\begin{Def}[H\"older class for a fixed atlas]\label{def:manifold-holder}
For $i=1,\ldots,J$,
write
\[
 D_i=\Omega_i\times[0,1]^{d_z},\qquad
 S_i(u,z)=S(P_i^{-1}(u),z).
\]
For a radius $K > 0$, output bound $F > 0$, and a smoothness parameter $\beta > 0$, the Hölder class on $\mathcal{M} \times [0, 1]^{d_z}$ is defined by
\[
 \Hcal^\beta_{\Mcal}(K, F)
 \coloneqq
 \left\{S:\mathcal{M} \times [0, 1]^{d_z}\to\Rbb: 0 \le S \le F,\;
 S_i\in C^{\beta}(D_i, K)
 \right\}.
\]
\end{Def}

For fixed $P_Z$, $\alpha\ge0$, and $C>0$, let
$\Hcal^\beta_{\Mcal,\alpha}(K,F,C)$ denote the subclass of
$\Hcal^\beta_{\Mcal}(K,F)$ consisting of all functions $S$ satisfying
\[
    \int_{[0,1]^{d_z}}\int_{\Mcal}
    \ind\{S(x\mid z)\le u\}\,d\nu(x)\,dP_Z(z)
    \le Cu^\alpha
    \qquad\text{for every }u>0.
\]

The next theorem derives the convergence rates for manifold-supported intensities.

\begin{theorem}[Convergence Rates for Manifold-Supported Intensities Using Deep Learning]
\label{thm: convergence rates for manifold}
Consider the problem of estimating the conditional intensity function
$S_0$ of the Poisson process $X\mid Z$, where $Z\sim P_Z$.  Fix the
manifold, atlas, partition, and reference measure specified above.
Suppose that $\beta>0$, $\alpha\ge0$, $K,F,C>0$, and $S_0\in\Hcal^\beta_{\Mcal,\alpha}(K,F,C).$
Let $\hat S_n$ be the NPMLE defined in \eqref{eq: NPMLE}, and let the
hypothesis space $\Fcal_n$ be the class of DNNs
$\Fcal_\sigma(L,r,B,F)$ satisfying
\begin{equation*}
 \mathrm{(i)}\; L\asymp\log(n),\qquad
 \mathrm{(ii)}\; r\asymp\sqrt{n\phi_n'},\qquad
 \mathrm{(iii)}\; B\asymp
 n^{\frac{2(\beta+1)}{(1+(\alpha\land1))\beta+d_\Mcal + d_z}}.
\end{equation*}
Then there exists a constant $C^\prime>0$, depending only on the fixed
geometry, $\beta,\alpha,K,F,C,P_Z$, and the chosen architecture
constants, such that, for all sufficiently large $n$,
\begin{equation*}
 \Ebb[\widetilde H^2(\hat S_n,S_0)]
 \le C^\prime\phi_n'\log(n)^3.
\end{equation*}
Here,
\[
 \phi_n'
 \coloneqq n^{-\frac{(1+(\alpha\land1))\beta}
                         {(1+(\alpha\land1))\beta+d_{\Mcal} + d_z}}.
\]
\end{theorem}

The exponent depends on the intrinsic input dimension $d_\Mcal + d_z$ rather than
$d_x + d_z$.  The networks nevertheless take ambient coordinates as inputs; the
charts are used to prove approximation within that network class.
The likelihood uses the exact volume measure.

\section{Simulation Studies}
\label{sec:numerical-experiments}

We conducted simulation studies to examine how estimation accuracy depends
on sample size, covariate dimension, and the structure of the intensity.
Three designs were considered: a compositional intensity, an
intensity with a deep valley whose location varies with the covariates,
and intensities supported on a circle or sphere. We compared DNN-based estimators with kernel estimators, including a pooled estimator to assess the effect of ignoring covariate dependence.

\subsection{Simulation design}
\label{subsec:numerical-design}

We generated independent pairs
$Z_i\sim\operatorname{Unif}([0,1]^{d_z})$ and
$X_i\mid Z_i\sim\operatorname{PP}\{S_0(\cdot\mid Z_i)\,d\nu\}$,
$i=1,\ldots,n$, with $d_z\in\{0,1,5\}$ and
$n\in\{100,316,1000,3162,10000\}$. Each configuration was repeated
20 times. Here $n$ is the number of point patterns. We scaled the
intensities to give approximately 30 events per pattern on average.
The Euclidean supports were $[0,1]$ and $[0,1]^2$; on the unit circle
and sphere, $\nu$ was arc length and surface area, respectively.

The compositional design combines low-dimensional functions of event
coordinates and covariates. The second design has a covariate-dependent
valley with a positive floor of $10^{-4}$ before rescaling, allowing us to
examine estimation in regions of low intensity. In both designs, at most
the first two covariates enter the intensity, so the five-covariate cases
also assess sensitivity to irrelevant coordinates. The manifold designs
introduce periodic variation with angle and covariates. The spherical
intensity includes a longitude-dependent term that is discontinuous at
the poles; this design extends the numerical comparison beyond the smooth
class in Theorem~\ref{thm: convergence rates for manifold}.
Section~\ref{subsec:simulation-generators} gives the generating functions
and simulation procedure.

Both DNN-based estimators have ReLU hidden activations. Their output functions
are Softplus and ReLU, each with $10^{-8}$ added to ensure positivity for an initial stage of the optimization;
we refer to them as DNN--Softplus and DNN--ReLU. We trained the networks
by minimizing the negative log-likelihood using AdamW and
stochastic quadrature. A random 80\% of the patterns was used for
parameter updates, and the remaining 20\% for early stopping. Each update
used a batch of 64 patterns and 128 integration points per pattern.
Training was limited to 200 epochs. Network depth and width increased with $n$
according to the rule in Section~\ref{subsec:simulation-implementation}.

For Euclidean supports, we used a Gaussian conditional kernel estimator
and a pooled kernel estimator. The pooled estimator ignores covariates;
when $d_z=0$, the two estimators coincide. Bandwidths were selected by
five-fold cross-validation, using at most 512 patterns and 128
events per pattern for this calculation, and the selected estimator was
then fitted to all $n$ patterns. For manifold supports, the kernel
estimator used the known geodesic distance and volume measure. The main
DNN comparison used a reconstructed integration mesh, obtained from
independent samples on the support together with event locations. Thus,
the mesh construction used support information in addition to the observed
point patterns. We also fitted DNNs using the known manifold measure;
Section~\ref{subsec:simulation-additional-results} compares these versions.

Accuracy was measured by the covariate-averaged squared Hellinger loss
\begin{equation*}
 \widetilde H^2(\widehat S,S_0)
 =\frac12\mathbb E_Z\left[\int
 \left(\sqrt{\widehat S(x\mid Z)}-\sqrt{S_0(x\mid Z)}\right)^2\,d\nu(x)\right].
\end{equation*}
We evaluated this loss using numerical integration and covariate draws
independent of the fitting sample. Within each repetition, all methods
used the same simulated data and evaluation points.
Table~\ref{tab:main-results} summarizes the results at $n=10000$.
Monte Carlo standard errors (MCSEs) are the sample standard deviations
divided by $\sqrt{20}$. Learning curves and results for all sample sizes
are reported below and in Section~\ref{sec:simulation-supplement}.

\begin{table}[t]
\centering
\caption{Mean squared Hellinger loss at $n=10000$, with Monte Carlo standard errors in parentheses, based on 20 repetitions. Here $d$ is the event-support dimension. ``Kernel'' denotes the conditional kernel for Euclidean supports and the geodesic kernel for manifold supports; the manifold DNNs use reconstructed integration meshes. A dash indicates that the pooled comparator was not used.}
\label{tab:main-results}
\small
\setlength{\tabcolsep}{4pt}
\begin{tabular}{lrrllll}
\toprule
Design/support & $d$ & $d_z$ & Softplus & ReLU & Kernel & Pooled \\
\midrule
Composition & 1 & 1 & 0.00926 (0.00073) & 0.0104 (0.00096) & 0.0329 (0.0048) & 0.474 (0.007) \\
Composition & 1 & 5 & 0.0194 (0.0011) & 0.0203 (0.0011) & 0.144 (0.01) & 0.544 (0.011) \\
Composition & 2 & 1 & 0.0205 (0.00071) & 0.0213 (0.00099) & 0.0441 (0.0055) & 0.398 (0.0052) \\
Composition & 2 & 5 & 0.0346 (0.0013) & 0.0329 (0.0012) & 0.182 (0.0062) & 0.466 (0.0067) \\
Near zero & 1 & 1 & 0.0092 (0.0012) & 0.0229 (0.0023) & 0.0442 (0.00036) & 2.61 (0.044) \\
Near zero & 1 & 5 & 0.0179 (0.001) & 0.0359 (0.0032) & 0.597 (0.0051) & 2.67 (0.031) \\
Near zero & 2 & 1 & 0.00788 (0.00056) & 0.0108 (0.00085) & 0.0455 (0.0035) & 0.516 (0.012) \\
Near zero & 2 & 5 & 0.0296 (0.0011) & 0.0408 (0.0022) & 0.408 (0.0082) & 0.951 (0.014) \\
Circle & 1 & 1 & 0.0057 (0.0005) & 0.00643 (0.00051) & 0.0052 (0.00028) & -- \\
Circle & 1 & 5 & 0.0105 (0.00066) & 0.0108 (0.00057) & 0.0244 (0.00093) & -- \\
Sphere & 2 & 1 & 0.0102 (0.00058) & 0.00997 (0.00057) & 0.0205 (0.00084) & -- \\
Sphere & 2 & 5 & 0.0184 (0.00077) & 0.0182 (0.00075) & 0.0392 (0.0014) & -- \\
\bottomrule
\end{tabular}
\end{table}

\subsection{Euclidean supports}
\label{subsec:numerical-compositional-results}

\begin{figure}[t]
 \centering
 \includegraphics[width=0.88\textwidth]{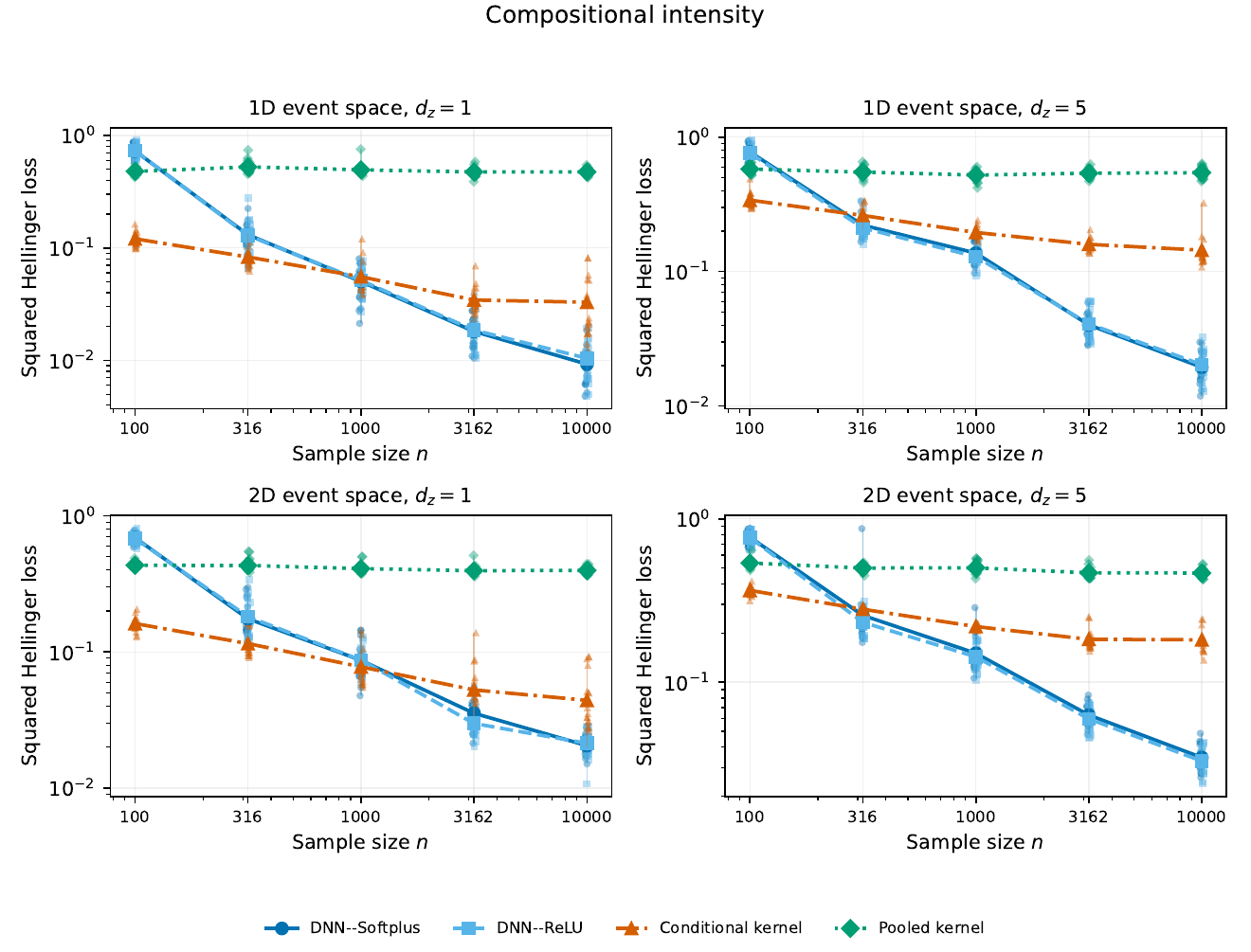}
 \caption{Squared Hellinger loss for the compositional intensity.
 Rows correspond to event dimension and columns to covariate dimension.
 Points represent individual repetitions, connected markers are means,
 and vertical segments span the minimum and maximum over 20 repetitions.}
 \label{fig:numerical-compositional}
\end{figure}

Figure~\ref{fig:numerical-compositional} shows the results for the
compositional intensity. The conditional kernel performed well
at small sample sizes. With one event coordinate, one covariate, and
$n=100$, its mean loss was $0.120$, compared with $0.728$ for
DNN--Softplus. At $n=1000$, the corresponding losses were $0.0554$ and
$0.0503$. With two event coordinates and one covariate, the kernel still
had a smaller mean loss at $n=1000$: $0.0777$, compared with $0.0869$
for DNN--Softplus. The unconditional results showed a similar
small-sample advantage for kernel smoothing: for example, with one event
coordinate and $n=100$, the pooled kernel and DNN--Softplus had mean
losses of $0.0682$ and $0.214$, respectively.

As $n$ increased, the DNN estimators became more competitive. At
$n=10000$, a DNN attained the smallest mean loss in each of the four
covariate-dependent compositional configurations. The difference was
largest with five covariates: the smaller of the two DNN mean losses was
$0.0194$ in one event dimension and $0.0329$ in two event dimensions,
compared with $0.144$ and $0.182$ for the conditional kernel. These
comparisons indicate that neural estimation can be advantageous when
several covariates are present, even though only a subset affects the
intensity. The relative performance of the two output functions varied
across configurations. The pooled kernel retained appreciable error at
large $n$, reflecting the covariate variation that it cannot represent.

For the intensity with a moving valley, DNN--Softplus had a smaller mean
loss than DNN--ReLU in 25 of the 30 configurations
(Figure~\ref{fig:numerical-near-zero}). It also attained the smallest mean
loss in all four covariate-dependent configurations at $n=10000$.
With two event coordinates and five covariates, its mean loss was
$0.0296$ (MCSE $0.0011$), compared with $0.0408$ ($0.0022$) for
DNN--ReLU and $0.408$ ($0.0082$) for the conditional kernel. Thus,
Softplus performed particularly well in this design with regions of very
low intensity. Because the positive floor was held fixed, these results
concern finite-sample estimation near zero rather than comparisons among
the small-value exponents appearing in the theory.

\subsection{Manifold supports}
\label{subsec:numerical-manifold-results}

The geodesic kernel remained competitive on manifold supports
(Figure~\ref{fig:numerical-manifold-statistical}). On the circle with
one covariate and $n=10000$, its mean loss was $0.00520$, slightly smaller
than $0.00570$ for the mesh-based DNN--Softplus. With five covariates,
the ordering reversed: the losses were $0.0244$ and $0.0105$,
respectively. On the sphere with five covariates, the kernel had a smaller
mean loss at $n=1000$ ($0.0480$ versus $0.0646$ for DNN--Softplus),
whereas at $n=10000$ the DNN losses were approximately $0.018$, compared
with $0.0392$ for the kernel. As in the Euclidean experiments, the DNN
advantage was most apparent at the largest sample size with several
covariates.

Using the known manifold measure did not uniformly improve the DNN
estimates. For example, with five covariates and $n=10000$, the
mesh-based and known-measure DNN--Softplus losses were $0.0105$ and
$0.0102$ on the circle, and $0.0184$ and $0.0185$ on the sphere.
These versions also use different dimensions in the network-width rule,
so their comparison reflects both architecture and integration measure.
The full results in Section~\ref{sec:simulation-supplement} document the
variation across sample sizes and include the no-covariate cases.

\providecommand{\AllenPublicationRoot}{allen}
\section{Stimulus-Dependent Neural Spike Intensity}
\label{sec:real-data}

Visual stimuli can affect both the number and timing of neural spikes.
We analyzed repeated Gabor presentations from the Allen Brain Observatory
Visual Coding Neuropixels project \citep{siegle2021survey} to estimate
how the temporal firing-rate profile varies with stimulus position and
orientation. The main question was whether a neural intensity estimator
could improve prediction beyond conditional kernel smoothing. A pooled
temporal estimator provided a second comparison to quantify the benefit
of accounting for the stimulus.

\subsection{Data and analysis design}
\label{subsec:allen-design}

We used session 715093703, selected as the smallest session identifier
among sessions with an available recording file and at least 20 eligible
primary visual cortex (VISp) units. After excluding 620 presentations
that overlapped invalid recording intervals and one with a different
duration, 3,024 of the 3,645 Gabor presentations remained. Each
presentation lasted $T=0.25$ seconds. There were 243 stimulus conditions,
formed by 81 screen positions and three axial orientations.

The covariate vector contained horizontal and vertical stimulus position,
scaled using development-data extrema, and $\sin(2\theta)$ and
$\cos(2\theta)$ for orientation $\theta$. The doubled-angle encoding
accounts for the equivalence of orientations separated by $\pi$.
We split presentations within stimulus condition, assigning 2,417 to
development and 607 to testing. Both subsets contained all 243 conditions.
This design evaluates prediction for new presentations of conditions
represented in the development sample.

Of the 62 VISp units satisfying the recording-quality criteria, we selected
the 12 with the highest development-sample spike counts and estimated a
separate intensity for each unit. Quality thresholds and preprocessing
are detailed in Appendix~\ref{sec:allen-supplement}. Test results for one
of these units, 950931315, had been examined in a preliminary analysis.
We therefore report results both for the 12-unit panel and for the
11 units excluding this previously examined unit. All models reported
here were refitted using the development-sample preprocessing described
below.

\subsection{Estimation and predictive assessment}
\label{subsec:allen-methods}

We compared DNN--Softplus, DNN--ReLU, a conditional Gaussian kernel
estimator, and a pooled temporal kernel estimator. Both networks had
three hidden layers of width 64 with ReLU activations and a positive
output offset of $10^{-8}$. We selected the number of training epochs
using 1,935 development presentations for fitting and 482 for validation,
then refitted each network from initialization on all 2,417 development
presentations for the selected number of epochs. Kernel bandwidths were
chosen by five-fold presentation-level Poisson cross-validation on the
development sample. The conditional kernel used the same stimulus
covariates as the networks, whereas the pooled kernel used only spike
times. Additional fitting details appear in
Appendix~\ref{sec:allen-supplement}.

For unit $v$, let $t_{vij}$ be spike $j$ in presentation $i$, with
$m_{vi}$ spikes and covariate vector $z_i$. We fitted intensities on
normalized time $u=t/T$ and converted them to spikes per second through
$\widehat S_{v,t}(t\mid z)=\widehat S_{v,u}(t/T\mid z)/T$.
Predictive performance was assessed by the Poisson negative log-likelihood
(NLL) on the test presentations,
\begin{equation}
 L_{vi}(\widehat S)
 =\int_0^T \widehat S_{v,t}(t\mid z_i)\,dt
  -\sum_{j=1}^{m_{vi}}\log\widehat S_{v,t}(t_{vij}\mid z_i),
 \label{eq:allen-nll}
\end{equation}
with parameter-independent terms omitted. Integrals were evaluated on a
1,024-point trapezoidal grid. Scores are in nats per presentation, with
time measured in seconds; paired score differences are invariant to the
choice of time unit.

For a comparison of methods $a$ and $b$, we averaged
$L_{vi}(a)-L_{vi}(b)$ over test presentations within each unit, then
averaged equally over units. A negative difference favors method $a$.
We also report the equal-unit mean information gain in bits per spike,
obtained by dividing each unit's total log-likelihood improvement by
its test spike count and $\log 2$. To quantify uncertainty across
stimulus conditions, we drew 2,000 paired bootstrap samples of the 243
conditions, retaining all presentations within a sampled condition and
using the same draws for all methods and units. The resulting percentile
intervals are pointwise 95\% intervals conditional on the fitted models
and selected units; they do not describe variation across neuronal
populations or recording sessions.

\subsection{Results}
\label{subsec:allen-results}

Accounting for the stimulus substantially improved prediction. Relative
to the pooled temporal kernel, mean NLL decreased by $0.38926$ nats
per presentation for the conditional kernel, $0.41827$ for
DNN--Softplus, and $0.41006$ for DNN--ReLU. Both neural estimators
improved on the pooled estimator in all 12 units. The fitted profiles
in Figure~\ref{fig:allen-intensities} illustrate this stimulus dependence:
for unit 950931315, the covariate-dependent estimators distinguish the
transient response at the preferred condition from the weaker responses
at the other two conditions. The pooled estimator gives the same temporal
profile for every stimulus.

The DNN estimators provided a smaller additional improvement over the
conditional kernel (Table~\ref{tab:allen-primary}). For the 12-unit panel,
the mean NLL differences were $-0.02901$ for DNN--Softplus and
$-0.02080$ for DNN--ReLU, with 95\% intervals
$[-0.04042,-0.01729]$ and $[-0.03092,-0.01026]$, respectively.
The corresponding information gains were $0.01073$ and $0.00600$ bits
per spike. Excluding the previously examined unit gave similar results:
the mean differences were $-0.02610$ and $-0.01954$, with both
intervals remaining below zero.

\begin{table}[tbp]
\centering\small
\caption{Predictive performance relative to the conditional kernel.
NLL differences and information gains are averaged equally over units;
brackets give pointwise 95\% paired condition-bootstrap intervals.
``Improved'' is the number of units with a negative NLL difference.
The 11-unit analysis excludes unit 950931315, examined in the pilot.}
\label{tab:allen-primary}
\setlength{\tabcolsep}{4pt}
\begin{tabular}{clrrr}
\toprule
Units & DNN & \shortstack{NLL difference\\(nats/presentation)} & Improved
& \shortstack{Information gain\\(bits/spike)}\\
\midrule
12 & Softplus & $-0.02901$ & 10/12 & $0.01073$\\
   &         & $[-0.04042,-0.01729]$ & & $[0.00606,0.01591]$\\
12 & ReLU    & $-0.02080$ & 9/12 & $0.00600$\\
   &         & $[-0.03092,-0.01026]$ & & $[0.00175,0.01001]$\\
11 & Softplus & $-0.02610$ & 9/11 & $0.00944$\\
   &         & $[-0.03791,-0.01378]$ & & $[0.00473,0.01461]$\\
11 & ReLU    & $-0.01954$ & 8/11 & $0.00526$\\
   &         & $[-0.02991,-0.00933]$ & & $[0.00101,0.00930]$\\
\bottomrule
\end{tabular}
\end{table}

\begin{figure}[tbp]
\centering
\includegraphics[width=0.94\linewidth]{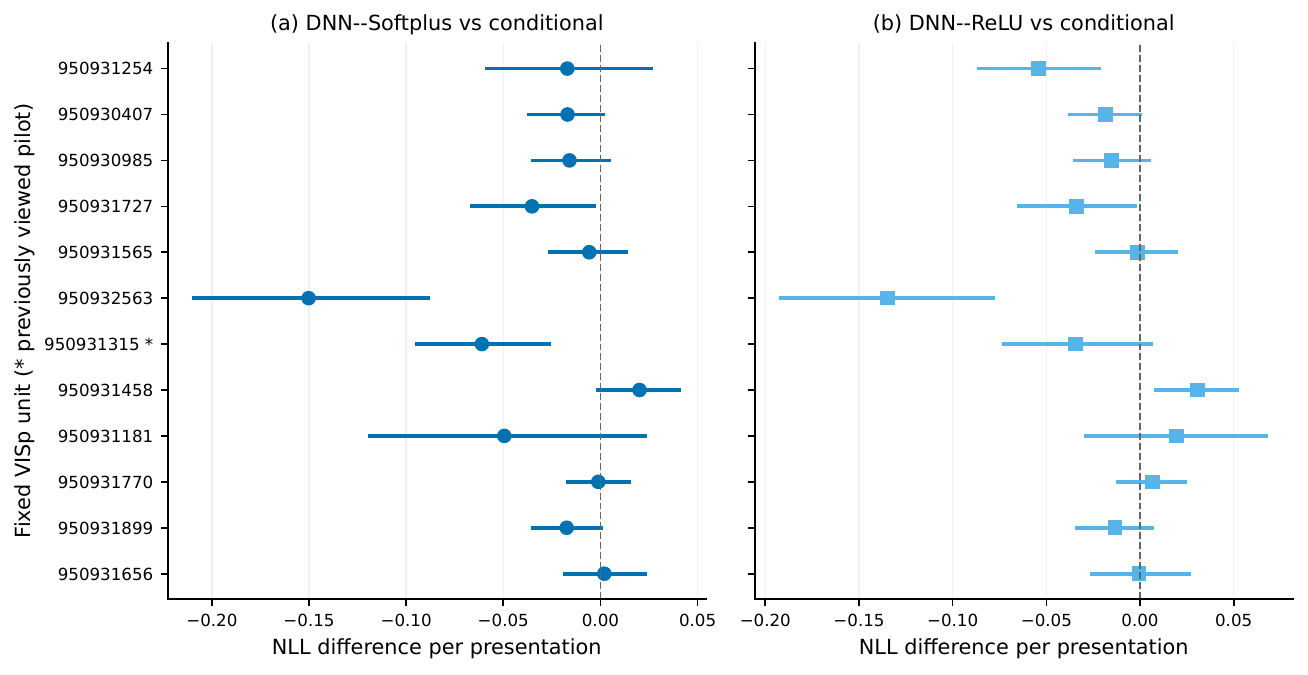}
\caption{Unit-specific mean test NLL differences from the conditional
kernel, in nats per presentation. Negative values favor the DNN.
Horizontal segments give pointwise 95\% paired condition-bootstrap
intervals. The asterisk denotes unit 950931315, examined in the pilot.}
\label{fig:allen-units}
\end{figure}

Performance varied across units (Figure~\ref{fig:allen-units}).
DNN--Softplus had a lower mean NLL than the conditional kernel for
10 of the 12 units, and DNN--ReLU for nine. Both had their largest gain
in unit 950932563, whereas the conditional kernel performed better in
unit 950931458. Many unit-level intervals included zero. Although
DNN--Softplus had the lower panel-average NLL, its difference from
DNN--ReLU was $-0.00821$, with interval $[-0.01728,0.00092]$.
Thus, the NLL comparison did not clearly distinguish the two output
functions.

Changing the integration grid had a negligible effect on the reported
score differences. Supplementary analyses examine network initialization,
width, running speed, prediction at held-out stimulus conditions, and
temporal variation. The main comparison indicates that most of the gain
over a pooled temporal estimate comes from modeling stimulus dependence,
with a further modest improvement from the DNNs. This finding concerns
repeated presentations in one session and a selected set of highly active
units. The Poisson likelihood serves as a criterion for first-order
intensity prediction; the time-rescaling diagnostics in
Appendix~\ref{sec:allen-supplement} show remaining departures from a full
Poisson description of the spike trains.

\section{Conclusion and Discussion}\label{sec:discussion}
Likelihood-based intensity estimation can be analyzed without requiring the entire candidate class to be uniformly bounded away from zero. The oracle inequality makes the relevant trade-off explicit: an approximating intensity must both dominate the truth and approximate it under covariate-averaged Hellinger loss, while the local entropy controls estimation error. The small-value condition describes how the behavior of the true intensity near zero affects this trade-off. For compositional functions, the resulting rate reflects the dimensions of the component maps. The minimax comparison applies to the unrestricted small-value class on the specified Euclidean domain; it does not establish optimality for every positive small-value exponent.

The simulations show that the relative performance of neural and kernel estimation depends on sample size and covariate dimension. Kernel estimators often perform better in small samples, whereas the DNNs are more accurate in the largest covariate-dependent Euclidean settings and in several higher-dimensional manifold settings. Softplus performs particularly well for the intensity with a moving low-intensity region, but neither output function is uniformly preferable. These findings complement the approximation theory by showing the finite-sample trade-offs under the fitting procedures used here.

In the neural-recording application, modeling stimulus dependence accounts for most of the improvement over a pooled temporal estimate. The DNNs provide an additional, more modest improvement over conditional kernel smoothing, with variation across units and no clear ordering of the two output functions by mean NLL. The analysis concerns repeated presentations within one session and a selected set of highly active units. Validation across sessions and animals would be needed to assess broader predictive performance.

Several extensions would strengthen the connection between theory and practice. First, the current theory concerns an exact constrained ReLU likelihood estimator with exact integration. A computational guarantee would need to account jointly for optimization error, quadrature error, and the output transformation used by the fitted networks. Second, learning an unknown support manifold requires control of the estimated integration measure in addition to function approximation on that support. The manifold theorem establishes approximation and risk bounds on a fixed smooth support for every $\beta>0$, with the chartwise H\"older convention of Definition~\ref{def:manifold-holder}; it does not solve support estimation or adaptation to unknown intrinsic dimension. Third, repeated-session validation and chronological holdout would address forms of generalization not measured by the current neural-recording split. Modeling spike-history effects would also permit a direct assessment of departures from the conditional Poisson working model.

\appendix
\section{Proof of lemmas in Section~\ref{sec: difficulty}}
\subsection{Proof of Lemma~\ref{lem: the excess risk}}
\begin{proof}
    We first compute $\inf_{S \in \Fcal_0} \mathbb{E}[\ell(S)]$.
    Suppose that, $\log(S(\cdot \mid z)) \in L^1(\lambda_0(\cdot \mid z)),\, P_Z\text{-}a.s.$
    Then, by Campbell's formula \cite[Proposition 2.7]{LastPenrose17},
    \begin{align*}
        \mathbb{E}[\ell(S)]
        &= \mathbb{E}\left[
            -\int_A \log(S(x \mid Z))\, dX(x \mid Z)
            + \int_A S(x \mid Z)\, d\nu(x)
        \right] \\
        &= \mathbb{E}_Z \left[
            -\mathbb{E}\left[
                \int_A \log(S(x \mid Z))\, dX(x \mid Z)
                \middle| Z
            \right]
            + \int_A S(x \mid Z)\, d\nu(x)
        \right] \\
        &= \mathbb{E}_Z \left[
            -\int_A \log(S(x \mid Z)) S_0(x \mid Z)\, d\nu(x)
            + \int_A S(x \mid Z)\, d\nu(x)
        \right] \\
        &= \mathbb{E}_Z \left[
            \int_A \left\{
                S(x \mid Z) - S_0(x \mid Z)\log(S(x \mid Z))
            \right\}\, d\nu(x)
        \right].
    \end{align*}
    Since, for $a \ge 0$ and $y \ge 0$, the function $y - a\log y$ is minimized at $y=a$, it follows that
    $\inf_{S \in \Fcal_0} \mathbb{E}[\ell(S)] = \mathbb{E}[\ell(S_0)]$.

    Indeed, for sufficiently small $\epsilon > 0$,
    \begin{align*}
        &\int_A \lvert \log(S_0(x \mid z)) \rvert\, d\lambda_0(x \mid z)
        = \int_A \lvert \log(S_0(x \mid z)) \rvert S_0(x \mid z)\, d\nu(x) \\
        &= \int_A \lvert \log(S_0(x \mid z)) \rvert S_0(x \mid z)\ind(S_0(x \mid z)\ge \epsilon)\, d\nu(x)
         + \int_A \lvert \log(S_0(x \mid z)) \rvert S_0(x \mid z)\ind(S_0(x \mid z)<\epsilon)\, d\nu(x) \\
        &\le \log\left(\frac{1}{\epsilon}\right)
            \int_A S_0(x \mid z)\ind(S_0(x \mid z)\ge \epsilon)\, d\nu(x)
         + \int_A \log\left(\frac{1}{S_0(x \mid z)}\right)
            S_0(x \mid z)\ind(S_0(x \mid z)<\epsilon)\, d\nu(x) \\
        &\le \log\left(\frac{1}{\epsilon}\right)
            \int_A S_0(x \mid z)\ind(S_0(x \mid z)\ge \epsilon)\, d\nu(x) 
         + \int_A (1-S_0(x \mid z))
            \ind(S_0(x \mid z)<\epsilon)\, d\nu(x).
    \end{align*}
    Since $S_0$ is bounded and $A$ is compact, we obtain
    $\log(S_0(\cdot \mid z)) \in L^1(\lambda_0(\cdot \mid z)),\, P_Z\text{-a.s.}$
    Here, in the first inequality, we used the boundedness of $S_0$ and the fact that $\epsilon$ is sufficiently small.
    In the last inequality, we used $\log x \le x-1$ for $x>0$.
    On the other hand, if $\log(S(\cdot \mid z)) \notin L^1(\lambda_0(\cdot \mid z))$ for some $\Zcal \subset [0, 1]^{d_z}$ with $P_Z(\Zcal) > 0$, then $\mathbb{E}[\ell(S)] = \infty$.
    To see this, first suppose that $S(\cdot \mid z)$ is a simple function, that is,
    \begin{equation*}
        S(x \mid z) = S^{(m)}(x) \coloneqq \sum_{j=1}^m S_j \ind(x \in A_j),
    \end{equation*}
    where $S_j \ge 0$ are constants, $m<\infty$, and $A_j \in \Bcal(A)$ are pairwise disjoint.
    Then
    \begin{align*}
        -\mathbb{E}\left[\int_A \log(S(x \mid z)) dX(x \mid z)\right]
        = -\sum_{j=1}^m \log(S_j)\lambda_0(A_j \mid z) 
        = -\int_A \log(S(x \mid z))\, d\lambda_0(x \mid z).
    \end{align*}
    Next, consider a general nonnegative, bounded, measurable function.
    By measurability, there exists a sequence of nonnegative simple functions $(S^{(m)})$ such that, for every $x \in A$,
    \begin{equation*}
        S^{(m)}(x \mid z) \to S(x \mid z)
        \; \text{ as } m \to \infty.
    \end{equation*}
    Hence, by the monotonicity and continuity of the logarithm, the boundedness of $S$, and the monotone convergence theorem, for any $\omega \in \Omega$,
    \begin{align*}
        -\int_A \log(S^{(m)}(x \mid z)) dX(x \mid z)
        &= \int_A \{-\log(S^{(m)}(x \mid z)) \lor 0\} dX(x \mid z) 
         - \int_A \{\log(S^{(m)}(x \mid z)) \lor 0\} dX(x \mid z) \\
        &\nearrow \int_A \{-\log(S(x \mid z)) \lor 0\} dX(x \mid z)
         - \int_A \{\log(S(x \mid z)) \lor 0\} dX(x \mid z) \\
        &= -\int_A \log(S(x \mid z)) dX(x \mid z)
        \;\text{ as }\; m \to \infty.
    \end{align*}
    Therefore, again by the monotonicity and continuity of the logarithm, the boundedness of $S$, and the monotone convergence theorem,
    \begin{align*}
        -\mathbb{E}\left[\int_A \log(S(x \mid z)) dX(x \mid z)\right]
        &= \mathbb{E}\left[
            \lim_{m\to\infty}
            \left\{
                -\int_A \log(S^{(m)}(x \mid z)) dX(x \mid z)
            \right\}
        \right] \\
        &= \lim_{m\to\infty}
        \mathbb{E}\left[
            -\int_A \log(S^{(m)}(x \mid z)) dX(x \mid z)
        \right] \\
        &= -\lim_{m\to\infty}
        \int_A \log(S^{(m)}(x \mid z))\, d\lambda_0(x \mid z) \\
        &= -\int_A \log(S(x \mid z))\, d\lambda_0(x \mid z).
    \end{align*}
    Since $\log(S(\cdot \mid z)) \notin L^1(\lambda_0(\cdot \mid z))$ and $S$ is bounded, we must have
    \[
    -\int_A \log(S(x \mid z))\, d\lambda_0(x \mid z) = \infty.
    \]
    Hence, together with the assumption that the covariate $Z$ has a positive density on its domain, this implies that if
    $\log(S(\cdot \mid z)) \notin L^1(\lambda_0(\cdot \mid z))$, then $\mathbb{E}[\ell(S)] = \infty$, and therefore
    \[
    \inf_{S \in \Fcal_0} \mathbb{E}[\ell(S)] = \mathbb{E}[\ell(S_0)].
    \]

    Next, we compute $R(\hat{S}_n)$.
    If $\lambda_0(\cdot \mid z) \ll \hat{\lambda}_n(\cdot \mid z),\, P_Z-a.s.$, then by the formula of the KL divergence \cite[Corollary 5.2]{leskelä24} and Campbell's formula, we have
    \begin{align*}
        R(\hat{S}_n)
        &= \mathbb{E}[\ell(\hat{S}_n)] - \mathbb{E}[\ell(S_0)] \\
        &= \mathbb{E}\left[
            \int_A \log\left(\frac{S_0(x \mid Z)}{\hat{S}_n(x \mid Z)}\right) dX(x \mid Z)
            + \int_A \left(\hat{S}_n(x \mid Z) - S_0(x \mid Z)\right)\, d\nu(x)
        \right] \\
        &= \mathbb{E}_Z \left[
            \mathbb{E}\left[
                \int_A \log\left(\frac{S_0(x \mid Z)}{\hat{S}_n(x \mid Z)}\right) dX(x \mid Z)
                \middle| Z
            \right]
            + \int_A \left(\hat{S}_n(x \mid Z) - S_0(x \mid Z)\right)\, d\nu(x)
        \right] \\
        &= \mathbb{E}_Z\!\left[
            \mathrm{KL}\!\left(P_{\lambda_0(\cdot \mid Z)} \parallel P_{\hat{\lambda}_n(\cdot \mid Z)}\right)
        \right].
    \end{align*}
    On the other hand, if there exists some $\Zcal \subset [0, 1]^{d_z}$ with $P_Z(\Zcal) > 0$ such that
    $\lambda_0(\cdot \mid z) \ll \hat{\lambda}_n(\cdot \mid z)$ is not satisfied for all $z \in \Zcal$, then there exists
    $B_z \in \Bcal(A)$ such that $\lambda_0(B_z \mid z) > 0$ and $\hat{\lambda}_n(B_z \mid z) = 0$.
    Hence, $\hat{S}_n(x \mid z) = 0$ almost everywhere on $B_z$, and there exists a subset $B_0 \subset B_z$ of positive measure such that
    $S_0(x \mid z) > 0$ for every $x \in B_0$.
    Using this fact and Proposition 3.8 of \cite{LastPenrose17}, we obtain
    \begin{align*}
        &\mathbb{E}\left[
            \int_A \log\left(\frac{S_0(x \mid z)}{\hat{S}_n(x \mid z)}\right) dX(x \mid z)
        \right] \\
        &= \mathbb{E}_N \left[
            \mathbb{E}_{x_1:x_N} \left[
                \sum_{j=1}^N
                \log\left(\frac{S_0(x_j \mid z)}{\hat{S}_n(x_j \mid z)}\right)
                \middle| N
            \right]
        \right] \\
        &= \mathbb{E}_N \left[
            \int_A \cdots \int_A
            \sum_{j=1}^N
            \log\left(\frac{S_0(x_j \mid z)}{\hat{S}_n(x_j \mid z)}\right)
            \prod_{j=1}^N \frac{S_0(x_j \mid z)}{\lambda_0(A \mid z)}
            \, d\nu(x_1) \cdots d\nu(x_N)
        \right] \\
        &\ge
        \mathbb{E}_N \left[
            \int_{B_0} \cdots \int_{B_0}
            \sum_{j=1}^N
            \log\left(\frac{S_0(x_j \mid z)}{\hat{S}_n(x_j \mid z)}\right)
            \prod_{j=1}^N \frac{S_0(x_j \mid z)}{\lambda_0(A \mid z)}
            \, d\nu(x_1) \cdots d\nu(x_N)
        \right]
        = \infty.
    \end{align*}
    Here, $N$ is a Poisson random variable with intensity parameter $\lambda_0(A \mid z)$.
    Therefore, together with the assumption that the covariate $Z$ has a positive density on its domain, the claim follows.
\end{proof}

\subsection{Proof of Lemma~\ref{lem: lower bound of the excess risk}}
\begin{proof}
    By Lemma~\ref{lem: the excess risk} and the inequality $\log x \le x-1$ for any $x>0$,
    \begin{align*}
        R(\hat{S}_n)
        &= \mathbb{E}_Z \left[
            \int_A S_0(x \mid Z)
            \log\left(\frac{S_0(x \mid Z)}{\hat{S}_n(x \mid Z)}\right)\, d\nu(
x)            + \int_A \bigl(\hat{S}_n(x \mid Z)-S_0(x \mid Z)\bigr)\, d\nu(x)
        \right] \\
        &= \mathbb{E}_Z \left[
            -2 \int_A S_0(x \mid Z)
            \log\left(\sqrt{\frac{\hat{S}_n(x \mid Z)}{S_0(x \mid Z)}}\right)\, d\nu(x)
            + \int_A \bigl(\hat{S}_n(x \mid Z)-S_0(x \mid Z)\bigr)\, d\nu(x)
        \right] \\
        &\ge \mathbb{E}_Z \left[
            -2 \int_A S_0(x \mid Z)
            \left(
                \sqrt{\frac{\hat{S}_n(x \mid Z)}{S_0(x \mid Z)}} - 1
            \right)\, d\nu(x)
            + \int_A \bigl(\hat{S}_n(x \mid Z)-S_0(x \mid Z)\bigr)\, d\nu(x)
        \right] \\
        &= \mathbb{E}_Z \left[
            \int_A
            \left(
                S_0(x \mid Z)
                - 2\sqrt{S_0(x \mid Z)\hat{S}_n(x \mid Z)}
                + \hat{S}_n(x \mid Z)
            \right)\, d\nu(x)
        \right]
        = 2\widetilde{H}^2(\hat{S}_n,S_0).
    \end{align*}
\end{proof}

\section{Proof of Theorem~\ref{thm: oracle inequality} and Corollary~\ref{cor: oracle inequality}}
In this section, we provide the proof of Theorem~\ref{thm: oracle inequality}, Corollary~\ref{cor: oracle inequality}, and technical lemmas to prove the theorem.

\subsection{Proof of Lemma~\ref{lem: basic inequality}}
\begin{proof}
    Let
    \[
        \bar{S}_n
        \coloneqq
        \frac{\hat{S}_n + S_n^*}{2}.
    \]
    First, by the concavity of the logarithm, for any
    $n \in \Nbb$, $z \in [0,1]^{d_z}$, and $x \in A_n^*(z)$,
    \begin{equation} \label{eq: basic-log-concavity}
        \frac{1}{2}
        \log\left(
            \frac{\bar{S}_n(x \mid z)}
            {S_n^*(x \mid z)}
        \right)
        =
        \frac{1}{2}
        \log\left(
            \frac{\hat{S}_n(x \mid z) + S_n^*(x \mid z)}
            {2S_n^*(x \mid z)}
        \right)
        \geq
        \frac{1}{4}
        \log\left(
            \frac{\hat{S}_n(x \mid z)}
            {S_n^*(x \mid z)}
        \right).
    \end{equation}

    Since $\hat{S}_n$ is the maximum likelihood estimator and $S_n^* \in \mathcal{F}_n$, we have
    \begin{align}
        0
        &\leq
        P_n
        \left\{
            \int_A
            \log\left(
                \frac{\hat{S}_n(x \mid Z)}
                {S_n^*(x \mid Z)}
            \right)
            dX(x \mid Z)
            -
            \int_A
            \left(
                \hat{S}_n(x \mid Z)
                -
                S_n^*(x \mid Z)
            \right)
            d\nu(x)
        \right\}.
        \label{eq: mle-likelihood-ratio}
    \end{align}
    By Assumption~\ref{asm: weak absolute continuity},
    \[
        S_n^*(x \mid z)=0
        \quad\Longrightarrow\quad
        S_0(x \mid z)=0
        \quad
        \text{for } \nu\text{-a.e. }x .
    \]
    Hence,
    \[
        \mathbb{E}
        \left[
            X_i\left(A \setminus A_n^*(z_i) \mid z_i\right)
            \mid Z_i=z_i
        \right]
        =
        \int_{A \setminus A_n^*(z_i)}
        S_0(x \mid z_i)\,d\nu(x)
        =
        0.
    \]
    Therefore, since $X_i(A \setminus A_n^*(z_i) \mid z_i)$ is nonnegative random variable,
    \[
        X_i\left(A \setminus A_n^*(z_i) \mid z_i\right)=0
        \quad
        \text{a.s.}
    \]
    Consequently, the stochastic integral in \eqref{eq: mle-likelihood-ratio} may be restricted to
    $A_n^*(Z)$ almost surely. Moreover, on $A \setminus A_n^*(Z)$, we have
    $S_n^*(x \mid Z)=0$. Thus,
    \begin{align}
        0
        &\leq
        P_n
        \left\{
            \int_{A_n^*(Z)}
            \log\left(
                \frac{\hat{S}_n(x \mid Z)}
                {S_n^*(x \mid Z)}
            \right)
            dX(x \mid Z)
        \right. \nonumber \\
        &\quad\left.
            +
            \int_{A_n^*(Z)}
            \left(
                S_n^*(x \mid Z)
                -
                \hat{S}_n(x \mid Z)
            \right)
            d\nu(x)
            -
            \int_{A \setminus A_n^*(Z)}
            \hat{S}_n(x \mid Z)\,d\nu(x)
        \right\}.
        \label{eq: restricted-likelihood-ratio}
    \end{align}
    Multiplying both sides by $1/4$ and using \eqref{eq: basic-log-concavity}, we obtain
    \begin{align}
        0
        &\leq
        P_n
        \left\{
            \int_{A_n^*(Z)}
            \frac{1}{2}
            \log\left(
                \frac{\bar{S}_n(x \mid Z)}
                {S_n^*(x \mid Z)}
            \right)
            dX(x \mid Z)
        \right. \nonumber \\
        &\quad\left.
            +
            \frac{1}{2}
            \int_{A_n^*(Z)}
            \left(
                S_n^*(x \mid Z)
                -
                \bar{S}_n(x \mid Z)
            \right)
            d\nu(x)
        \right.
        \left.
            -
            \frac{1}{2}
            \int_{A \setminus A_n^*(Z)}
            \bar{S}_n(x \mid Z)\,d\nu(x)
        \right\}.
        \label{eq: basic-pn}
    \end{align}
    Here we used
    \[
        \frac{1}{4}
        \left(
            S_n^*-\hat{S}_n
        \right)
        =
        \frac{1}{2}
        \left(
            S_n^*-\bar{S}_n
        \right)
        \quad
        \text{on } A_n^*(Z),
    \]
    and, since $S_n^*=0$ on $A \setminus A_n^*(Z)$,
    \[
        \frac{1}{4}\hat{S}_n
        =
        \frac{1}{2}\bar{S}_n
        \quad
        \text{on } A \setminus A_n^*(Z).
    \]

    Define
    \begin{align*}
        G_n(Z,X)
        &\coloneqq
        \int_{A_n^*(Z)}
        \frac{1}{2}
        \log\left(
            \frac{\bar{S}_n(x \mid Z)}
            {S_n^*(x \mid Z)}
        \right)
        dX(x \mid Z) \\
        &\quad
        +
        \frac{1}{2}
        \int_{A_n^*(Z)}
        \left(
            S_n^*(x \mid Z)
            -
            \bar{S}_n(x \mid Z)
        \right)
        d\nu(x)
        -
        \frac{1}{2}
        \int_{A \setminus A_n^*(Z)}
        \bar{S}_n(x \mid Z)\,d\nu(x).
    \end{align*}
    Then \eqref{eq: basic-pn} gives
    \[
        0
        \leq
        P_n G_n
        =
        (P_n-P)G_n + PG_n.
    \]

    By Campbell's formula,
    \begin{align*}
        PG_n
        &=
        \mathbb{E}_Z
        \left\{
            \int_{A_n^*(Z)}
            \frac{1}{2}
            \log\left(
                \frac{\bar{S}_n(x \mid Z)}
                {S_n^*(x \mid Z)}
            \right)
            S_0(x \mid Z)\,d\nu(x)
        \right. \\
        &\quad\left.
            +
            \frac{1}{2}
            \int_{A_n^*(Z)}
            \left(
                S_n^*(x \mid Z)
                -
                \bar{S}_n(x \mid Z)
            \right)
            d\nu(x)
        \right.
        \left.
            -
            \frac{1}{2}
            \int_{A \setminus A_n^*(Z)}
            \bar{S}_n(x \mid Z)\,d\nu(x)
        \right\}.
    \end{align*}
    The application of Campbell's formula is justified as follows. Since
    $\log(\bar{S}_n / S_n^*) \geq -\log 2$
    on $A_n^*(z)$, the inequality $\log x \leq x-1$ for $x>0$, together with
    Assumption~\ref{asm: weak absolute continuity}, implies that, for any
    $z \in [0,1]^{d_z}$,
    \begin{align*}
        &\int_{A_n^*(z)}
        \left\lvert
            \log\left(
                \frac{\bar{S}_n(x \mid z)}
                {S_n^*(x \mid z)}
            \right)
        \right\rvert
        S_0(x \mid z)\,d\nu(x) \\
        &=
        \int_{A_n^*(z)}
        \left\lvert
            \log\left(
                \frac{\bar{S}_n(x \mid z)}
                {S_n^*(x \mid z)}
            \right)
        \right\rvert
        \frac{S_0(x \mid z)}
        {S_n^*(x \mid z)}
        S_n^*(x \mid z)\,d\nu(x) \\
        &\leq
        c_0^2
        \int_{A_n^*(z)}
        \left[
            \log(2)S_n^*(x \mid z)
            +
            \left\lvert
                \bar{S}_n(x \mid z)
                -
                S_n^*(x \mid z)
            \right\rvert
        \right]d\nu(x)
        <
        \infty.
    \end{align*}

    Next, for any $z \in [0,1]^{d_z}$, using $\log x \leq x-1$ with
    $x=\sqrt{\bar{S}_n/S_n^*}$ gives
    \begin{align*}
        &\int_{A_n^*(z)}
        \frac{1}{2}
        \log\left(
            \frac{\bar{S}_n(x \mid z)}
            {S_n^*(x \mid z)}
        \right)
        S_0(x \mid z)\,d\nu(x)
        \leq
        -
        \int_{A_n^*(z)}
        \left(
            1
            -
            \sqrt{
                \frac{\bar{S}_n(x \mid z)}
                {S_n^*(x \mid z)}
            }
        \right)
        S_0(x \mid z)\,d\nu(x).
    \end{align*}
    We decompose $S_0$ as $S_0=S_n^*+(S_0-S_n^*)$ on $A_n^*(z)$. Then
    \begin{align*}
        &-
        \int_{A_n^*(z)}
        \left(
            1
            -
            \sqrt{
                \frac{\bar{S}_n(x \mid z)}
                {S_n^*(x \mid z)}
            }
        \right)
        S_0(x \mid z)\,d\nu(x) \\
        &=
        -
        \int_{A_n^*(z)}
        \left(
            S_n^*(x \mid z)
            -
            \sqrt{S_n^*(x \mid z)}
            \sqrt{\bar{S}_n(x \mid z)}
        \right)d\nu(x) \\
        &\quad
        -
        \int_{A_n^*(z)}
        \left(
            \sqrt{S_n^*(x \mid z)}
            -
            \sqrt{\bar{S}_n(x \mid z)}
        \right)
        \left(
            \sqrt{
                \frac{S_0(x \mid z)}
                {S_n^*(x \mid z)}
            }
            +
            1
        \right)
        \left(
            \sqrt{S_0(x \mid z)}
            -
            \sqrt{S_n^*(x \mid z)}
        \right)d\nu(x).
    \end{align*}
    By Assumption~\ref{asm: weak absolute continuity},
    \[
        \sqrt{
            \frac{S_0(x \mid z)}
            {S_n^*(x \mid z)}
        }
        \leq
        c_0
        \quad
        \text{on } A_n^*(z).
    \]
    Hence, by the Cauchy--Schwarz inequality,
    \begin{align*}
        &-
        \int_{A_n^*(z)}
        \left(
            \sqrt{S_n^*(x \mid z)}
            -
            \sqrt{\bar{S}_n(x \mid z)}
        \right)
        \left(
            \sqrt{
                \frac{S_0(x \mid z)}
                {S_n^*(x \mid z)}
            }
            +
            1
        \right) 
        \left(
            \sqrt{S_0(x \mid z)}
            -
            \sqrt{S_n^*(x \mid z)}
        \right)d\nu(x) \\
        &\leq
        2(c_0+1)
        H\left(
            \bar{S}_n(\cdot \mid z),
            S_n^*(\cdot \mid z)
        \right)
        H\left(
            S_n^*(\cdot \mid z),
            S_0(\cdot \mid z)
        \right).
    \end{align*}
    Therefore,
    \begin{align*}
        &\int_{A_n^*(z)}
        \frac{1}{2}
        \log\left(
            \frac{\bar{S}_n(x \mid z)}
            {S_n^*(x \mid z)}
        \right)
        S_0(x \mid z)\,d\nu(x) \\
        &\leq
        -
        \int_{A_n^*(z)}
        \left(
            S_n^*(x \mid z)
            -
            \sqrt{S_n^*(x \mid z)}
            \sqrt{\bar{S}_n(x \mid z)}
        \right)d\nu(x) \\
        &\quad
        +
        2(c_0+1)
        H\left(
            \bar{S}_n(\cdot \mid z),
            S_n^*(\cdot \mid z)
        \right)
        H\left(
            S_n^*(\cdot \mid z),
            S_0(\cdot \mid z)
        \right).
    \end{align*}

    Combining this bound with the second integral term in the definition of \(G_n\), we obtain
    \begin{align*}
        &-
        \int_{A_n^*(z)}
        \left(
            S_n^*(x \mid z)
            -
            \sqrt{S_n^*(x \mid z)}
            \sqrt{\bar{S}_n(x \mid z)}
        \right)d\nu(x) 
        +
        \frac{1}{2}
        \int_{A_n^*(z)}
        \left(
            S_n^*(x \mid z)
            -
            \bar{S}_n(x \mid z)
        \right)d\nu(x) \\
        &=
        -
        \frac{1}{2}
        \int_{A_n^*(z)}
        \left(
            \sqrt{\bar{S}_n(x \mid z)}
            -
            \sqrt{S_n^*(x \mid z)}
        \right)^2d\nu(x).
    \end{align*}
    Moreover, since $S_n^*(x \mid z)=0$ on $A \setminus A_n^*(z)$,
    \begin{align*}
        -
        \frac{1}{2}
        \int_{A \setminus A_n^*(z)}
        \bar{S}_n(x \mid z)\,d\nu(x)
        &=
        -
        \frac{1}{2}
        \int_{A \setminus A_n^*(z)}
        \left(
            \sqrt{\bar{S}_n(x \mid z)}
            -
            \sqrt{S_n^*(x \mid z)}
        \right)^2d\nu(x).
    \end{align*}
    Hence,
    \begin{align*}
        PG_n
        &\leq
        -
        \mathbb{E}_Z
        \left[
            H^2\left(
                \bar{S}_n(\cdot \mid Z),
                S_n^*(\cdot \mid Z)
            \right)
        \right]
        +
        2(c_0+1)
        \mathbb{E}_Z
        \left[
            H\left(
                \bar{S}_n(\cdot \mid Z),
                S_n^*(\cdot \mid Z)
            \right)
            H\left(
                S_n^*(\cdot \mid Z),
                S_0(\cdot \mid Z)
            \right)
        \right] \\
        &=
        -
        \widetilde{H}^2\left(
            \bar{S}_n,
            S_n^*
        \right)
        +
        2(c_0+1)
        \mathbb{E}_Z
        \left[
            H\left(
                \bar{S}_n(\cdot \mid Z),
                S_n^*(\cdot \mid Z)
            \right)
            H\left(
                S_n^*(\cdot \mid Z),
                S_0(\cdot \mid Z)
            \right)
        \right].
    \end{align*}
    By Hölder's inequality,
    \begin{equation*}
        \mathbb{E}_Z
        \left[
            H\left(
                \bar{S}_n(\cdot \mid Z),
                S_n^*(\cdot \mid Z)
            \right)
            H\left(
                S_n^*(\cdot \mid Z),
                S_0(\cdot \mid Z)
            \right)
        \right]
        \leq
        \widetilde{H}\left(
            \bar{S}_n,
            S_n^*
        \right)
        \widetilde{H}\left(
            S_n^*,
            S_0
        \right).
    \end{equation*}
    Since \(S_n^*(x\mid Z)=0\) on \(A\setminus A_n^*(Z)\), the integral terms in \(G_n\) can be combined as
\[
    \frac{1}{2}
    \int_{A_n^*(Z)}
    \left(
        S_n^*(x\mid Z)-\bar S_n(x\mid Z)
    \right)d\nu(x)
    -
    \frac{1}{2}
    \int_{A\setminus A_n^*(Z)}
    \bar S_n(x\mid Z)d\nu(x)
    =
    \frac{1}{2}
    \int_A
    \left(
        S_n^*(x\mid Z)-\bar S_n(x\mid Z)
    \right)d\nu(x).
\]
This completes the proof.
\end{proof}

\subsection{Technical lemmas}
In this section, we provide several definitions, technical lemmas, and their proof.

\begin{lemma} \label{lem: hellinger of convex combination}
    For any conditional intensity functions $S_1$ and $S_2$, we have
    \begin{equation*}
        \widetilde{H}^2(S_1, S_2)
        \leq
        16 \widetilde{H}^2 \left( \frac{S_1 + S_2}{2}, S_2 \right).
    \end{equation*}
\end{lemma}

\begin{proof}
    By the subadditivity of the square-root function, for any $a, b \ge 0$, we have
    \begin{equation*}
        \sqrt{\frac{a + b}{2}}
        + \sqrt{b}
        \leq
        2\left(\sqrt{a} + \sqrt{b}\right).
    \end{equation*}
    Therefore, by simple calculation, we have
    \begin{equation*}
        |\sqrt{a} - \sqrt{b}| \le \frac{2|(a+b)/2 - b|}{\sqrt{a} + \sqrt{b}}
        = 2 \cdot \frac{\sqrt{(a+b)/2} + \sqrt{b}}{\sqrt{a} + \sqrt{b}} \left|\sqrt{\frac{a + b}{2}} - \sqrt{b}\right|
        \le 4\left|\sqrt{\frac{a + b}{2}} - \sqrt{b}\right|.
    \end{equation*}
    Substituting $a=S_1,\, b=S_2$, squaring both sides, integrating over $A$, and taking expectation with respect to $Z$ gives the claim.
\end{proof}

\begin{Def}[Bernstein size]
    Let $P$ be a probability measure on a measurable space $(S,\Scal)$, and let $K>0$ be a constant.
    For a measurable function $f$, define
    \begin{equation*}
        \rho_K^2(f)
        =
        2K^2 \Ebb \left[
            e^{\lvert f(X) \rvert/K}
            - 1
            - \frac{\lvert f(X) \rvert}{K}
        \right],
    \end{equation*}
    where $X$ is a random variable with distribution $P$.
    The positive square root of $\rho_K^2(f)$ is called the Bernstein size of $f$.
\end{Def}

Using this Bernstein size, we define the generalized bracketing number as follows.

\begin{Def}[Generalized bracketing number]
    Let $\Fcal$ be a class of functions, let $\delta>0$, and let $P$ be a probability measure.
    Let $\Ncal_{B,K}(\delta,\Fcal,P)$ be the smallest integer $N$ for which there exists a collection of pairs of functions $\{(f_j^L,f_j^U)\}_{j=1}^N$ satisfying the following conditions:
    \begin{itemize}
        \item For every $j=1,\dots,N$, we have $\rho_K(f_j^U-f_j^L)\leq \delta$.
        \item For every $f \in \Fcal$, there exists an index $j=j(f)\in\{1,\dots,N\}$ such that
        $
            f_j^L \leq f \leq f_j^U.
        $
    \end{itemize}
    If no such $N$ exists, we set $\Ncal_{B,K}(\delta,\Fcal,P)=\infty$.
    We call $\Ncal_{B,K}(\delta,\Fcal,P)$ the $\delta$-generalized bracketing number of $\Fcal$.
\end{Def}

The following Bernstein-type uniform inequality will be used to control the empirical process.

\begin{theorem}[Theorem 5.11 of \citealp{vandeGeer00}; Theorem 3.5.21 of \citealp{GineNickl2021}] \label{thm: uniform bernstein inequality}
    Let $\Gcal$ be a class of functions.
    Suppose that, for some constants $K$ and $R$, the class $\Gcal$ satisfies
    \begin{equation}
        \sup_{g \in \Gcal} \rho_K(g) \leq R. \label{eq: bernstein condition}
    \end{equation}
    If
    \begin{align}
        a &\leq C_1\sqrt{n}R^2/K, \label{eq: condition 1 for bernstein} \\
        a &\leq 8\sqrt{n}R, \label{eq: condition 2 for bernstein} \\
        a &\geq C_0\left(
            \int_{a/(2^6\sqrt{n})}^R
            \sqrt{\log\mathcal{N}_{B,K}(u,\Gcal,P)} \, du
            \lor R
        \right), \label{eq:condition 3 for bernstein} \\
        C_0^2 &\geq C^2(C_1+1), \label{eq: condition 4 for bernstein}
    \end{align}
    then
    \begin{equation*}
        \Pbb\left(
            \sqrt{n}\sup_{g \in \Gcal}
            \lvert (P_n-P)g \rvert
            \geq a
        \right)
        \leq
        C\exp\left(
            -\frac{a^2}{C^2(C_1+1)R^2}
        \right).
    \end{equation*}
\end{theorem}

\begin{lemma}[Lemma 7.1 of \citealp{vandeGeer00}] \label{lem: bound for bernstein semi-norm}
    For any $x \geq -\Gamma$ with $\Gamma > 0$, we have
    \begin{equation*}
        2(e^{\lvert x \rvert} - 1 - \lvert x \rvert)
        \leq
        c_\Gamma^2(e^x - 1)^2,
        \quad
        c_\Gamma =
        \frac{2(e^\Gamma - 1 - \Gamma)}{(e^{-\Gamma} - 1)^2}.
    \end{equation*}
\end{lemma}

\begin{lemma} \label{lem: upper bound for difference of product}
    For any $n \in \Nbb$ and any bounded real numbers
    $x_1,\dots,x_n,y_1,\dots,y_n$ satisfying
    $\lvert x_i \rvert \leq M$ and $\lvert y_i \rvert \leq M$ for some $M \geq 0$, we have
    \begin{equation*}
        \left(
            \prod_{i=1}^n x_i
            -
            \prod_{i=1}^n y_i
        \right)^2
        \leq
        M^{2n-2}n
        \sum_{i=1}^n (x_i-y_i)^2.
    \end{equation*}
\end{lemma}

\begin{proof}
    We first prove, by induction, that
    \begin{equation*}
        \left\lvert
            \prod_{i=1}^n x_i
            -
            \prod_{i=1}^n y_i
        \right\rvert
        \leq
        M^{n-1}
        \sum_{i=1}^n \lvert x_i-y_i\rvert.
    \end{equation*}
    The claim is obvious when $n=1$.
    Suppose that $n \geq 2$ and that the claim holds for $n-1$.
    Then, by the boundedness of $x_i$ and $y_i$, the triangle inequality, and the induction hypothesis,
    \begin{align*}
        \left\lvert
            \prod_{i=1}^{n} x_i
            -
            \prod_{i=1}^{n} y_i
        \right\rvert
        &=
        \left\lvert
            x_n \prod_{i=1}^{n-1} x_i
            -
            y_n \prod_{i=1}^{n-1} x_i
            +
            y_n \prod_{i=1}^{n-1} x_i
            -
            y_n \prod_{i=1}^{n-1} y_i
        \right\rvert \\
        &\leq
        \lvert x_n-y_n\rvert
        \prod_{i=1}^{n-1} \lvert x_i\rvert
        +
        \lvert y_n\rvert
        \left\lvert
            \prod_{i=1}^{n-1} x_i
            -
            \prod_{i=1}^{n-1} y_i
        \right\rvert 
        \leq
        M^{n-1}\lvert x_n-y_n\rvert
        +
        M^{n-1}
        \sum_{i=1}^{n-1} \lvert x_i-y_i\rvert \\
        &\leq
        M^{n-1}
        \sum_{i=1}^{n} \lvert x_i-y_i\rvert.
    \end{align*}
    This proves the desired bound.
    Therefore, using the inequality
    \[
    \left(\sum_{j=1}^n a_j\right)^2
    \leq
    n\sum_{j=1}^n a_j^2
    \]
    for arbitrary nonnegative real numbers $a_j \geq 0$, $j=1,\dots,n$, we obtain
    \begin{equation*}
        \left(
            \prod_{i=1}^n x_i
            -
            \prod_{i=1}^n y_i
        \right)^2
        \leq
        M^{2n-2}n
        \sum_{i=1}^n (x_i-y_i)^2.
    \end{equation*}
\end{proof}

\begin{lemma} \label{lem: hellinger upper bound for bernstein semi-norm}
    Suppose that Assumption \ref{asm: weak absolute continuity} holds and \(\sup_{S \in \Fcal_n} \lVert S \rVert_\infty \le F < \infty\).
    Then, for any $g_S \in \Gcal_n$, we have
    \begin{equation*}
        \rho_1^2(g_S)
        \leq
        2^6 C_{\rho}^2
        \widetilde{H}^2\left(
            \frac{S+S_n^*}{2},
            S_n^*
        \right).
    \end{equation*}
    Here, the constant $C_{\rho}$ is given by
    \begin{equation*}
        C_{\rho}
        =
        \exp \left(
            2c_0^2
            (\nu(A) \lor 1)
            \left(
                F
                \lor 1
            \right)
        \right).
    \end{equation*}
\end{lemma}

\begin{proof}
    Fix $z$ and put $A^*=A_n^*(z)$, $\lambda=\lambda_0(A^*\mid z)$,
$V=\nu(A^*)$, and
$H_z^2=H^2(\bar S(\cdot\mid z),S_n^*(\cdot\mid z))$.
If $\lambda=0$, then $X(A^*)=0$ almost surely and $g_S=0$.
Suppose henceforth that $\lambda>0$. Conditional on $N=X(A^*)\ge1$,
the locations have density $S_0/\lambda$ and
$g_S\ge-N\log2/2$. With $\Gamma=N\log2/2$, the preceding two lemmas
and $S_0\le c_0^2S_n^*$ therefore give
\begin{align*}
 &\Ebb[2(e^{|g_S|}-1-|g_S|)\mid Z=z,N]\\
 &\le \frac{c_\Gamma^2c_0^{2N}}{\lambda^N}
   \int_{(A^*)^N}
   \left(\prod_{j=1}^N\sqrt{\bar S(x_j\mid z)}
                  -\prod_{j=1}^N\sqrt{S_n^*(x_j\mid z)}\right)^2
       d\nu^{\otimes N}\\
 &\le \frac{2c_\Gamma^2c_0^{2N}N^2(FV)^{N-1}}{\lambda^N}H_z^2.
\end{align*}
For $b=e^\Gamma>1$, $\log b\ge1-b^{-1}$ implies
\[
 c_\Gamma=
 \frac{2(b-1-\log b)}{(1-b^{-1})^2}\le2b,
 \qquad c_\Gamma^2\le4\,2^N.
\]
Put $M=(\nu(A)\vee1)(F\vee1)\ge1$ and $x=2c_0^2M$.
The contribution of $N=0$ is zero. Summing the preceding bound over
$N\sim\operatorname{Poisson}(\lambda)$ yields
\begin{align*}
 \Ebb[2(e^{|g_S|}-1-|g_S|)\mid Z=z]
 &\le \frac8M e^{-\lambda}H_z^2
                  \sum_{N\ge1}\frac{N^2x^N}{N!}\\
 &=\frac8M(x^2+x)e^{x-\lambda}H_z^2
 \le8e^{2x}H_z^2\le64C_\rho^2H_z^2.
\end{align*}
Here $(x^2+x)\le e^x$ for $x\ge0$; for example this follows from
$e^x\ge1+x+x^2/2+x^3/6$ and
$1-x^2/2+x^3/6\ge1/3$. Also $e^{2x}=C_\rho^2$.
The last coefficient is independent of $z$. Taking expectation over
$Z$ now proves the result without separating dependent factors.
\end{proof}

\begin{lemma} \label{lem: relation of bracketing between bernstein semi-norm and hellinger distance}
    Suppose Assumption \ref{asm: weak absolute continuity} holds and \(\sup_{S \in \Fcal_n} \lVert S \rVert_\infty = F < \infty\).
    Let $S_L$ and $S_U$ satisfy
    $
    0 \leq S_L(x \mid z)
    \leq
    S_U(x\mid z)
    \leq
    F.
    $
    pointwise.
    Define
    \begin{align*}
        g_L
        &=
        \int_{A_n^*(Z)}
        \frac{1}{2}
        \log\left(
            \frac{
                S_L(x \mid Z)
                +
                S_n^*(x \mid Z)
            }{
                2S_n^*(x \mid Z)
            }
        \right)
        dX(x \mid Z), \\
        g_U
        &=
        \int_{A_n^*(Z)}
        \frac{1}{2}
        \log\left(
            \frac{
                S_U(x \mid Z)
                +
                S_n^*(x \mid Z)
            }{
                2S_n^*(x \mid Z)
            }
        \right)
        dX(x \mid Z).
    \end{align*}
    Then
    \begin{equation*}
        \rho_1^2(g_U-g_L)
        \leq
        2C_\rho^2
        \widetilde{H}^2\left(
            \frac{S_U+S_n^*}{2},
            \frac{S_L+S_n^*}{2}
        \right).
    \end{equation*}
    Here, $C_\rho$ is the same constant as in Lemma
    \ref{lem: hellinger upper bound for bernstein semi-norm}.
\end{lemma}

\begin{proof}
Put $\bar S_j=(S_j+S_n^*)/2$ for $j\in\{L,U\}$ and
$h=g_U-g_L\ge0$.
Fix $z$ and write
\[
    A^*=A_n^*(z),\qquad
    \lambda=\lambda_0(A^*\mid z),\qquad
    V=\nu(A^*),\qquad
    H_z^2
    =H^2\bigl(\bar S_U(\cdot\mid z),\bar S_L(\cdot\mid z)\bigr).
\]
If $\lambda=0$, then $h=0$ almost surely conditional on $Z=z$.
Suppose that $\lambda>0$ and condition on $N=X(A^*)\ge1$.
Since $2(e^u-1-u)\le(e^u-1)^2$ for $u\ge0$ and
$S_0/\bar S_L\le2c_0^2$ on $A^*$, the conditional distribution
of the points and
Lemma~\ref{lem: upper bound for difference of product} give
\begin{align*}
    &\Ebb\bigl[2(e^h-1-h)\mid Z=z,N\bigr]\\
    &\le
    \frac{(2c_0^2)^N}{\lambda^N}
    \int_{(A^*)^N}
    \left(
        \prod_{j=1}^N\sqrt{\bar S_U(x_j\mid z)}
        -
        \prod_{j=1}^N\sqrt{\bar S_L(x_j\mid z)}
    \right)^2\,d\nu^{\otimes N}\\
    &\le
    \frac{2(2c_0^2)^N N^2(FV)^{N-1}}{\lambda^N}H_z^2.
\end{align*}
Let $M=(\nu(A)\vee1)(F\vee1)$ and $x=2c_0^2M$.
The contribution of $N=0$ is zero. Since
\[
    \sum_{N\ge1}\frac{N^2x^N}{N!}
    =(x^2+x)e^x
    \quad\text{and}\quad
    x^2+x\le e^x\quad(x\ge0),
\]
summing over $N\sim\operatorname{Poisson}(\lambda)$ yields
\begin{align*}
    \Ebb\bigl[2(e^h-1-h)\mid Z=z\bigr]
    &\le
    \frac{2}{M}e^{-\lambda}H_z^2
    \sum_{N\ge1}\frac{N^2x^N}{N!}\\
    &=
    \frac{2}{M}(x^2+x)e^{x-\lambda}H_z^2\\
    &\le
    2e^{2x}H_z^2
    =
    2C_\rho^2H_z^2.
\end{align*}
The coefficient is independent of $z$.
Taking expectation over $Z$ proves the result.
\end{proof}

\begin{lemma}[Lemma 3.5.18 of \citealp{GineNickl2021}] \label{lem: properties of Bernstein size}
    The Bernstein size $\rho_K(\cdot)$ satisfies the following properties:
    \begin{itemize}
        \item $\rho_K(f)$ is nonincreasing in $K$.
        \item If $\lVert f \rVert_\infty \leq K$ and $\Ebb[f(X)^2] \leq R^2$, then $\rho_{2K}^2(f) \leq 2R^2.$
    \end{itemize}
\end{lemma}

\begin{lemma} \label{lem: Hellinger upper bound of Bernstein size 2}
    Suppose that \(\sup_{S \in \Fcal_n} \lVert S \rVert_\infty = F < \infty\).
    For any \(g_S^\prime \in \Gcal_n^\prime\), we have
    \[
        \rho_{F\nu(A)}^2(g_S^\prime)
        \leq
        4F\nu(A)
        \widetilde H^2
        \left(
            S_n^*,
            \frac{S+S_n^*}{2}
        \right).
    \]
\end{lemma}

\begin{proof}
    Since
    \(0 \leq S_n^* \leq F\) and \(0 \leq (S + S_n^*) / 2 \leq F\), we have
    \[
        \lvert g_S^\prime \rvert
        \leq
        \frac{1}{2}
        \int_A
        \left\lvert S_n^*(x \mid Z)-\frac{S(x \mid Z) + S_n^*(x \mid Z)}{2} \right\rvert
        d\nu(x)
        \leq
        \frac{F\nu(A)}{2}.
    \]
    Moreover, by Jensen's inequality,
    \begin{align*}
        \mathbb E[(g_S^\prime)^2]
        &\leq
        \frac{\nu(A)}{4}
        \mathbb E_Z
        \left[
            \int_A
            \left(
                S_n^*(x \mid Z)-\frac{S(x \mid Z) + S_n^*(x \mid Z)}{2}
            \right)^2
            d\nu(x)
        \right].
    \end{align*}
    Since \(0 \leq S_n^*, (S(x \mid Z) + S_n^*(x \mid Z))/2 \leq F\), we have
    \begin{align*}
        \left(
            S_n^*-\frac{S(x \mid Z) + S_n^*(x \mid Z)}{2}
        \right)^2
        &=
        \left(
            \sqrt{S_n^*}-\sqrt{\frac{S(x \mid Z) + S_n^*(x \mid Z)}{2}}
        \right)^2
        \left(
            \sqrt{S_n^*}+\sqrt{\frac{S(x \mid Z) + S_n^*(x \mid Z)}{2}}
        \right)^2 \\
        &\leq
        4F
        \left(
            \sqrt{S_n^*}-\sqrt{\frac{S(x \mid Z) + S_n^*(x \mid Z)}{2}}
        \right)^2.
    \end{align*}
    Hence,
    \[
        \mathbb E[(g_S^\prime)^2]
        \leq
        2F\nu(A)
        \widetilde H^2(S_n^*,\bar S).
    \]
    By Lemma~\ref{lem: properties of Bernstein size}, with \(K=F\nu(A) / 2\), we obtain
    \[
        \rho_{F\nu(A)}^2(g_S^\prime)
        \leq
        4F\nu(A)
        \widetilde H^2(S_n^*,\bar S).
    \]
\end{proof}

\begin{lemma} \label{lem: relation of bracketing between bernstein size and hellinger distance 2}
    Suppose that \(\sup_{S \in \Fcal_n}\lVert S\rVert_\infty \le F < \infty\).
    For any \(0 \leq S_L \leq S_U \leq F\), define
    \[
        g_L^\prime
        \coloneqq
        \frac{1}{2}
        \int_A
        \left(
            S_n^*(x \mid Z)-\frac{S_U(x \mid Z) + S_n^*(x \mid Z)}{2}
        \right)d\nu(x),
    \]
    and
    \[
        g_U^\prime
        \coloneqq
        \frac{1}{2}
        \int_A
        \left(
            S_n^*(x \mid Z)-\frac{S_L(x \mid Z) + S_n^*(x \mid Z)}{2}
        \right)d\nu(x).
    \]
    Then
    \[
        \rho_{F\nu(A)}^2(g_U^\prime-g_L^\prime)
        \leq
        4F\nu(A)
        \widetilde H^2
        \left(
            \frac{S_U+S_n^*}{2},
            \frac{S_L+S_n^*}{2}
        \right).
    \]
\end{lemma}

\begin{proof}
    We have
    \[
        g_U^\prime-g_L^\prime
        =
        \frac{1}{4}
        \int_A
        \left(
             S_U(x \mid Z)- S_L(x \mid Z)
        \right)d\nu(x).
    \]
    Hence,
    \[
        0
        \leq
        g_U^\prime-g_L^\prime
        \leq
        \frac{F\nu(A)}{4}.
    \]
    Moreover, by Jensen's inequality and the inequality
    \[
        \left(\frac{S_U+S_n^*}{2}-\frac{S_L+S_n^*}{2}\right)^2
        \leq
        4F
        \left(
            \sqrt{\frac{S_U+S_n^*}{2}}-\sqrt{\frac{S_L+S_n^*}{2}}
        \right)^2,
    \]
    we obtain
    \[
        \mathbb E[(g_U^\prime-g_L^\prime)^2]
        \leq
        2F\nu(A)
        \widetilde H^2
        \left(
            \frac{S_U+S_n^*}{2},
            \frac{S_L+S_n^*}{2}
        \right).
    \]
    Therefore, Lemma~\ref{lem: properties of Bernstein size} gives
    \[
        \rho_{F\nu(A)}^2(g_U^\prime-g_L^\prime)
        \leq
        4F\nu(A)
        \widetilde H^2
        \left(
            \frac{S_U+S_n^*}{2},
            \frac{S_L+S_n^*}{2}
        \right).
    \]
\end{proof}

\subsection{Main proof}
\begin{proof}{(Proof of Theorem~\ref{thm: oracle inequality})}
    Put
    \[
        \bar S_n
        \coloneqq
        \frac{\hat S_n+S_n^*}{2}.
    \]
    For \(S\in\Fcal_n\), define
    \[
        \bar S
        \coloneqq
        \frac{S+S_n^*}{2},
    \]
    and let
    \[
        g_S
        \coloneqq
        \int_{A_n^*(Z)}
        \frac{1}{2}
        \log\left(
            \frac{\bar S(x\mid Z)}
            {S_n^*(x\mid Z)}
        \right)
        dX(x\mid Z),
    \]
    and
    \[
        g_S^\prime
        \coloneqq
        \frac{1}{2}
        \int_A
        \left(
            S_n^*(x\mid Z)-\bar S(x\mid Z)
        \right)d\nu(x).
    \]
    By Lemma~\ref{lem: basic inequality}, we have
    \begin{align*}
        \widetilde H^2(\bar S_n,S_n^*)
        &\leq
        (P_n-P)g_{\hat S_n}
        +
        (P_n-P)g_{\hat S_n}^\prime
        +
        2(c_0+1)
        \widetilde H(\bar S_n,S_n^*)
        \widetilde H(S_n^*,S_0).
    \end{align*}
    Denote
        $r_n
        \coloneqq
        \widetilde H(\bar S_n,S_n^*),\,
        h_n
        \coloneqq
        \widetilde H(S_n^*,S_0)$.

    First consider the case in which
    \[
        \max\{(P_n-P)g_{\hat S_n},(P_n-P)g_{\hat S_n}^\prime,2(c_0+1)r_nh_n\}
        =
        2(c_0+1)r_nh_n.
    \]
    Then, we obtain
    \[
        r_n^2
        \leq
        6(c_0+1)r_nh_n.
    \]
    If \(r_n>0\), this implies
    \[
        r_n^2
        \leq
        36(c_0+1)^2h_n^2.
    \]
    If \(r_n=0\), the desired bound is immediate. Therefore, by the triangle inequality and Lemma~\ref{lem: hellinger of convex combination},
    \begin{align*}
        \widetilde H^2(\hat S_n,S_0)
        &\leq
        2\left\{
            \widetilde H^2(\hat S_n,S_n^*)
            +
            \widetilde H^2(S_n^*,S_0)
        \right\} \\
        &\leq
        2\left\{
            16\widetilde H^2(\bar S_n,S_n^*)
            +
            \widetilde H^2(S_n^*,S_0)
        \right\} \\
        &\leq
        1154(c_0+1)^2
        \widetilde H^2(S_n^*,S_0).
    \end{align*}
    Hence the theorem follows in this case.

    It remains to consider the cases in which either \(\mathrm{I}\) or \(\mathrm{II}\) is the largest term. If
    \[
        \max\{(P_n-P)g_{\hat S_n},(P_n-P)g_{\hat S_n}^\prime,2(c_0+1)r_nh_n\}
        =
        (P_n-P)g_{\hat S_n},
    \]
    then
    \[
        r_n^2
        \leq
        3(P_n-P)g_{\hat S_n}.
    \]
    If
    \[
        \max\{(P_n-P)g_{\hat S_n},(P_n-P)g_{\hat S_n}^\prime,2(c_0+1)r_nh_n\}
        =
        (P_n-P)g_{\hat S_n}^\prime,
    \]
    then
    \[
        r_n^2
        \leq
        3(P_n-P)g_{\hat S_n}^\prime.
    \]

    Let
\[
    E_\delta
    =
    \left\{
        \widetilde H^2(\hat S_n,S_0)
        >
        1154(c_0+1)^2(\delta^2+h_n^2)
    \right\},
    \qquad
    D(S)
    =
    \widetilde H^2\left(
        \frac{S+S_n^*}{2},S_n^*
    \right).
\]
By the triangle inequality and
Lemma~\ref{lem: hellinger of convex combination},
\[
    \widetilde H^2(\hat S_n,S_0)
    \le32r_n^2+2h_n^2.
\]
Consequently, on $E_\delta$ we have $r_n^2>\delta^2/16$.
Moreover, the preceding case analysis shows that the approximation
term cannot be the largest term on $E_\delta$.
Hence, on this event, at least one of
\[
    (P_n-P)g_{\hat S_n}\ge r_n^2/3,
    \qquad
    (P_n-P)g_{\hat S_n}^\prime\ge r_n^2/3
\]
holds.

Since $D(S)\le F\nu(A)$ for every $S\in\Fcal_n$,
the event $E_\delta$ is empty if $\delta^2\ge16F\nu(A)$.
Otherwise, set
\[
    T
    =
    \min\left\{
        t\in\{0,1,2,\ldots\}:
        2^{2t-2}\delta^2\ge F\nu(A)
    \right\}.
\]
Partitioning the range $D(S)>\delta^2/16$ into the intervals
$(2^{2t-4}\delta^2,2^{2t-2}\delta^2]$, $t=0,\ldots,T$,
and using $2^{2t-4}/3\ge2^{2t-6}$, we obtain
\begin{align}
    \Pbb(E_\delta)
    &\le
    \sum_{t=0}^T
    \Pbb\left(
        \sup_{\substack{
            S\in\Fcal_n:\\
            D(S)\le2^{2t-2}\delta^2
        }}
        \sqrt n(P_n-P)g_S
        \ge2^{2t-6}\sqrt n\,\delta^2
    \right)\notag\\
    &\quad+
    \sum_{t=0}^T
    \Pbb\left(
        \sup_{\substack{
            S\in\Fcal_n:\\
            D(S)\le2^{2t-2}\delta^2
        }}
        \sqrt n(P_n-P)g_S^\prime
        \ge2^{2t-6}\sqrt n\,\delta^2
    \right).
    \label{eq: oracle-peeling}
\end{align}

    We first control the terms involving \(g_S\). For each \(t=0,\ldots,T\), set
\[
    a_t
    \coloneqq
    2^{2t-6}\sqrt n\,\delta^2 .
\]
This quantity will be used as the parameter \(a\) in
Theorem~\ref{thm: uniform bernstein inequality}. We verify the conditions of the theorem with
\[
    a=a_t,
    \qquad
    K=1,
    \qquad
    R=2^{t+2}C_\rho\delta.
\]
Let
\[
    \Gcal_{n,t}
    \coloneqq
    \left\{
        g_S:
        S\in\Fcal_n,\ 
        \widetilde H^2\left(
            \frac{S+S_n^*}{2},
            S_n^*
        \right)
        \leq
        2^{2t-2}\delta^2
    \right\}.
\]

We first check the Bernstein-size condition. By
Lemma~\ref{lem: hellinger upper bound for bernstein semi-norm}, for every
\(g_S\in\Gcal_{n,t}\),
\[
    \rho_1^2(g_S)
    \leq
    2^6C_\rho^2
    \widetilde H^2\left(
        \frac{S+S_n^*}{2},
        S_n^*
    \right)
    \leq
    2^{2t+4}C_\rho^2\delta^2.
\]
Hence, we have
\[
    \sup_{g\in\Gcal_{n,t}}\rho_1(g)
    \leq
    2^{t+2}C_\rho\delta
    =
    R,
\]
which verifies condition~\eqref{eq: bernstein condition}.

Next we check condition~\eqref{eq: condition 1 for bernstein}. Since
\[
    \frac{\sqrt n R^2}{K}
    =
    \sqrt n
    \left(
        2^{t+2}C_\rho\delta
    \right)^2
    =
    2^{2t+4}C_\rho^2\sqrt n\,\delta^2,
\]
we have
\[
    a_t
    =
    2^{2t-6}\sqrt n\,\delta^2
    \leq
    C_1
    2^{2t+4}C_\rho^2\sqrt n\,\delta^2
    =
    C_1\frac{\sqrt n R^2}{K},
\]
provided \(C_1C_\rho^2\geq 2^{-10}\). Since \(C_\rho^2\geq 1\), this is satisfied, for instance, by taking \(C_1=15\). Thus
condition~\eqref{eq: condition 1 for bernstein} holds.

We next verify condition~\eqref{eq: condition 2 for bernstein}. We need
    $a_t
    \leq
    8\sqrt n R.$
Substituting the definitions of \(a_t\) and \(R\), this is equivalent to
\[
    2^{t-11}\delta
    \leq
    C_\rho.
\]
It is enough to show that
\[
    2^{2t-4}\delta^2
    \leq
    2^{18}C_\rho^2.
\]
Since \(t\leq T\), by the definition of \(T\), we have
\[
    2^{2t-4}\delta^2
    \leq
    2^{2(T-1)-2}\delta^2
    \le
    F\nu(A).
\]
Moreover, by enlarging \(c_0\) if necessary, we may assume
    $C_\rho^2\geq F\nu(A)\lor 1$.
Therefore condition~\eqref{eq: condition 2 for bernstein} holds after changing only universal constants.

It remains to verify the entropy condition~\eqref{eq:condition 3 for bernstein}. By
Lemma~\ref{lem: relation of bracketing between bernstein semi-norm and hellinger distance}, if
\[
    \left\|
        \sqrt{\frac{S_U+S_n^*}{2}}
        -
        \sqrt{\frac{S_L+S_n^*}{2}}
    \right\|_{L^2(\nu\otimes P_Z)}
    \leq u,
\]
then
\[
    \rho_1(g_U-g_L)
    \leq
    C_\rho u.
\]
Consequently, for every \(u>0\),
\[
    \Ncal_{B,1}
    \left(
        C_\rho u,
        \Gcal_{n,t},
        P
    \right)
    \leq
    N_{2,B}
    \left(
        u,
        \bar{\Fcal}_n^{1/2}(2^{t-1}\delta),
        \nu\otimes P_Z
    \right).
\]
Equivalently,
\[
    \log
    \Ncal_{B,1}
    \left(
        v,
        \Gcal_{n,t},
        P
    \right)
    \leq
    \log
    N_{2,B}
    \left(
        \frac{v}{C_\rho},
        \bar{\Fcal}_n^{1/2}(2^{t-1}\delta),
        \nu\otimes P_Z
    \right).
\]

We now bound the entropy integral in
condition~\eqref{eq:condition 3 for bernstein}. Since
    $a_t/(2^6\sqrt n)
    =
    2^{2t-12}\delta^2$
and
    $R=2^{t+2}C_\rho\delta$,
we have
\begin{align*}
    &\int_{a_t/(2^6\sqrt n)}^{R}
    \sqrt{
        \log
        \Ncal_{B,1}
        \left(
            u,
            \Gcal_{n,t},
            P
        \right)
    }\,du \\
    &\leq
    C_\rho
    \int_{2^{2t-12}\delta^2/C_\rho}^{2^{t+2}\delta}
    \sqrt{
        \log
        N_{2,B}
        \left(
            v,
            \bar{\Fcal}_n^{1/2}(2^{t-1}\delta),
            \nu\otimes P_Z
        \right)
    }\,dv,
\end{align*}
where we used the change of variables \(v=u/C_\rho\).
By the definition of \(J_B\) and the assumption
    $\Psi(s)\geq J_B(s,\bar{\Fcal}_n^{1/2}(s))$,
we have
\begin{align*}
    &C_\rho
    \int_{2^{2t-12}\delta^2/C_\rho}^{2^{t+2}\delta}
    \sqrt{
        \log
        N_{2,B}
        \left(
            v,
            \bar{\Fcal}_n^{1/2}(2^{t-1}\delta),
            \nu\otimes P_Z
        \right)
    }\,dv, \\
    &\le
    C_\rho
    \int_{2^{2t-12}\delta^2/C_\rho}^{2^{t+2}\delta}
    \sqrt{
        \log
        N_{2,B}
        \left(
            v,
            \bar{\Fcal}_n^{1/2}(2^{t+2}\delta),
            \nu\otimes P_Z
        \right)
    }\,dv,
    \le C_\rho\Psi(2^{t+2}\delta).
\end{align*}
Since \(\Psi(s)/s^2\) is nonincreasing and \(\delta\geq\delta_n\), the critical inequality
    $\sqrt n\,\delta_n^2
    \geq
    c\Psi(\delta_n)$
implies
    $\sqrt n\,\delta^2
    \geq
    c\Psi(\delta).$
Moreover, again using the monotonicity of \(\Psi(s)/s^2\), we have
    $\Psi(2^{t+2}\delta)
    \leq
    2^{2t+4}\Psi(\delta).$
Hence,
\[
    C_\rho\Psi(2^{t+2}\delta)
    \leq
    \frac{C_\rho}{c}
    2^{2t+4}\sqrt n\,\delta^2.
\]
Then, we obtain
\[
    \int_{a_t/(2^6\sqrt n)}^{R}
    \sqrt{
        \log
        \Ncal_{B,1}
        \left(
            u,
            \Gcal_{n,t},
            P
        \right)
    }\,du
    \leq
    \frac{2^{10}C_\rho}{c}a_t.
\]
Furthermore, by the same critical inequality and the fact that \(\Psi(\delta)\geq J_B(\delta,\bar{\Fcal}_n^{1/2}(\delta)) \geq \delta\), enlarging \(c\) if necessary gives
\[
    R
    =
    2^{t+2}C_\rho\delta
    \leq
    \frac{2^{10}C_\rho}{c}a_t.
\]
Therefore, with $C_0 = c / (2^{10}C_\rho)$,
\[
    a_t
    \geq
    C_0
    \left(
        \int_{a_t/(2^6\sqrt n)}^{R}
        \sqrt{
            \log
            \Ncal_{B,1}
            \left(
                u,
                \Gcal_{n,t},
                P
            \right)
        }\,du
        \lor R
    \right),
\]
which verifies condition~\eqref{eq:condition 3 for bernstein}.

Finally, condition~\eqref{eq: condition 4 for bernstein} is satisfied by choosing
\(C = C_0/4\) sufficiently large so that
\[
    C_0^2\geq C^2(C_1+1).
\]
Thus all the conditions of Theorem~\ref{thm: uniform bernstein inequality} are verified for the class \(\Gcal_{n,t}\) with
\[
    a=a_t,
    \qquad
    K=1,
    \qquad
    R=2^{t+2}C_\rho\delta.
\]
Therefore,
\[
    \Pbb\left(
        \sqrt n
        \sup_{g\in\Gcal_{n,t}}
        |(P_n-P)g|
        \geq
        a_t
    \right)
    \leq
    C\exp\left(
        -
        \frac{a_t^2}{C^2(C_1+1)R^2}
    \right).
\]
Since
\[
    \frac{a_t^2}{R^2}
    =
    \frac{
        2^{4t-12}n\delta^4
    }{
        2^{2t+4}C_\rho^2\delta^2
    }
    =
    \frac{2^{2t-16}n\delta^2}{C_\rho^2},
\]
we obtain
\[
    \Pbb\left(
        \sup_{g\in\Gcal_{n,t}}
        (P_n-P)g
        \geq
        2^{2t-6}\delta^2
    \right)
    \leq
    C\exp\left(
        -
        \frac{2^{2t}n\delta^2}{2^{20}C^2C_\rho^2}
    \right),
\]
where the universal constants have been enlarged if necessary.
Summing over \(t\) gives
\[
    \sum_{t=0}^T
    \Pbb\left(
        \sup_{g\in\Gcal_{n,t}}
        (P_n-P)g
        \geq
        2^{2t-6}\delta^2
    \right)
    \leq
    \sum_{t=0}^\infty
    \Pbb\left(
        \sup_{g\in\Gcal_{n,t}}
        (P_n-P)g
        \geq
        2^{2t-6}\delta^2
    \right)
    \leq
    c\exp\left(
        -\frac{n\delta^2}{c^2}
    \right).
\]

    We next control the terms involving \(g_S^\prime\). For each \(t=0,\ldots,T\), we use the same value
\[
    a_t
    =
    2^{2t-6}\sqrt n\,\delta^2
\]
as the parameter \(a\) in Theorem~\ref{thm: uniform bernstein inequality}. Define
\[
    \Gcal_{n,t}^\prime
    \coloneqq
    \left\{
        g_S^\prime:
        S\in\Fcal_n,\ 
        \widetilde H^2\left(
            \frac{S+S_n^*}{2},
            S_n^*
        \right)
        \leq
        2^{2t-2}\delta^2
    \right\}.
\]
We verify the conditions of Theorem~\ref{thm: uniform bernstein inequality} with
\[
    a=a_t,
    \qquad
    K=F\nu(A),
    \qquad
    R=2^t\sqrt{F\nu(A)}\,\delta.
\]

First, by Lemma~\ref{lem: Hellinger upper bound of Bernstein size 2}, for every
\(g_S^\prime\in\Gcal_{n,t}^\prime\),
\[
    \rho_{F\nu(A)}^2(g_S^\prime)
    \leq
    4F\nu(A)
    \widetilde H^2\left(
        \frac{S+S_n^*}{2},
        S_n^*
    \right)
    \leq
    F\nu(A)2^{2t}\delta^2
    = R^2
\]
This verifies condition~\eqref{eq: bernstein condition}.

Next, we verify condition~\eqref{eq: condition 1 for bernstein}. Since
\[
    \frac{\sqrt nR^2}{K}
    =
    \frac{
        \sqrt n
        \left(
            2^t\sqrt{F\nu(A)}\,\delta
        \right)^2
    }{
        F\nu(A)
    }
    =
    2^{2t}\sqrt n\,\delta^2,
\]
we have
\[
    a_t
    =
    2^{2t-6}\sqrt n\,\delta^2
    \leq
    C_1 2^{2t}\sqrt n\,\delta^2
    =
    C_1\frac{\sqrt nR^2}{K}
\]
as long as \(C_1\geq 2^{-6}\). Thus condition~\eqref{eq: condition 1 for bernstein} holds, for example with \(C_1=15\).

We now check condition~\eqref{eq: condition 2 for bernstein}. We need
\[
    2^{t-9}\delta
    \leq
    \sqrt{F\nu(A)}.
\]
It is enough to show
\[
    2^{2t-18}\delta^2
    \leq
    F\nu(A).
\]
By the definition of \(T\), for \(t\leq T\), the quantity
\(2^{2t-4}\delta^2\) is bounded by \(F\nu(A)\).
Therefore condition~\eqref{eq: condition 2 for bernstein} follows.

We next verify the entropy condition~\eqref{eq:condition 3 for bernstein}. By
Lemma~\ref{lem: relation of bracketing between bernstein size and hellinger distance 2}, if
\[
    \left\|
        \sqrt{\frac{S_U+S_n^*}{2}}
        -
        \sqrt{\frac{S_L+S_n^*}{2}}
    \right\|_{L^2(\nu\otimes P_Z)}
    \leq u,
\]
then
\[
    \rho_{F\nu(A)}(g_U^\prime-g_L^\prime)
    \leq
    \sqrt{2F\nu(A)}\,u.
\]
Hence, for every \(u>0\),
\[
    \Ncal_{B,F\nu(A)}
    \left(
        \sqrt{2F\nu(A)}\,u,
        \Gcal_{n,t}^\prime,
        P
    \right)
    \leq
    N_{2,B}
    \left(
        u,
        \bar{\Fcal}_n^{1/2}(2^{t-1}\delta),
        \nu\otimes P_Z
    \right).
\]
Equivalently,
\[
    \log
    \Ncal_{B,F\nu(A)}
    \left(
        v,
        \Gcal_{n,t}^\prime,
        P
    \right)
    \leq
    \log
    N_{2,B}
    \left(
        \frac{v}{\sqrt{2F\nu(A)}},
        \bar{\Fcal}_n^{1/2}(2^{t-1}\delta),
        \nu\otimes P_Z
    \right).
\]

Since
    $a_t/(2^6\sqrt n)
    =
    2^{2t-12}\delta^2$
and
    $R=2^t\sqrt{F\nu(A)}\,\delta$,
we obtain, by the change of variables
\(v=u/\sqrt{2F\nu(A)}\),
\begin{align*}
    &\int_{a_t/(2^6\sqrt n)}^{R}
    \sqrt{
        \log
        \Ncal_{B,F\nu(A)}
        \left(
            u,
            \Gcal_{n,t}^\prime,
            P
        \right)
    }\,du \\
    &\leq
    \sqrt{2F\nu(A)}
    \int_{2^{2t-12}\delta^2/\sqrt{2F\nu(A)}}^{2^{t-1/2}\delta}
    \sqrt{
        \log
        N_{2,B}
        \left(
            v,
            \bar{\Fcal}_n^{1/2}(2^{t-1}\delta),
            \nu\otimes P_Z
        \right)
    }\,dv \\
    &\leq
    C_\rho
    \int_{2^{2t-12}\delta^2/C_\rho}^{2^{t-1/2}\delta}
    \sqrt{
        \log
        N_{2,B}
        \left(
            v,
            \bar{\Fcal}_n^{1/2}(2^{t-1}\delta),
            \nu\otimes P_Z
        \right)
    }\,dv
    .
\end{align*}
Here, we used that we may take $C_\rho \ge \sqrt{2F\nu(A)}$ by setting $c_0$ sufficiently large if necessary.
Remember the definition of $J_B$ and the condition $\Psi(s) \ge J_B(s, \bar{\Fcal}_n(s))$, we have
\begin{align*}
    &C_\rho
    \int_{2^{2t-12}\delta^2/C_\rho}^{2^{t-1/2}\delta}
    \sqrt{
        \log
        N_{2,B}
        \left(
            v,
            \bar{\Fcal}_n^{1/2}(2^{t-1}\delta),
            \nu\otimes P_Z
        \right)
    }\,dv, \\
    &\le
    C_\rho
    \int_{2^{2t-12}\delta^2/C_\rho}^{2^{t+2}\delta}
    \sqrt{
        \log
        N_{2,B}
        \left(
            v,
            \bar{\Fcal}_n^{1/2}(2^{t+2}\delta),
            \nu\otimes P_Z
        \right)
    }\,dv,
    \le C_\rho\Psi(2^{t+2}\delta).
\end{align*}
We note that, by setting $c_0$ sufficiently large if necessary, $2^{t+2}C_\rho \ge 2^t\sqrt{F\nu(A)}$.
Therefore, by repeating almost the same argument as controlling the terms involving $g_S$, we have
\[
    \sum_{t=0}^T
    \Pbb\left(
        \sup_{g\in\Gcal_{n,t}}
        (P_n-P)g
        \geq
        2^{2t-6}\delta^2
    \right)
    \leq
    \sum_{t=0}^\infty
    \Pbb\left(
        \sup_{g\in\Gcal_{n,t}}
        (P_n-P)g
        \geq
        2^{2t-6}\delta^2
    \right)
    \leq
    c\exp\left(
        -\frac{n\delta^2}{c^2}
    \right).
\]

    Combining the two bounds in \eqref{eq: oracle-peeling}, we obtain
    \[
        \Pbb\left(
            \widetilde H^2(\hat S_n,S_0)
            >
            1154(c_0+1)^2
            \left(
                \delta^2+\widetilde H^2(S_n^*,S_0)
            \right)
        \right)
        \leq
        c\exp\left(
            -\frac{n\delta^2}{c^2}
        \right).
    \]
    This completes the proof.
\end{proof}

\begin{proof}{(Proof of Corollary~\ref{cor: oracle inequality})}
    For any nonnegative random variable $Y$, we have
    \[
    \Ebb[Y] = \int_0^\infty \Pbb(Y \geq t)\,dt.
    \]
    Using this identity together with Theorem~\ref{thm: oracle inequality}, we obtain
    \begin{align*}
        \Ebb\left[\widetilde{H}^2(\hat{S}_n, S_0)\right]
        &=
        \int_0^\infty
        \Pbb\left(
            \widetilde{H}^2(\hat{S}_n, S_0) \geq t
        \right) dt \\
        &=
        \int_0^{1154(c_0+1)^2(\delta^2+\widetilde{H}^2(S_n^*,S_0))}
        \Pbb\left(
            \widetilde{H}^2(\hat{S}_n, S_0) \geq t
        \right) dt \\
        &\quad
        +
        \int_{1154(c_0+1)^2(\delta^2+\widetilde{H}^2(S_n^*,S_0))}^{\infty}
        \Pbb\left(
            \widetilde{H}^2(\hat{S}_n, S_0) \geq t
        \right) dt \\
        &\leq
        1154(c_0+1)^2
        \left(
            \delta^2
            +
            \widetilde{H}^2(S_n^*,S_0)
        \right) \\
        &\quad
        +
        1154(c_0+1)^2
        \int_{\delta^2}^{\infty}
        \Pbb\left(
            \widetilde{H}^2(\hat{S}_n,S_0)
            \geq
            1154(c_0+1)^2
            \left(
                t+\widetilde{H}^2(S_n^*,S_0)
            \right)
        \right) dt \\
        &\leq
        1154(c_0+1)^2
        \left(
            \delta^2
            +
            \widetilde{H}^2(S_n^*,S_0)
            +
            c\int_0^\infty
            \exp\left(-\frac{nt}{c^2}\right)dt
        \right) \\
        &=
        1154(c_0+1)^2
        \left(
            \delta^2
            +
            \widetilde{H}^2(S_n^*,S_0)
            +
            \frac{c^3}{n}
        \right).
    \end{align*}
\end{proof}

\section{Proof of Theorem~\ref{thm: convergence rates with SVB}} \label{app: proof of thm svb}

Lemma~\ref{lem: integral of the reciprocal of an intensity function} gives an upper bound on the integral of the reciprocal of a conditional intensity function, which appears in the proof of Theorem~\ref{thm: convergence rates with SVB}.

\begin{lemma} \label{lem: integral of the reciprocal of an intensity function}
    Suppose that there exist $\alpha \geq 0$ and a constant $C>0$ such that, for a conditional intensity function $S$,
    \begin{equation*}
        \int_{[0, 1]^{d_z}} \int_A \ind \left(S(x \mid z) \leq t\right)\, d\nu(x)dP_Z(z)
        \leq
        Ct^\alpha,
        \quad \text{for all } t \in (0,\infty).
    \end{equation*}
    Let $H\in[0,1]$.
    Then
    \begin{align*}
        \int_{[0, 1]^{d_z}} \int_A
        \frac{1}{S(x \mid z)}
        \ind\left(S(x \mid z) \geq H\right)
        \, d\nu(x)dP_Z(z)
        \leq
        \begin{cases}
        (C \lor \nu(A)) \dfrac{H^{\alpha-1}}{1-\alpha}, & \text{if } \alpha \in [0,1), \\
        (C \lor 1)\left(\nu(A)-\log(H)\right), & \text{if } \alpha \geq 1.
        \end{cases}
    \end{align*}
\end{lemma}

\begin{proof}
    For any nonnegative function $h$ and random variable $Y$, we have
    \[
    \Ebb[h(Y)] = \int_0^\infty \Pbb_Y(h(Y)\geq u)\,du.
    \]
    Combining this identity with Tonelli's theorem gives
    \begin{align*}
        &\int_{[0, 1]^{d_z}} \int_A
        \frac{1}{S(x \mid z)}
        \ind\left(S(x \mid z)\geq H\right)
        \, d\nu(x)dP_Z(z) 
        =
        \int_A \int_0^\infty
        \Pbb_Z\left(
            \frac{1}{S(x \mid Z)}
            \ind\left(S(x \mid Z)\geq H\right)
            \geq u
        \right)\,du\,d\nu(x) \\
        &\leq
        \int_A \int_0^\infty
        \Pbb_Z\left(
            H \leq S(x \mid Z)\leq \frac{1}{u}
        \right)\,du\,d\nu(x) 
        =
        \int_A \int_0^{1/H}
        \Pbb_Z\left(
            H \leq S(x \mid Z)\leq \frac{1}{u}
        \right)\,du\,d\nu(x) \\
        &\leq
        \int_0^{1/H}
        \int_A
        \Pbb_Z\left(
            S(x \mid Z)\leq \frac{1}{u}
        \right)\,d\nu(x)\,du.
    \end{align*}

    If $\alpha=0$, then $\Pbb_Z(S(x \mid Z)\leq t)\leq 1$ for any $t\geq 0$, and hence
    \begin{equation*}
        \int_0^{1/H}
        \int_A
        \Pbb_Z\left(S(x \mid Z)\leq \frac{1}{u}\right)\,d\nu(x)\,du
        \leq
        \frac{\nu(A)}{H}.
    \end{equation*}
    If $0<\alpha<1$, then the assumption gives
    \begin{equation*}
        \int_0^{1/H}
        \int_A
        \Pbb_Z\left(S(x \mid Z)\leq \frac{1}{u}\right)\,d\nu(x)\,du
        \leq
        C\int_0^{1/H} u^{-\alpha}\,du
        =
        \frac{C H^{\alpha-1}}{1-\alpha}.
    \end{equation*}
    If $\alpha\geq 1$ and $C\lor 1\leq H^{-1}$, then using
    \[
    \int_A \Pbb_Z(S(x \mid Z)\leq t)\,d\nu(x) \leq Ct
    \]
    for any $t\in[0,1]$, we obtain
    \begin{align*}
        &\int_0^{1/H}
        \int_A
        \Pbb_Z\left(S(x \mid Z)\leq \frac{1}{u}\right)\,d\nu(x)\,du \\
        &\leq
        \int_0^{C\lor 1}
        \int_A
        \Pbb_Z\left(S(x \mid Z)\leq \frac{1}{u}\right)\,d\nu(x)\,du
        +
        \int_{C\lor 1}^{1/H}
        \int_A
        \Pbb_Z\left(S(x \mid Z)\leq \frac{1}{u}\right)\,d\nu(x)\,du \\
        &\leq
        \int_0^{C\lor 1} \nu(A)\,du
        +
        C\int_{C\lor 1}^{1/H} \frac{1}{u}\,du \\
        &=
        (C\lor 1)\nu(A)
        +
        C\left(-\log(H)-\log(C\lor 1)\right).
    \end{align*}
    On the other hand, if $\alpha\geq 1$ and $C\lor 1\geq H^{-1}$, then
    \begin{equation*}
        \int_0^{1/H}
        \int_A
        \Pbb_Z\left(S(x \mid Z)\leq \frac{1}{u}\right)\,d\nu(x)\,du
        \leq
        \int_0^{C\lor 1} \nu(A)\,du
        =
        (C\lor 1)\nu(A).
    \end{equation*}
    Combining these bounds and simplifying the resulting expression yields the claim.
\end{proof}

Now, we prove Theorem~\ref{thm: convergence rates with SVB}.
\begin{proof}{(Proof of Theorem~\ref{thm: convergence rates with SVB})}
    From the condition $\lVert S-S_0\rVert_\infty \leq C_1/M$, we have, for any $x\in A$ and $z\in [0, 1]^{d_z}$,
    \[
    S(x \mid z)
    \geq
    S_0(x \mid z)-\frac{C_1}{M}.
    \]
    Combining this with
    \[
    \inf_{z\in [0, 1]^{d_z}}\inf_{x\in A} S(x \mid z)\geq \frac{1}{M},
    \]
    we obtain, for any $x\in A$ and $z\in [0, 1]^{d_z}$,
    \begin{equation*}
        S(x \mid z)
        \geq
        \left(S_0(x \mid z)-\frac{C_1}{M}\right)\lor \frac{1}{M}
        \geq
        \frac{S_0(x \mid z)}{C_1+1}\lor \frac{1}{M}.
    \end{equation*}
    In the last inequality, we used the fact that if
    $S_0(x \mid z)\geq (C_1+1)/M$, then
    \begin{equation*}
        S_0(x \mid z)-\frac{C_1}{M}
        \geq
        S_0(x \mid z)
        \left(1-\frac{C_1}{C_1+1}\right)
        =
        \frac{S_0(x \mid z)}{C_1+1}.
    \end{equation*}
    Using this bound, we obtain
    \begin{align*}
        \left(
            \sqrt{S(x \mid z)}
            -
            \sqrt{S_0(x \mid z)}
        \right)^2
        &\leq
        \frac{
            \left(
                S(x \mid z)-S_0(x \mid z)
            \right)^2
        }{
            \left(
                \sqrt{S(x \mid z)}
                +
                \sqrt{S_0(x \mid z)}
            \right)^2
        } \\
        &\leq
        \frac{C_1^2}{M}
        \ind\left(
            S_0(x \mid z)
            \leq
            \frac{C_1+1}{M}
        \right) \\
        &\quad
        +
        \frac{C_1^2}{M^2}
        \cdot
        \frac{C_1+1}{(\sqrt{C_1+1}+1)^2}
        \cdot
        \frac{1}{S_0(x \mid z)}
        \ind\left(
            S_0(x \mid z)
            \geq
            \frac{C_1+1}{M}
        \right).
    \end{align*}
    Integrating both sides with respect to $\nu$ and taking expectation with respect to $Z$ gives
    \begin{align*}
        \widetilde{H}^2(S,S_0)
        &\leq
        \frac{C_1^2}{M}
        \int_{[0, 1]^{d_z}} \int_A
        \ind\left(
            S_0(x \mid z)
            \leq
            \frac{C_1+1}{M}
        \right)
        d\nu(x)dP_Z(z) \\
        &\quad
        +
        \frac{C_1^2}{M^2}
        \int_{[0, 1]^{d_z}} \int_A
        \frac{
            \ind\left(
                S_0(x \mid z)
                \geq
                (C_1+1)/M
            \right)
        }{
            S_0(x \mid z)
        }
        d\nu(x)dP_Z(z).
    \end{align*}
    For the first term, the assumption gives
    \begin{equation*}
        \frac{C_1^2}{M}
        \int_{[0, 1]^{d_z}} \int_A
        \ind\left(
            S_0(x \mid z)
            \leq
            \frac{C_1+1}{M}
        \right)
        d\nu(x)dP_Z(z)
        \leq
        \frac{C_1^2 C}{M}
        \left(\frac{C_1+1}{M}\right)^{\alpha\land 1}.
    \end{equation*}
    For the second term, applying Lemma~\ref{lem: integral of the reciprocal of an intensity function} with
    $H=(C_1+1)/M$, we obtain, for $\alpha<1$,
    \begin{equation*}
        \frac{C_1^2}{M^2}
        \int_{[0, 1]^{d_z}} \int_A
        \frac{
            \ind\left(
                S_0(x \mid z)
                \geq
                (C_1+1)/M
            \right)
        }{
            S_0(x \mid z)
        }
        d\nu(x)dP_Z(z)
        \leq
        \frac{(C_1+1)^2}{M^{1+\alpha}}
        \cdot
        \frac{C\lor \nu(A)}{1-\alpha},
    \end{equation*}
    while for $\alpha\geq 1$,
    \begin{equation*}
        \frac{C_1^2}{M^2}
        \int_{[0, 1]^{d_z}} \int_A
        \frac{
            \ind\left(
                S_0(x \mid z)
                \geq
                (C_1+1)/M
            \right)
        }{
            S_0(x \mid z)
        }
        d\nu(x)dP_Z(z)
        \leq
        \frac{C_1^2}{M^2}
        (C\lor 1)
        \left(\nu(A)+\log(M)\right).
    \end{equation*}
    Combining these bounds and simplifying the expression yields the conclusion.
\end{proof}

\section{Proof of Theorem~\ref{thm: convergence rates for composition structured functions}}

In this section, we provide the proof of Theorem~\ref{thm: convergence rates for composition structured functions}.
We first state several definitions and lemmas used to evaluate the approximation error $\widetilde{H}^2(S_n^*, S_0)$ and the bracketing entropy integral in \eqref{eq: bracketing entropy integral}, which appear in Theorem~\ref{thm: oracle inequality} and Corollary~\ref{cor: oracle inequality}.

\begin{theorem}[Theorem 2 of \citealp{KohlerLanger21}] 
\label{thm: approximation error for Hölder class}
    Let $d \in \Nbb$, and let $f \in \Ccal^\beta([0,1]^d,K)$ for some $\beta,K \in \Rbb$.
    Let $M \in \Nbb$ be sufficiently large, and suppose that
    \begin{equation*}
        M \geq 2
        \quad \text{and} \quad
        M^{2\beta}
        \geq
        c_1(1\lor K)^{4(\lfloor\beta\rfloor+1)}
    \end{equation*}
    hold for some sufficiently large constant $c_1 \geq 1$.

    Let $L,r \in \Nbb$ and $B \in \Rbb$ satisfy
    \begin{align*}
        &\mathrm{(i)}\quad
        L
        \geq
        5
        +
        \left\lceil \log_4(M^{2\beta}) \right\rceil
        \left(
            \left\lceil \log_2(\lfloor\beta\rfloor \lor d) \right\rceil
            +1
        \right), \\
        &\mathrm{(ii)}\quad
        r
        \geq
        64\cdot 2^d
        \binom{d+\lfloor\beta\rfloor}{d}
        d^2(\lfloor\beta\rfloor+1)M^d, \\
        &\mathrm{(iii)}\quad
        B
        \ge
        M^{4\beta+4}\lor 2Ke^d \lor 4.
    \end{align*}
    Then there exists a neural network $\hat{f}\in\Fcal_{\mathrm{id}}(L,r,B)$ such that
    \begin{equation*}
        \lVert f-\hat{f}\rVert_\infty
        \leq
        c_2(1\lor K)^{4(\lfloor\beta\rfloor+1)}M^{-2\beta},
    \end{equation*}
    Here, $\mathrm{id}$ denotes the identity map, and $c_2$ is a constant depending only on $\beta, K$, and $d$.
\end{theorem}

\begin{remark}
    In the original statement of Theorem 2 of \cite{KohlerLanger21}, an upper bound on the network parameters is not explicitly stated.
    In this study, however, such a bound is needed to evaluate the covering number of the DNN class.
    Therefore, following the proof of Theorem 2 of \cite{KohlerLanger21}, we calculate and include the parameter bound $B$.
    For simplicity, we also take the original parameter $a$ to be $a=1$.
\end{remark}

\begin{lemma}[Lemma 3 of \citealp{Schmidt-Hieber20}] 
\label{lem: sup norm for composition of Hölder functions}
    Let $f=g_q\circ\cdots\circ g_0$, and suppose that, for each
    $g_i=(g_{ij})_{j=1,\dots,d^{(i+1)}}$, we have
    $g_{ij}\in \Ccal^{\beta_i}([a_i,b_i]^{t_i},K_i)$.
    Assume further that $K_i\geq 1$.
    For $i=1,\dots,q-1$, define
    \begin{equation*}
        h_0(x)
        \coloneqq
        \frac{g_0(x)}{2K_0}+\frac{1}{2},
        \qquad
        h_i(x)
        \coloneqq
        \frac{g_i(2K_{i-1}x-K_{i-1})}{2K_i}+\frac{1}{2},
        \qquad
        h_q(x)
        \coloneqq
        g_q(2K_{q-1}x-K_{q-1}).
    \end{equation*}
    Then
    \begin{equation*}
        f
        =
        g_q\circ\cdots\circ g_0
        =
        h_q\circ\cdots\circ h_0.
    \end{equation*}
    Moreover, $h_{0j}$ takes values in $[0,1]$, and
    \[
        h_{0j}\in \Ccal^{\beta_0}([0,1]^{t_0},1),
        \qquad
        h_{ij}\in \Ccal^{\beta_i}([0,1]^{t_i},(2K_{i-1})^{\beta_i}),
    \]
    for $i=1,\dots,q-1$, while
    \[
        h_{q1}\in
        \Ccal^{\beta_q}([0,1]^{t_q},K_q(2K_{q-1})^{\beta_q}).
    \]
    Furthermore, for any
    $\tilde{h}_i=(\tilde{h}_{ij})_{j=1,\dots,d^{(i+1)}}$ such that
    $\tilde{h}_{ij}:[0,1]^{t_i}\to[0,1]$ for $i=0,\dots,q-1$, we have
    \begin{equation*}
        \left\lVert
            h_q\circ\cdots\circ h_0
            -
            \tilde{h}_q\circ\cdots\circ\tilde{h}_0
        \right\rVert_\infty
        \leq
        K_q
        \prod_{l=0}^{q-1}(2K_l)^{\beta_{l+1}}
        \sum_{i=0}^q
        \left\lVert h_i-\tilde{h}_i \right\rVert_\infty^{
            \prod_{k=i+1}^q(\beta_k\land 1)
        }.
    \end{equation*}
\end{lemma}

\begin{Def}[Covering number with respect to the uniform norm]
    Let $\Fcal$ be a class of functions, and let $\delta>0$.
    We define $N(\delta,\Fcal,\lVert\cdot\rVert_\infty)$ as the smallest integer $N$ for which there exists a collection of functions $\{f_j\}_{j=1}^N$ such that
    \begin{equation*}
        \sup_{f\in\Fcal}
        \min_{j=1,\dots,N}
        \lVert f-f_j\rVert_\infty
        \leq
        \delta.
    \end{equation*}
    If no such $N$ exists, we set
    $N(\delta,\Fcal,\lVert\cdot\rVert_\infty)=\infty$.
    We call $N(\delta,\Fcal,\lVert\cdot\rVert_\infty)$ the $\delta$-covering number of $\Fcal$ with respect to the uniform norm.
\end{Def}

\begin{lemma}[Lemma A.7 of \citealp{YaraTerada2026}] 
\label{lem: relation between bracketing number and covering number}
    Let $\Fcal$ be a class of functions on $\Xcal$, and let $\mu$ be a finite measure on $\Xcal$.
    Then, for any $\delta>0$,
    \begin{equation*}
        N_{p,B}(\delta,\Fcal,\mu)
        \leq
        N\left(
            \frac{\delta}{2\mu(\Xcal)^{1/p}},
            \Fcal,
            \lVert\cdot\rVert_\infty
        \right).
    \end{equation*}
\end{lemma}

\begin{lemma}[Lemma 3 of \citealp{Suzuki19}] 
\label{lem: covering number of DNN}
    For any $\delta>0$,
    \begin{equation*}
        \log N\left(
            \delta,
            \mathcal{F}_{\mathrm{id}}(L-1,r,B,\infty,s),
            \lVert\cdot\rVert_\infty
        \right)
        \leq
        2sL\log((B\lor 1)(r+1))
        +
        s\log(\delta^{-1}L).
    \end{equation*}
    In particular, for fully connected (non-sparse) neural networks,
    \begin{equation*}
        \log N\left(
            \delta,
            \mathcal{F}_{\mathrm{id}}(L-1,r,B,\infty),
            \lVert\cdot\rVert_\infty
        \right)
        \leq
        2r^2L^2\log((B\lor 1)(r+1))
        +
        r^2L\log(\delta^{-1}L),
    \end{equation*}
    where $\mathrm{id}$ denotes the identity map.
\end{lemma}

\begin{remark}
    Note that the statement of Lemma~\ref{lem: covering number of DNN} differs slightly from that in \cite{Suzuki19}, because the convention for counting the number of layers in our definition of DNNs differs by one from that used in \cite{Suzuki19}.
\end{remark}

We are now ready to derive the convergence rate of the DNN-based maximum likelihood estimator when the true intensity function belongs to
$\Gcal_\alpha(q,\mathsf{t},\beta,K,C)$.

\begin{proof}{(Proof of Theorem\ref{thm: convergence rates for composition structured functions})}
    For each $i = 0, \dots, q$, let
    \[
        M_i
        =
        \left\lceil
        n^{
            \frac{1}{
                2\left((1+(\alpha\land 1))\beta_i^*+t_i\right)
            }
        }
        \right\rceil,
        \qquad
        \xi_n
        =
        \max_{i=0,\dots,q}
        n^{
            -
            \frac{
                \beta_i^*
            }{
                (1+(\alpha\land 1))\beta_i^*+t_i
            }
        }
        \qquad
        \phi_n
        =
        \max_{i=0,\dots,q}
        n^{
            -
            \frac{
                (1+(\alpha\land 1))\beta_i^*
            }{
                (1+(\alpha\land 1))\beta_i^*+t_i
            }
        }.
    \]
    We write the hypothesis space as
        $\Fcal_n \coloneqq \Fcal_\sigma(L,r,B,F)$.

    We first evaluate the approximation error
    $\widetilde{H}^2(S_n^*,S_0)$.
    Since
        $S_0 \in \Gcal(q,\mathsf{d},\mathsf{t},\beta,K)$,
    we can write
    \[
        S_0 = g_q\circ \cdots \circ g_0,
        \qquad
        g_{ij}\in \Ccal^{\beta_i}([a_i,b_i]^{t_i},K_i),
    \]
    for $i=0,\dots,q$ and $j=1,\dots,d^{(i+1)}$.
    Without loss of generality, we may assume that $K_i\geq 1$.
    By the first part of Lemma~\ref{lem: sup norm for composition of Hölder functions}, we can rewrite
    \[
        S_0 = h_q\circ\cdots\circ h_0,
    \]
    where
    \[
        h_{0j}\in \Ccal^{\beta_0}([0,1]^{t_0},1),
        \qquad
        h_{ij}\in \Ccal^{\beta_i}([0,1]^{t_i},(2K_{i-1})^{\beta_i}),
        \quad i=1,\dots,q-1,
    \]
    and
    \[
        h_{q1}\in
        \Ccal^{\beta_q}([0,1]^{t_q},K_q(2K_{q-1})^{\beta_q}).
    \]
    Hence, by Theorem~\ref{thm: approximation error for Hölder class}, for sufficiently large $n$, there exist DNNs
    $\tilde{h}_{ij}\in \Fcal_{\mathrm{id}}(L_i,r_i,B_i,\infty)$ such that
    \begin{align}
        \lVert \tilde{h}_{0j}-h_{0j}\rVert_\infty
        &\leq
        c_2 M_0^{-2\beta_0},
        \nonumber \\
        \lVert \tilde{h}_{ij}-h_{ij}\rVert_\infty
        &\leq
        c_2(2K_{i-1})^{4\beta_i(\lfloor\beta_i\rfloor+1)}
        M_i^{-2\beta_i},
        \quad i=1,\dots,q-1,
        \label{eq: approximation error of composition elements} \\
        \lVert \tilde{h}_{q1}-h_{q1}\rVert_\infty
        &\leq
        c_2
        \left(
            K_q(2K_{q-1})^{\beta_q}
        \right)^{4(\lfloor\beta_q\rfloor+1)}
        M_q^{-2\beta_q},
        \nonumber
    \end{align}
    where $\mathrm{id}$ denotes the identity map.
    Here,
    \begin{align*}
        L_i
        &=
        5
        +
        \left\lceil \log_4(M_i^{2\beta_i}) \right\rceil
        \left(
            \left\lceil \log_2(\lfloor\beta_i\rfloor \lor t_i) \right\rceil
            +1
        \right), \\
        r_i
        &=
        64\cdot 2^{t_i}
        \binom{t_i+\lfloor\beta_i\rfloor}{t_i}
        t_i^2(\lfloor\beta_i\rfloor+1)M_i^{t_i}, \\
        B_i
        &=
        M_i^{4\beta_i+4}\lor 2Ke^{t_i}.
    \end{align*}

    To apply the second part of
Lemma~\ref{lem: sup norm for composition of Hölder functions},
we restrict the range of $\tilde{h}_{ij}$, $i=0,\ldots,q-1$,
to $[0,1]$.
Define
\[
    h_{ij}^*
    =
    \sigma(\tilde{h}_{ij})
    -
    \sigma(\tilde{h}_{ij}-1).
\]
Then
\[
    h_{ij}^*
    =
    (\tilde{h}_{ij}\vee0)\wedge1,
    \qquad
    \sigma(h_{ij}^*)=h_{ij}^*.
\]
This clipping operation is represented by one additional hidden layer
with two active neurons, output weights $(1,-1)$, and biases $(0,1)$.
Since $r_i\ge2$ and $B_i\ge1$, zero padding gives
\[
    h_{ij}^*
    \in
    \Fcal_{\mathrm{id}}(L_i+1,r_i,B_i,\infty).
\]
Because $h_{ij}$ takes values in $[0,1]$,
\begin{equation}
    \lVert \sigma(h_{ij}^*)-h_{ij}\rVert_\infty
    \leq
    \lVert \tilde{h}_{ij}-h_{ij}\rVert_\infty.
    \label{eq: approximation error of clipped DNN}
\end{equation}
    Let
        $h_i^*=(h_{ij}^*)_{j=1,\dots,d^{(i+1)}}$.
    By computing the DNNs $h_{ij}^*$ in parallel, we obtain
        $h_i^*\in \Fcal_{\mathrm{id}}(L_i+1,d^{(i+1)}r_i,B_i,\infty)$.
    Define
    \[
        \widetilde{S}_n
        =
        \tilde{h}_{q1}
        \circ
        \sigma(h_{q-1}^*)
        \circ
        \cdots
        \circ
        \sigma(h_0^*).
    \]
    Then,
        $\widetilde{S}_n
        \in
        \Fcal_{\mathrm{id}}(\widetilde{L},\tilde{r},\widetilde{B},\infty)$,
    where
    \[
        \widetilde{L}
        \coloneqq
        2q+1+\sum_{i=0}^q L_i,
        \qquad
        \tilde{r}
        \coloneqq
        \left(
            \max_{i=0,\dots,q-1}d^{(i+1)}r_i
        \right)
        \lor r_q,
        \qquad
        \widetilde{B}
        \coloneqq
        \max_{i=0,\dots,q}B_i.
    \]
    By Lemma~\ref{lem: sup norm for composition of Hölder functions},
    \eqref{eq: approximation error of composition elements}, and
    \eqref{eq: approximation error of clipped DNN}, we obtain
    \begin{align}
        \lVert \widetilde{S}_n-S_0\rVert_\infty
        &\leq
        K_q
        \prod_{l=0}^{q-1}(2K_l)^{\beta_{l+1}}
        \left(
            \sum_{i=0}^{q-1}
            \lVert \sigma(h_i^*)-h_i\rVert_\infty^{
                \prod_{k=i+1}^q(\beta_k\land 1)
            }
            +
            \lVert \tilde{h}_{q1}-h_{q1}\rVert_\infty
        \right)
        \nonumber \\
        &\leq
        K_q
        \prod_{l=0}^{q-1}(2K_l)^{\beta_{l+1}}
        \left(
            \sum_{i=0}^{q-1}
            \lVert \tilde{h}_i-h_i\rVert_\infty^{
                \prod_{k=i+1}^q(\beta_k\land 1)
            }
            +
            \lVert \tilde{h}_{q1}-h_{q1}\rVert_\infty
        \right)
        \nonumber \\
        &\leq
        C^\prime
        \max_{i=0,\dots,q}
        M_i^{-2\beta_i\prod_{k=i+1}^q(\beta_k\land 1)}
        \leq
        C^\prime\xi_n.
        \label{eq: approximation error before clip}
    \end{align}
    Here, $C^\prime$ is a constant depending only on
    $q,t,\beta,K$, and $d^{(0)}, \dots, d^{(q)}$.

    For all sufficiently large $n$, we have
$0<\xi_n\le F$ and $\widetilde B\ge F\vee1$.
Define
\[
    S_n^*
    =
    \xi_n
    +
    \sigma(\widetilde S_n-\xi_n)
    -
    \sigma(\widetilde S_n-F)
    =
    (\widetilde S_n\vee\xi_n)\wedge F.
\]
The additive constant is represented by a hidden neuron with zero
incoming weights and bias $-1$, followed by output weight $\xi_n$.
Together with the two clipping neurons, this adds one hidden layer
and uses three active neurons.
Since $\tilde r\ge3$, zero padding preserves the prescribed width.
The resulting output lies in $[\xi_n,F]$, so the prescribed final
ReLU leaves it unchanged.
Consequently,
\[
    S_n^*
    \in
    \Fcal_\sigma(\widetilde L+1,\tilde r,\widetilde B,F).
\]
    From the definition of $L_i$, there exists a constant $E$ independent of $M_i$'s such that
    \[
        2q+2+\sum_{i=0}^q L_i
        =
        E
        +
        \sum_{i=0}^q
        \left(
            \left\lceil \log_2(t_i\lor \beta_i)\right\rceil
            +1
        \right)
        \left\lceil 2\beta_i\log_4(M_i)\right\rceil.
    \]
    Since $\lceil x\rceil < x+1$, for sufficiently large $n$, we have
    \[
        \widetilde{L}+1
        \leq
        \sum_{i=0}^q
        2\beta_i
        \left(
            \log_2(\beta_i\lor t_i)+2
        \right)
        \log_4(M_i).
    \]
    Therefore, by the choices of $L,r$, and $B$ in the statement of the theorem, we have
    $S_n^*\in \Fcal_n=\Fcal_\sigma(L,r,B,F)$.

    We now evaluate the approximation error.
    If
        $S_n^*(x\mid z)=\widetilde{S}_n(x\mid z) \land F$,
    then \eqref{eq: approximation error before clip} and $F \ge K$ gives
    \[
        \lvert S_n^*(x\mid z)-S_0(x\mid z)\rvert
        \leq
        C^\prime\xi_n.
    \]
    If
        $S_n^*(x\mid z)=\xi_n$,
    then $\widetilde{S}_n(x\mid z)\leq \xi_n$.
    Hence, by \eqref{eq: approximation error before clip} and the nonnegativity of functions in
    $\Gcal_\alpha(q,t,\beta,K,C)$, we obtain
    \[
        0\leq S_0(x\mid z)
        \leq
        (C^\prime + 1)\xi_n.
    \]
    Thus, in this case as well,
    \[
        \lvert S_n^*(x\mid z)-S_0(x\mid z)\rvert
        \leq
        (C^\prime+2)\xi_n.
    \]
    Consequently,
    \[
        \lVert S_n^*-S_0\rVert_\infty
        \leq
        (C^\prime+2)\xi_n.
    \]
    Applying Theorem~\ref{thm: convergence rates with SVB} with
    $C_1=C^\prime+2$, we obtain
    \begin{equation*}
        \widetilde{H}^2(S_n^*,S_0)
        \lesssim
        \phi_n\log(n).
    \end{equation*}
    More explicitly, the hidden constant depends only on
    $q,t,\beta,K,\alpha, C,F$, and $d^{(0)}, \dots, d^{(q)}$.

    Moreover, for any $x\in A$ and $z\in [0, 1]^{d_z}$,
    \begin{equation*}
        \frac{S_0(x\mid z)}{S_n^*(x\mid z)}
        \leq
        \frac{
            (C^\prime+2)\xi_n
        }{
            S_n^*(x\mid z)
        }
        +1
        \leq
        C^\prime+3.
    \end{equation*}
    Hence, Assumption~\ref{asm: weak absolute continuity} holds with
    \[
        c_0 \geq \sqrt{C^\prime+3}.
    \]

    We next evaluate the bracketing entropy integral
    \eqref{eq: bracketing entropy integral} in order to apply
    Corollary~\ref{cor: oracle inequality}.
    For the above $S_n^*\in\Fcal_n$, define
    \[
        \bar{\Fcal}_n^{1/2}
        \coloneqq
        \left\{
            \sqrt{\frac{S+S_n^*}{2}}
            :
            S\in \Fcal_n
        \right\}.
    \]
    Since
        $\bar{\Fcal}_n^{1/2}(\delta)
        \subset
        \bar{\Fcal}_n^{1/2}$,
    for any $u>0$,
    \[
        N_{2,B}\left(
            u,\bar{\Fcal}_n^{1/2}(\delta),\nu\otimes P_Z
        \right)
        \leq
        N_{2,B}\left(
            u,\bar{\Fcal}_n^{1/2},\nu\otimes P_Z
        \right).
    \]
    Let $\{(f_j^L,f_j^U)\}_{j=1}^N$ be a minimal
    $2\sqrt{2\xi_n}\,u$-bracketing set for $\Fcal_n$ with respect to the
    $L^2(\nu\otimes P_Z)$ norm.
    Since $S_n^*\geq \xi_n$, we have
    \begin{align*}
        &\Ebb_Z\left[
            \int_A
            \left(
                \sqrt{\frac{f_j^U(x\mid Z)+S_n^*(x\mid Z)}{2}}
                -
                \sqrt{\frac{f_j^L(x\mid Z)+S_n^*(x\mid Z)}{2}}
            \right)^2
            d\nu(x)
        \right] \\
        &\leq
        \Ebb_Z\left[
            \int_A
            \frac{
                \left(f_j^U(x\mid Z)-f_j^L(x\mid Z)\right)^2
            }{
                8S_n^*(x\mid Z)
            }
            d\nu(x)
        \right]
        \leq
        u^2.
    \end{align*}
    Therefore, using Lemma~\ref{lem: relation between bracketing number and covering number} together with the fact that ReLU is 1-Lipschitz, we obtain
    \begin{align*}
        N_{2,B}\left(
            u,\bar{\Fcal}_n^{1/2}(\delta),\nu\otimes P_Z
        \right)
        &\leq
        N_{2,B}\left(
            2\sqrt{2\xi_n}\,u,
            \Fcal_n,
            \nu\otimes P_Z
        \right) \\
        &\leq
        N\left(
            \frac{\sqrt{2\xi_n}}{\sqrt{\nu(A)}}u,
            \Fcal_n,
            \lVert\cdot\rVert_\infty
        \right) \\
        &\leq
        N\left(
            \frac{\sqrt{2\xi_n}}{\sqrt{\nu(A)}}u,
            \Fcal_{\mathrm{id}}(L,r,B,\infty),
            \lVert\cdot\rVert_\infty
        \right).
    \end{align*}
    By Lemma~\ref{lem: covering number of DNN},
    \begin{equation*}
        \log
        N_{2,B}\left(
            u,\bar{\Fcal}_n^{1/2}(\delta),\nu\otimes P_Z
        \right)
        \leq
        2s(L+1)
        \log\left(
            \frac{
                (L+1)(r+1)(B\lor 1)\sqrt{\nu(A)}
            }{
                \sqrt{2\xi_n}
            }
            u^{-1}
        \right),
    \end{equation*}
    where $s$ denotes the total number of parameters of the DNN.

    Set
    \[
        D
        \coloneqq
        \frac{
            (L+1)(r+1)(B\lor 1)\sqrt{\nu(A)}
        }{
            \sqrt{2\xi_n}
        }.
    \]
    Then
    \begin{align*}
        \int_{\delta^2/(2^{16}C_\rho)}^\delta
        \sqrt{
            \log
            N_{2,B}\left(
                u,\bar{\Fcal}_n^{1/2}(\delta),\nu\otimes P_Z
            \right)
        }\;du 
        \leq
        \sqrt{2s(L+1)}
        \int_0^\delta
        \sqrt{\log\left(\frac{D}{u}\right)}\,du.
    \end{align*}
    By the change of variables $x=\sqrt{\log(D/u)}$, we obtain
\begin{align*}
    \int_0^\delta
    \sqrt{\log\left(\frac{D}{u}\right)}\,du
    &=
    2D
    \int_{\sqrt{\log(D/\delta)}}^\infty
    x^2e^{-x^2}\,dx
    =
    \delta \sqrt{\log(D/\delta)}
    +
    D
    \int_{\sqrt{\log(D/\delta)}}^\infty
    e^{-x^2}\,dx.
\end{align*}
For the second term, let $Y$ be a centered normal random variable with variance $1/2$.
Then its density is $\pi^{-1/2}e^{-x^2}$, and therefore,
\begin{align*}
    \int_{\sqrt{\log(D/\delta)}}^\infty e^{-x^2}\,dx
    &=
    \frac{\sqrt{\pi}}{2}
    \Pbb\left(
        |Y| \geq \sqrt{\log(D/\delta)}
    \right).
\end{align*}
Using the standard Gaussian tail bound,
\[
    \Pbb\left(
        |Y| \geq a
    \right)
    \leq
    2e^{-a^2},
    \qquad a>0,
\]
with $a=\sqrt{\log(D/\delta)}$, we have
\begin{align*}
    D
    \int_{\sqrt{\log(D/\delta)}}^\infty e^{-x^2}\,dx
    &\leq
    D\frac{\sqrt{\pi}}{2}
    \cdot
    2\exp\left\{-\log(D/\delta)\right\}
    =
    \sqrt{\pi}\,\delta.
\end{align*}
Combining the above estimates yields
\begin{align*}
    \int_0^\delta
    \sqrt{\log\left(\frac{D}{u}\right)}\,du
    &\leq
    \delta
    \left(
        \sqrt{\log(D/\delta)}
        +
        \sqrt{\pi}
    \right).
\end{align*}
    Therefore,
    \begin{equation*}
        J_B(\delta,\bar{\Fcal}_n^{1/2}(\delta))
        \leq
        \delta\sqrt{2s(L+1)}
        \left(
            \sqrt{\log(D/\delta)}
            +
            \sqrt{\pi}
        \right).
    \end{equation*}

    In Corollary~\ref{cor: oracle inequality}, take
    \[
        \Psi(\delta)
        =
        \delta\sqrt{2s(L+1)}
        \left(
            \sqrt{\log(D/\delta)}
            +
            \sqrt{\pi}
        \right).
    \]
    To satisfy the critical inequality \eqref{eq: critical inequality}, set
    \begin{equation*}
        \delta_n
        =
        \frac{
            \sqrt{2s(L+1)}
            \left(
                \sqrt{\log(\sqrt{n}D)}
                +
                \sqrt{\pi}
            \right)
        }{
            \sqrt{n}
        }.
    \end{equation*}
    Then Corollary~\ref{cor: oracle inequality} yields
    \begin{align*}
        \Ebb\left[\widetilde{H}^2(\hat{S}_n,S_0)\right]
        &\leq
        1154(1+c_0)^2
        \left(
            \frac{
                2s(L+1)
                \left(
                    \sqrt{\log(\sqrt{n}D)}
                    +
                    \sqrt{\pi}
                \right)^2
            }{n}
            +
            \widetilde{H}^2(S_n^*,S_0)
            +
            \frac{1}{n}
        \right).
    \end{align*}
    Combining this bound with
    \[
        \widetilde{H}^2(S_n^*,S_0)\lesssim \phi_n\log(n)
    \]
    and using $s\lesssim r^2L$, $L\asymp\log(n)$, and
    $r\asymp \sqrt{n\phi_n}$, we obtain that there exists a constant
    $C^\prime>0$, depending only on
    $q,t,\beta,K,\alpha,C,c_0$, and $d^{(0)}, \dots, d^{(q)}$, such that
    \begin{equation*}
        \Ebb\left[\widetilde{H}^2(\hat{S}_n,S_0)\right]
        \leq
        C^\prime \phi_n\log(n)^3.
    \end{equation*}
\end{proof}

\section{Proof of Theorem~\ref{thm: minimax optimality}}

\begin{proof}
    Let $v$ and $V$ denote lower and upper bounds, respectively, of the density function $p_Z$ of $Z$ with respect to the Lebesgue measure.
    Moreover, we denote $\epsilon_n \coloneqq \max_{i=0, \dots, q} n^{-\beta_i^*/(\beta_i^* + t_i)}$
    By Theorem 2.7 of \cite{Tsybakov2008}, it suffices to construct functions
    \[
        S_0,\dots,S_M \in \Gcal(q,t,\beta,K)
    \]
    for some integer $M>1$ and some $\kappa>0$ such that
    \begin{equation*}
        \mathrm{(i)}\quad
        \widetilde{H}^2(S_j,S_k)\geq \kappa\epsilon_n
        \quad
        \text{for all } 0\leq j<k\leq M,
    \end{equation*}
    and
    \begin{equation*}
        \mathrm{(ii)}\quad
        \sum_{j=1}^M
        \mathrm{KL}(P_j\parallel P_0)
        \leq
        \frac{M\log M}{9}.
    \end{equation*}
    Indeed, if such functions exist, then there exists a positive constant
    $c=c(\kappa,\gamma)$ such that
    \begin{equation*}
        \inf_{\hat{S}_n}
        \sup_{S_0\in\Gcal(q,t,\beta,K)}
        \Ebb\left[
            \widetilde{H}^2(\hat{S}_n,S_0)
        \right]
        \geq
        c\epsilon_n.
    \end{equation*}
    Here, the infimum is taken over all estimators, the expectation is taken with respect to the observed data, and $P_j$ denotes the distribution of the data when the conditional intensity function is $S_j$.

    We now construct functions
    $S_0,\dots,S_M\in\Gcal(q,t,\beta,K)$
    satisfying $\mathrm{(i)}$ and $\mathrm{(ii)}$.
    Let
    \[
        i^*
        \in
        \argmin_{i=0,\dots,q}
        \frac{\beta_i^*}{\beta_i^*+t_i}.
    \]
    For notational convenience, write
        $\beta^* \coloneqq \beta_{i^*},
        \,
        \beta^{**} \coloneqq \beta_{i^*}^*,
        \,
        t^* \coloneqq t_{i^*}$.
    Choose a nonzero nonnegative function
$\xi\in C_c^\infty((0,1))$.
Multiplying it by a sufficiently small positive constant, we may
assume that $\xi\in\Ccal^{\beta^*}(\Rbb,1)$.
Then $\xi\in L^1(\Rbb)$, and
$0<\lVert\xi^B\rVert_1<\infty$ for the positive exponent $B$
defined below.
    Set
    $
        m_n
        \coloneqq
        \left\lfloor
            \rho n^{1/(\beta^{**}+t^*)}
        \right\rfloor,
        \,
        h_n
        \coloneqq
        1/m_n,
    $
    where the constant $\rho\geq 1$ is chosen so that
    \begin{equation*}
        n h_n^{\beta^{**}+t^*}
        \leq
        \frac{\log 2}{72V \lVert \xi^B\rVert_1^{t^*}}.
    \end{equation*}
    Here,
    \[
        B
        \coloneqq
        \prod_{l=i^*+1}^q(\beta_l\land 1).
    \]

    Let $\Ucal_n
        \coloneqq
        \{
            (u_1,\dots,u_{t^*})
            :
            u_i\in
            \{0,h_n,2h_n,\dots,(m_n-1)h_n\}
        \}$.
    For each
    $
        u=(u_1,\dots,u_{t^*})
        \in
        \Ucal_n,
    $
    define
    \begin{equation*}
        \psi_{u}(x_1,\dots,x_{t^*})
        \coloneqq
        h_n^{\beta^*}
        \prod_{j=1}^{t^*}
        \xi\left(
            \frac{x_j-u_j}{h_n}
        \right).
    \end{equation*}
    Since $\xi\in \Ccal^{\beta^*}(\Rbb,1)$, for any multi-index $\gamma$ satisfying $\lVert\gamma\rVert_1<\beta^*$,
    \begin{equation*}
        \lVert \partial^{\gamma}\psi_{u}\rVert_\infty
        =
        h_n^{\beta^*}
        \prod_{j=1}^{t^*}
        \frac{
            \lVert \xi^{(\gamma_j)}\rVert_\infty
        }{
            h_n^{\gamma_j}
        }
        \leq
        1.
    \end{equation*}
    Moreover, for any multi-index $\gamma$ with
    $\lVert\gamma\rVert_1=\lfloor\beta^*\rfloor$, the triangle inequality and the fact that
    $\xi\in \Ccal^{\beta^*}(\Rbb,1)$ imply
    \begin{equation*}
        h_n^{\beta^*-\lfloor\beta^*\rfloor}
        \frac{
            \left|
            \prod_{j=1}^{t^*}
            \xi^{(\gamma_j)}
            \left(
                \frac{x_j-u_j}{h_n}
            \right)
            -
            \prod_{j=1}^{t^*}
            \xi^{(\gamma_j)}
            \left(
                \frac{y_j-u_j}{h_n}
            \right)
            \right|
        }{
            \max_{i=1,\dots,t^*}
            |x_i-y_i|^{\beta^*-\lfloor\beta^*\rfloor}
        }
        \leq
        t^*.
    \end{equation*}
    Hence,
\[
    \psi_u
    \in
    \Ccal^{\beta^*}([0,1]^{t^*},C_{\beta^*,t^*}).
\]
Here $C_{\beta^*,t^*}>0$ is a fixed constant, independent of
$n$, $u$, and $w$.
    For
        $w=(w_{u})_{u\in\Ucal_n}
        \in
        \{0,1\}^{|\Ucal_n|}$,
    define
    \begin{equation*}
        \varphi_{w}
        \coloneqq
        \sum_{u\in\Ucal_n}
        w_{u}\psi_{u}.
    \end{equation*}
    By construction, for distinct
    $u,u^\prime\in\Ucal_n$, the functions
    $\psi_{u}$ and $\psi_{u^\prime}$ have disjoint supports.
    Therefore,
\[
    \varphi_w
    \in
    \Ccal^{\beta^*}([0,1]^{t^*},2C_{\beta^*,t^*}).
\]

    For $i<i^*$, define
    \[
        g_i(x)
        \coloneqq
        (x_1,\dots,x_{d_i\land d_{i+1}},0,\dots,0)^\top.
    \]
    For $i=i^*$, define
    \begin{equation*}
        g_{i^*,w}(x)
        \coloneqq
        \left(
            \varphi_{w}(x_1,\dots,x_{t^*}),
            0,\dots,0
        \right)^\top.
    \end{equation*}
    For $i>i^*$, define
    \[
        g_i(x)
        \coloneqq
        (x_1^{\beta_i\land 1},0,\dots,0)^\top.
    \]
    Since
        $t_i\leq \min\{d_0,\dots,d_{i-1}\}$,
    and since the functions $\psi_{u}$ have disjoint supports, noting that
    $\beta^{**}=\beta^*B$, we define
    \begin{align*}
        f_{w}(x)
        &\coloneqq
        g_q\circ\cdots\circ g_{i^*+1}
        \circ g_{i^*,w}
        \circ g_{i^*-1}
        \circ\cdots\circ g_0(x) \\
        &=
        \varphi_{w}^B(x_1,\dots,x_{t^*}) \\
        &=
        \sum_{u\in\Ucal_n}
        w_{u}\psi_{u}^B(x_1,\dots,x_{t^*}).
    \end{align*}
    Taking $K$ sufficiently large ensures that
        $f_{w}\in \Gcal(q,t,\beta,K)$.
    For any $u\in\Ucal_n$, we have
    \begin{align*}
        \lVert \psi_{u}^B\rVert_1
        &=
        h_n^{\beta^*B}
        \prod_{j=1}^{t^*}
        \int
        \xi^B\left(
            \frac{x_j-u_j}{h_n}
        \right)
        dx_j 
        =
        h_n^{\beta^{**}+t^*}
        \prod_{j=1}^{t^*}
        \int
        \xi^B(x_j)\,dx_j 
        =
        h_n^{\beta^{**}+t^*}
        \lVert \xi^B\rVert_1^{t^*}.
    \end{align*}

    Let
    \[
        \mathrm{Ham}(w,w^\prime)
        \coloneqq
        \sum_{u\in\Ucal_n}
        \ind(w_{u}\neq w_{u}^\prime)
    \]
    denote the Hamming distance.
    Since the supports of $\psi_{u}$ are pairwise disjoint,
    \begin{align*}
        \widetilde{H}^2(f_{w},f_{w^\prime})
        &=
        \Ebb_Z
        \left[
            \frac{1}{2}
            \int_A
            \left(
                \sqrt{f_{w}(x\mid Z)}
                -
                \sqrt{f_{w^\prime}(x\mid Z)}
            \right)^2
            dx
        \right] \\
        &=
        \Ebb_Z
        \left[
            \frac{1}{2}
            \int_A
            \left(
                \sqrt{
                    \sum_{u\in\Ucal_n}
                    w_{u}\psi_{u}^B(x\mid Z)
                }
                -
                \sqrt{
                    \sum_{u\in\Ucal_n}
                    w_{u}^\prime\psi_{u}^B(x\mid Z)
                }
            \right)^2
            dx
        \right] \\
        &=
        \Ebb_Z
        \left[
            \frac{1}{2}
            \int_A
            \sum_{u\in\Ucal_n}
            \left(
                \sqrt{
                    w_{u}\psi_{u}^B(x\mid Z)
                }
                -
                \sqrt{
                    w_{u}^\prime\psi_{u}^B(x\mid Z)
                }
            \right)^2
            dx
        \right] \\
        &=
        \Ebb_Z
        \left[
            \frac{1}{2}
            \int_A
            \sum_{u\in\Ucal_n}
            \ind(w_{u}\neq w_{u}^\prime)
            \psi_{u}^B(x\mid Z)
            dx
        \right] \\
        &\geq
        \frac{v}{2}
        \mathrm{Ham}(w,w^\prime)
        h_n^{\beta^{**}+t^*}
        \lVert \xi^B\rVert_1^{t^*}.
    \end{align*}
    Here, among the variables $x_1,\dots,x_{t^*}$, we have relabeled those corresponding to the support of the Poisson process as $x$ and those corresponding to the covariates as $Z$.
    By Lemma 2.9 of \cite{Tsybakov2008} and $m_n^{t^*}=|\Ucal_n|$,
    there exists a subset $\Wcal\subset \{0,1\}^{m_n^{t^*}}$
    containing the all-one vector such that $|\Wcal| \geq 2^{m_n^{t^*}/8}$,
    and, for any distinct $w,w^\prime\in\Wcal$,
    \[
        \mathrm{Ham}(w,w^\prime)
        \geq
        \frac{m_n^{t^*}}{8}.
    \]
    Therefore, taking
    \[
        \kappa
        =
        \frac{
            v^{1/2} \lVert \xi^B\rVert_1^{t^*/2}
        }{
            4\rho^{\beta^{**}/2}
        },
    \]
    we obtain, for any distinct
    $w,w^\prime\in\Wcal$,
    \begin{equation*}
        \widetilde{H}^2(f_{w},f_{w^\prime})
        \geq
        \frac{v}{16}
        h_n^{\beta^{**}}
        \lVert \xi^B\rVert_1^{t^*}
        \geq
        \kappa^2\epsilon_n.
    \end{equation*}

    We next bound the KL divergence for the data distributions induced by the intensity functions constructed above.
    Let $P_{w}$ denote the population distribution when the intensity function is $f_{w}$, and let
    $f_1, P_{1}$ denote the conditional intensity function and population distribution corresponding to the all-one vector $(1,\dots,1)^\top$, respectively.
    Using the disjointness of the supports of $\psi_{u}$, the KL-divergence formula for Poisson processes \cite[Corollary 5.2]{leskelä24}, and the definitions of $h_n$ and $\rho$, we obtain
    \begin{align*}
        \mathrm{KL}
        \left(
            P_{w}^{1:n}
            \parallel
            P_{1}^{1:n}
        \right)
        &=
        n\Ebb_Z
        \left[
            \int
            \left\{
                f_{w}(x\mid Z)
                \log\left(
                    \frac{f_{w}(x\mid Z)}{f_{1}(x\mid Z)}
                \right)
                +
                f_{1}(x\mid Z)-f_{w}(x\mid Z)
            \right\} dx
        \right] \\
        &\leq
        n\Ebb_Z
        \left[
            \int
            \frac{
                (f_{w}(x\mid Z)-f_{1}(x\mid Z))^2
            }{
                f_{1}(x\mid Z)
            } dx
        \right] \\
        &\leq
        nV
        m_n^{t^*}
        \| \psi_{u}^B \|_1 \\
        &\leq
        nV
        m_n^{t^*}
        h_n^{\beta^{**}+t^*}
        \lVert \xi^B\rVert_1^{t^*} \\
        &\leq
        \frac{\log |\Wcal|}{9}.
    \end{align*}
    Here, for a probability distribution $P$, the notation $P^{1:n}$ denotes the joint distribution of $n$ i.i.d.\ observations from $P$.
    This proves the claim.
\end{proof}

\section{Proof of Theorem~\ref{thm: convergence rates for manifold}}
\label{app:manifold-proof}

We first supply uniform local extensions, then construct the reference
network and apply the oracle inequality.  Constants in this appendix do
not depend on the member $S$ of the fixed H\"older class or on the
approximation index.  Throughout, $q=\lceil\beta\rceil-1$ and
$\theta=\beta-q\in(0,1]$, including the Lipschitz endpoint.

\subsection{Charts, covariate boundaries, and local extensions}

The required geometric choices exist for every fixed compact smooth
embedded manifold without boundary.  At a point of $\Mcal$, projection
onto its tangent space has invertible derivative on the tangent space.
The inverse function theorem gives a chart; shrinking its image gives an
open ball.  Compactness gives a finite cover by smaller chart domains
whose closures lie inside such charts.
For completeness, the partition can be chosen with smooth ambient
extensions.  Choose nonnegative ambient bumps $b_i$ supported in
neighborhoods whose intersections with $\Mcal$ lie in $V_i$, with
$\sum_i b_i>0$ on $\Mcal$.  On a neighborhood $O$ of $\Mcal$ where this
sum is positive, take an ambient smooth cutoff $\zeta$ supported
compactly in $O$, equal to one near $\Mcal$.  The functions
\[
 R_i=\zeta\,\frac{b_i}{\sum_j b_j}\quad\text{on }O,
 \qquad R_i=0\quad\text{outside }O
\]
are smooth, take values in $[0,1]$, and restrict to the required
$\rho_i$ on $\Mcal$.  Refining the finite cover if necessary ensures
$\operatorname{supp}\rho_i\Subset V_i$.  In the product domain we use
$\rho_i(x,z)=\rho_i(x)$.

\begin{lemma}[Uniform local Euclidean extensions]
\label{lem:manifold-extension}
Under Definition~4.2, there exist fixed affine
bijections $A_i:[0,1]^k\to Q_i$, where $Q_i$ are rectangular boxes, and
functions $f_{i,S}\in C^{q,\theta}([0,1]^k)$ such that
\[
 f_{i,S}\bigl(A_i^{-1}(P_i x,z)\bigr)=S(x,z)
 \quad\text{if }x\in\operatorname{supp}\rho_i,\ z\in[0,1]^{d_z},
 \qquad
 \max_i\|f_{i,S}\|_{C^{q,\theta}([0,1]^k)}\le C_{\rm ext}K.
\]
The constant is uniform over $S\in\Hcal^\beta_{\Mcal}(K)$.
\end{lemma}

\begin{proof}
Write $u_i(v,z)=S(P_i^{-1}v,z)$ and
$E_i=P_i(\operatorname{supp}\rho_i)\subset\Omega_i$.
The covariate faces prevent $E_i\times[0,1]^{d_z}$ from being an interior
subset of an open product chart, so we first extend in those coordinates.
Let $J_q=q+1$ and choose the unique coefficients $a_1,\ldots,a_{J_q}$
satisfying the Vandermonde system
\begin{equation}\label{eq:reflection-jets}
 \sum_{j=1}^{J_q}a_j(-j)^r=1,\qquad r=0,\ldots,q.
\end{equation}
Set $\delta=1/(2J_q)$.  In one covariate coordinate, with all remaining
coordinates denoted by $w$, use
\[
 (Tu)(w,t)=
 \begin{cases}
  \sum_{j=1}^{J_q}a_j u(w,-jt),&-\delta<t<0,\\
  u(w,t),&0\le t\le1,\\
  \sum_{j=1}^{J_q}a_j u(w,1-j(t-1)),&1<t<1+\delta.
 \end{cases}
\]
Every argument on the right belongs to $[0,1]$.  Equation
\eqref{eq:reflection-jets} matches derivatives of order $r$ in $t$ at
both endpoints.  It also matches every mixed derivative of total order
at most $q$, because differentiating the other coordinates leaves the
same identities.  On each of the three pieces, the chain rule bounds
all derivative suprema and highest-order H\"older seminorms by a
constant depending only on $q,\theta$ times the original norm.
For points in different pieces, split the line segment at the intervening
faces.  The derivatives agree at those faces, and
$\sum_{\ell=1}^N h_\ell^\theta\le
N^{1-\theta}(\sum_{\ell=1}^N h_\ell)^\theta$ bounds the resulting
seminorm for every $0<\theta\le1$.  At $\theta=1$ the right-hand
side is simply the sum of the segment lengths, so the same estimate
is Lipschitz and has no singular endpoint constant.  The domain is
convex, so the segment stays in the domain.
Applying this construction successively to the $d_z$ coordinates gives
$\bar u_i$ on
$W_i=\Omega_i\times(-\delta,1+\delta)^{d_z}$, with
\[
 \bar u_i=u_i\ \text{on }\Omega_i\times[0,1]^{d_z},
 \qquad \|\bar u_i\|_{C^{q,\theta}(W_i)}\le C_{\beta,d_z}K.
\]
The same segment argument across at most $2d_z$ faces controls corners.
For $\beta=1$, $q=0$ and $a_1=1$: this is ordinary reflection of a
Lipschitz function, and only continuity is matched at each face.
When $d_z=0$, this step is the identity.

Choose fixed smooth cutoffs $\chi_i\in C_c^\infty(\Omega_i)$ and
$\eta\in C_c^\infty((-\delta,1+\delta)^{d_z})$, equal to one near
$E_i$ and $[0,1]^{d_z}$, respectively; omit $\eta$ when $d_z=0$.
Define $F_{i,S}=\chi_i\eta\bar u_i$ on $W_i$ and zero outside.
Leibniz's rule and
$[fg]_\theta\le\|f\|_\infty[g]_\theta+
\|g\|_\infty[f]_\theta$ give the $C^{q,\theta}$ bound on $W_i$.
Lower-order derivatives are Lipschitz on this bounded convex domain,
with bounds supplied by the next derivatives, and hence have bounded
$\theta$-seminorms.  The product vanishes in a fixed neighborhood of
$\partial W_i$.  Its zero extension is consequently $C^{q,\theta}$:
for a nonzero derivative value inside and a point outside $W_i$, the
distance is at least the fixed distance of the cutoff support from the
boundary.  This bounds the cross-boundary seminorm also at $\theta=1$;
no derivative of order $q+1$ of the target is required.  Therefore
\[
 \|F_{i,S}\|_{C^{q,\theta}(\Rbb^k)}\le C_iK,
 \qquad F_{i,S}(P_i x,z)=S(x,z)
 \quad(x\in\operatorname{supp}\rho_i).
\]
Take a fixed box $Q_i$ containing $E_i\times[0,1]^{d_z}$ in its interior
and a fixed affine parametrization $A_i$ of that box.  Then
$f_{i,S}=F_{i,S}\circ A_i$ has the asserted properties by the chain rule.
Taking the maximum of the finitely many constants proves uniformity.
The possibly signed extension is an approximation target; clipping below
imposes the intensity bounds.  No general extension theorem is used.
\end{proof}

\subsection{A reference network with intrinsic width}

We use Theorem~\ref{thm: approximation error for Hölder class} in its
asymptotic form: for fixed dimension $d$, smoothness $s>0$, and H\"older
radius, an error $O(m^{-2s})$ is attainable with depth $O(\log m)$,
width $O(m^d)$, and parameters $O(m^{4s+4})$.  The constants are uniform
on that H\"older ball.  In fixed dimension, the maximum-based norm and
Euclidean distance in Definition~4.2 control the
sum-based norm and maximum distance used in that theorem up to a fixed
constant, which is absorbed into the H\"older radius.  At integer
smoothness this is exactly the $C^{q,1}$ class used in the imported
approximation result; no smoothing or loss of smoothness is needed.

We record the elementary network operations used below.  For fixed
$U>0$,
\[
 \operatorname{clip}_{[0,U]}(v)=\sigma(v)-\sigma(v-U),\qquad
 \operatorname{clip}_{[\varepsilon,U]}(v)
 =\varepsilon+\sigma(v-\varepsilon)-\sigma(v-U).
\]
An additive constant is represented by a hidden neuron with zero incoming
weights and a negative bias.  Thus the absence of an output bias in
(4.1) causes no restriction.  To compose subnetworks
without multiplying two dense weight matrices, propagate a scalar linear
output $v$ as the two nonnegative channels $\sigma(v),\sigma(-v)$ and
recover $v$ in the next affine map.  Parallel branches are padded using
identity maps on nonnegative hidden activations.  Each scalar connection adds a fixed number of layers and channels.
Padding equalizes branch depths at their maximum, and the finite number
of successive stages keeps the total depth of order $\log m$ below.
Widths and parameter bounds change by fixed factors only.

Multiplication on $[0,F]\times[0,1]$ can be implemented with error
$\varepsilon$, constant width, depth $O(\log(1/\varepsilon))$, and
parameters bounded by a constant depending only on $F$.
Here is an explicit construction, following the squaring argument of
\citet{Yarotsky2017}.  Put
\[
 T(t)=2\sigma(t)-4\sigma(t-1/2)+2\sigma(t-1),\qquad
 Q_R(t)=t-\sum_{\ell=1}^R4^{-\ell}T^{\circ\ell}(t).
\]
On $[0,1]$, $Q_R$ is the linear interpolant of $t^2$ on the dyadic mesh
of spacing $2^{-R}$, as follows by induction on its midpoint values.
Hence $\|Q_R-t^2\|_\infty\le4^{-R}$.
A network propagating $t,T^{\circ\ell}(t)$, and the accumulated partial
sum has constant width, $O(R)$ layers, and bounded parameters.  For
$a=h/F$ and $b=v$, define
\[
 \operatorname{Mult}_{F,\varepsilon}(h,v)
 =\frac F2\{4Q_R((a+b)/2)-Q_R(a)-Q_R(b)\}.
\]
The error from $hv$ is at most $3F4^{-R}$, so
$R\ge\max\{1,\lceil\log_4(3F/\varepsilon)\rceil\}$ suffices.
Signed channels as above implement its possibly signed output.

\begin{lemma}[Uniform manifold approximation]\label{lem:manifold-network}
For all sufficiently large integers $m$ and every
$S\in\Hcal^\beta_{\Mcal}(K)$ with $0\le S\le F$, there is a ReLU
network $S_m^*$ with
\[
 m^{-2\beta}\le S_m^*\le F\ \text{on }\Rbb^D,\qquad
 \|S_m^*-S\|_{\infty,K_0}\le C_A m^{-2\beta},
\]
and depth at most $C_L\log m$, width at most $C_r m^k$, and parameter
bound at most $C_Bm^{4\beta+4}$.  All constants are uniform over $S$.
\end{lemma}

\begin{proof}
Set $\varepsilon=m^{-2\beta}$.  Lemma~\ref{lem:manifold-extension}
and the Euclidean approximation theorem give networks $g_i$ approximating
$f_{i,S}$ to error $C\varepsilon$, with depth $O(\log m)$, width
$O(m^k)$, and parameter bound $O(m^{4\beta+4})$.  The map
$y=(x,z)\mapsto A_i^{-1}(P_i x,z)$ is affine on the ambient space.
Its input dimension is fixed, so absorbing it into the first affine
layer changes the parameter bound by a fixed factor.  Set
\[
 h_i(y)=\operatorname{clip}_{[0,F]}
          \{g_i(A_i^{-1}(P_i x,z))\}.
\]
Then $0\le h_i\le F$ everywhere and
$|h_i(y)-S(y)|\le C\varepsilon$ whenever
$x\in\operatorname{supp}\rho_i$.  No approximation outside that support
is needed.

Fix a noninteger $\gamma>\beta D/k$ and a bounded ambient box containing
$K_0$.  Approximate the fixed smooth ambient extensions $R_i(x)$ on this
box in dimension $D$, after a fixed affine rescaling, using
$m_\rho=\lceil m^{\beta/\gamma}\rceil$.  Clip to $[0,1]$ to obtain
$\widehat\rho_i$ with error $C\varepsilon$ on $K_0$ and budgets
\[
 L_\rho\lesssim\log m,\qquad
 r_\rho\lesssim m^{\beta D/\gamma}\le m^k,\qquad
 B_\rho\lesssim m^{4\beta+4\beta/\gamma}\le m^{4\beta+4}.
\]
The last inequality uses $\beta/\gamma<k/D\le1$.
Now form
\[
 \widetilde S_m(y)=\sum_{i=1}^J
 \operatorname{Mult}_{F,\varepsilon}(h_i(y),\widehat\rho_i(y)).
\]
Since the partition is applied exactly once, for $y\in K_0$,
\begin{align*}
 |\widetilde S_m(y)-S(y)|
 &\le J\varepsilon+
       \sum_i|h_i(y)|\,|\widehat\rho_i(y)-\rho_i(y)|
       +\sum_i\rho_i(y)|h_i(y)-S(y)|\\
 &\le C\varepsilon.
\end{align*}
The middle sum is controlled by global clipping.  In the last sum, a
nonzero partition weight guarantees that the local approximation holds.
Finite parallelization, the signed-channel connections, multiplication,
and summation have the claimed budgets.  Finally set
$S_m^*=\operatorname{clip}_{[\varepsilon,F]}(\widetilde S_m)$, for
$m$ sufficiently large that $\varepsilon\le F$.  For $S(y)\in[0,F]$,
clipping increases the approximation error by at most $\varepsilon$.
The final output is nonnegative, so applying the prescribed output ReLU
leaves it unchanged.  This proves all assertions, including the global
output bound and uniformity.
\end{proof}

\subsection{Main proof}

\begin{proof}[Proof of Theorem~4.4]
Set $\xi_n=n^{-\beta/(a\beta+k)}$ and take
$m_n=\lceil n^{1/[2(a\beta+k)]}\rceil$ in
Lemma~\ref{lem:manifold-network}.  For the final clipping, use $\xi_n$
in place of $m_n^{-2\beta}$.  Since
$m_n^{-2\beta}\le\xi_n\le2^{2\beta}m_n^{-2\beta}$, the resulting
network, denoted by $S_n^*$, satisfies
\[
 S_n^*\ge\xi_n,\qquad
 \|S_n^*-S_0\|_{\infty,K_0}\le C_A\xi_n.
\]
Its depth, width, and parameter bound are dominated by $L,r,B$ when
the implicit lower constants in (i)--(iii) are chosen sufficiently large.
Zero-neuron padding
and identity layers on nonnegative activations embed it into the exact
architecture without changing its output or parameter orders.  Thus
$S_n^*\in\Fcal_\sigma(L,r,B,F)$ uniformly over the asserted class.

\begin{samepage}
Let $\mu=\nu\otimes P_Z$ and use the notation $a,k,\phi_n'$ from
Theorem~4.4.  Since
$S_0/S_n^*\le 1+|S_0-S_n^*|/S_n^*\le1+C_A$,
Assumption~3.1 holds with a constant
$c_0$ independent of $n$.  Apply
Theorem~4.1 with
$M=\xi_n^{-1}$ and $C_1=\max\{1,C_A\}$.  Its lower-threshold condition holds for all
sufficiently large $n$, and it gives
\[
    \widetilde H^2(S_n^*,S_0)
    \lesssim \xi_n^{1+(\alpha\land1)}(1+\log(1/\xi_n))
    \lesssim \phi_n'\log n.
\]
\end{samepage}

For completeness, consider the entropy term.  Every network in
$\Fcal_n$ takes values in $[0,F]$.  Bracket endpoints may therefore be
clipped to this interval without increasing their width.  Since
$S_n^*\ge\xi_n$, for $u\ge v\ge0$,
\[
    \left|\sqrt{\frac{u+S_n^*}{2}}-
                  \sqrt{\frac{v+S_n^*}{2}}\right|
    \le\frac{|u-v|}{2\sqrt{2\xi_n}}.
\]
All covering norms below are taken on a fixed bounded ambient box containing
$K_0$.  An affine normalization to $[0,1]^D$ changes the input-layer
parameter bound by a fixed factor, so Lemma~\ref{lem: covering number of DNN}
applies with fixed-domain constants.  Combining its bound with
covering-to-bracketing conversion and the preceding inequality yields,
for $0<t\le1$,
\[
    \log N_{2,B}\left(t,\bar\Fcal_n^{1/2}(\delta),\mu\right)
    \le C s_n(L+1)\log(D_n/t),
\]
where the total number of network parameters is
$s_n=(L-1)r^2+(D+L+1)r$, and one may choose
\[
    D_n=e\lor
    \left\{C(L+1)(r+1)(B\lor1)
    \sqrt{\nu(\Mcal)/\xi_n}\right\}.
\]
Here and below constants may depend on the fixed input dimension and the
fixed domain.  The architecture assumptions imply
\[
    s_n\lesssim r^2L\lesssim n\phi_n'\log n,
    \qquad \log D_n\lesssim\log n.
\]
As in the entropy calculation for the composition theorem, an admissible
majorant is
\[
    \Psi_n(\delta)=C\delta\sqrt{s_n(L+1)}
        \left\{\sqrt{\log\!\left(\frac{D_n}{\delta\land1}\right)}+1\right\},
    \qquad\delta>0.
\]
For radii greater than one, extend the localized class by the same defining formula and use monotonicity of bracketing numbers in the
bracket width.  Integrating the bound with $t$ replaced by $t\land1$
from zero to $\delta$ verifies this majorant on the whole positive
axis.  Its ratio to $\delta^2$ is nonincreasing.  For a sufficiently large
fixed $C_0$, the radius
$\delta_n=C_0\{s_n(L+1)(\log n)/n\}^{1/2}$ satisfies the critical
inequality for all sufficiently large $n$.  In particular,
$\delta_n^2\lesssim\phi_n'\log^3 n$ and $\delta_n\to0$.
The constants and the threshold for $n$ are uniform in $S_0$.
Corollary~3.3 therefore yields
\[
 \sup_{S_0\in\Hcal^\beta_{\Mcal,\alpha}(K,F,C)}
 \Ebb[\widetilde H^2(\hat S_n,S_0)]
 \lesssim\delta_n^2+\phi_n'\log n+n^{-1}
 \lesssim\phi_n'\log^3 n.
\]
\end{proof}

\section{Supplementary Simulation Studies}
\label{sec:simulation-supplement}

\subsection{Generating intensities}
\label{subsec:simulation-generators}

We used Lebesgue measure on the Euclidean supports, arc length on the
unit circle, and surface area on the unit sphere. In the formulas below,
$c$ is a scaling constant chosen separately for each intensity, support,
and covariate dimension to give
$\mathbb E_Z\int S_0(x\mid Z)\,d\nu(x)\approx30$.
We calculated $c$ by numerical integration with resolution 1024 and
128 uniform covariate draws, using seed 314159. When $d_z=0$, the
covariate average reduces to a single integral. Covariate coordinates
absent from a design are set to zero.

For the compositional design, we set
\begin{equation}
 \begin{split}
 S_0(x\mid z)&=c\left[0.2+
 \{0.6+0.4\sin(2\pi h_1(x,z))\}^2
 +0.5\exp\{-10(h_2(x,z)-0.5)^2\}\right],\\
 h_1(x,z)&=x_1+0.3z_1+0.2z_2,\qquad
 h_2(x,z)=x_2+0.5\sin(2\pi z_1).
 \end{split}
 \label{eq:numerical-compositional-intensity}
\end{equation}
On $[0,1]$, we replace $x_2$ by $x_1/2$. Although the input dimension
varies, the intensity is formed by composing functions of a small number
of arguments. Only $z_1$ and $z_2$ enter the formula, leaving three
irrelevant covariates when $d_z=5$.

The second design is
\begin{equation}
 S_0(x\mid z)=c\left[10^{-4}+
 \{1-\exp(-18D^2(x,z))\}^2\right],
 \label{eq:simulation-near-zero}
\end{equation}
where
\[
 D^2(x,z)=
 \begin{cases}
 (x_1-0.25-0.5z_1)^2,&d_x=1,\\
 (x_1-0.25-0.5z_1)^2+
 \{x_2-0.5-0.25\sin(2\pi z_2)\}^2,&d_x=2.
 \end{cases}
\]
The intensity is small near $D(x,z)=0$, and the center of this region
moves with $z$. The additive constant gives a strictly positive minimum
of $10^{-4}c$ in every experiment.

On the unit circle, parameterized by angle $\theta$, we use
\begin{equation}
 S_0(\theta\mid z)=c\left[
 \{1+0.5\sin(3\theta+z_1)+0.3\cos(2\theta-z_2)\}_++10^{-4}\right].
 \label{eq:simulation-circle}
\end{equation}
On the unit sphere, with colatitude $\theta$ and longitude $\phi$, we use
\begin{equation}
 S_0(\theta,\phi\mid z)=c\left[
 \{1+0.5\sin(2\theta)\cos(3\phi+z_1)
       +0.3\sin(\phi+z_2)\}_++10^{-4}\right].
 \label{eq:simulation-sphere}
\end{equation}
The expressions inside the positive-part operators are bounded below by
$0.2$, so the operators do not alter these intensities. The term
$\sin(\phi+z_2)$ in \eqref{eq:simulation-sphere} has no common limit
at either pole. Any choice of values at the poles defines the same
intensity measure, but the resulting function is not globally continuous.
The sphere example therefore concerns estimation of a bounded measurable
intensity outside the smooth manifold class used in the rate theorem.

We simulated point patterns by thinning uniform candidate locations.
For each covariate vector, we evaluated the intensity on a grid with
resolution 1024 and set an envelope $B$ to 1.25 times the grid maximum,
with a small additive safeguard. We drew a Poisson number of candidates
with mean $B\nu(A)$ and retained a candidate at $x$ with probability
$\min\{S_0(x\mid z)/B,1\}$. This numerical envelope introduces a
possible simulation approximation if the intensity between grid points
exceeds $B$.

\subsection{DNN and kernel estimators}
\label{subsec:simulation-implementation}

For a candidate intensity $S$, the neural fitting criterion is
\[
 \ell_n(S)=\frac1n\sum_{i=1}^n
 \left\{\int S(x\mid Z_i)\,d\nu(x)
             -\sum_{x\in X_i}\log S(x\mid Z_i)\right\}.
\]
The two output functions are
$u\mapsto\log(1+e^u)+10^{-8}$ and $u\mapsto u_++10^{-8}$.
We initialized hidden layers by Kaiming initialization and output weights
near zero, with the output bias chosen to give approximately 30 expected
events on the support. AdamW was used with learning rate $10^{-3}$,
weight decay $10^{-4}$, and gradient-norm clipping at 5. Each update
used 64 patterns and 128 newly sampled integration points per pattern.
The epoch limits were 150 for Euclidean supports and 300 for manifold
supports. We stopped after 15 epochs without improvement in validation
loss and restored the best validation fit. The random 80\%/20\% split
was made at the pattern level; the networks were not subsequently refitted
on the validation patterns. Network training and prediction used single precision. We did not impose hard
parameter bounds or an upper bound on the fitted intensity.

All hidden layers have common width $W_n$, with depth and width given by
\begin{equation}
 L_n=\left\lceil 3\frac{\log n}{\log1000}\right\rceil,
 \qquad
 W_n=8\left\lceil
 \frac{64(n/1000)^{(1-\gamma_{\mathrm{arch}})/2}}8
 \right\rceil.
 \label{eq:simulation-architecture}
\end{equation}
For the compositional and moving-valley designs, respectively, we set
\[
 \gamma_{\mathrm{arch}}=\frac{2}{2+\min(d_x+d_z,4)},
 \qquad
 \gamma_{\mathrm{arch}}=\frac{4}{4+d_x+d_z}.
\]
For mesh-based integration on a manifold, we used
$\gamma_{\mathrm{arch}}=2/(2+d_{\mathrm{amb}}+d_z)$, where
$d_{\mathrm{amb}}=2$ for the circle and 3 for the sphere. For integration
with the known manifold measure, we replaced $d_{\mathrm{amb}}$ by
$d_{\mathcal M}=1$ or 2. In both versions, the networks take ambient
coordinates as inputs. These architecture choices were specified using
the generating design rather than selected from the data.

To construct an integration mesh, we combined 4,096 independent uniform
support samples with a reproducible subset of at most 4,096 event
locations. On the circle, convex-hull edges define a polygon; on the
sphere, triangular convex-hull faces define a polyhedral surface.
Integration weights are geometric edge lengths or planar face areas,
rather than event frequencies. The known-measure version instead uses
uniform samples on the true support and its volume. The comparison
between these versions consequently changes both the integration measure
and the network architecture.

For Euclidean supports, the conditional kernel estimator is
\[
 \widehat S(x\mid z)=
 \frac{\sum_{i=1}^n K_{h_z}(z-Z_i)
              \sum_{u\in X_i}K_{h_x}^{A}(x,u)}
      {\max\{\sum_{i=1}^n K_{h_z}(z-Z_i),10^{-6}\}},
\]
where $K_{h_z}$ is a Gaussian covariate kernel and $K_{h_x}^{A}$ is a
Gaussian event kernel normalized for the support boundary.
This estimator follows the conditional-kernel approach to Poisson
intensity estimation \citep{KlutchnikoffMassiot16}. The pooled estimator
uses equal weights for all patterns. We truncated Gaussian contributions
beyond four bandwidths. Five-fold Poisson cross-validation used pilot
bandwidths proportional to
$N_{\mathrm{events,train}}^{-1/(4+d_x)}$ for event coordinates and
$n_{\mathrm{train}}^{-1/(4+d_z)}$ for covariates, followed by a
five-value coarse search and a three-value refinement. To limit the
cost of bandwidth selection, this calculation used at most 512 patterns
and 128 events per pattern. The final kernel fit used all $n$ patterns.

On manifold supports, we used a geodesic kernel estimator motivated by
\citet{WardEtAl23}, with bandwidth selected by five-fold Poisson
cross-validation over a geodesic bandwidth grid. The estimator uses the
known geometry and does not include local or global shape corrections.

\subsection{Loss evaluation and Monte Carlo summaries}
\label{subsec:numerical-evaluation}

For each repetition, we generated one dataset and one set of evaluation
points shared by all methods, with separate random seeds for fitting.
The loss in (5.1) was averaged over 64 independent
uniform covariate vectors when $d_z>0$. For Euclidean supports, we used
512 integration points generated by scrambled Sobol quadrature when
available. On the circle, we used 64 arc-length-weighted points. On the
sphere, eight colatitude subdivisions and 16 longitude subdivisions
formed triangles; their projected centroids were the evaluation locations
and their planar areas the weights. Neural predictions were calculated
in single precision and converted to double precision for loss evaluation;
kernel predictions and loss accumulation used double precision.

The supplementary tables give the mean, median, standard deviation,
minimum, and maximum of the 20 repetitions for each method and
configuration. MCSEs are $\mathrm{SD}/\sqrt{20}$ and describe variation
across repetitions. They do not include systematic quadrature error;
the simulation losses were evaluated at a single integration resolution.

\begin{figure}[p]
 \centering
 \includegraphics[width=0.88\textwidth]{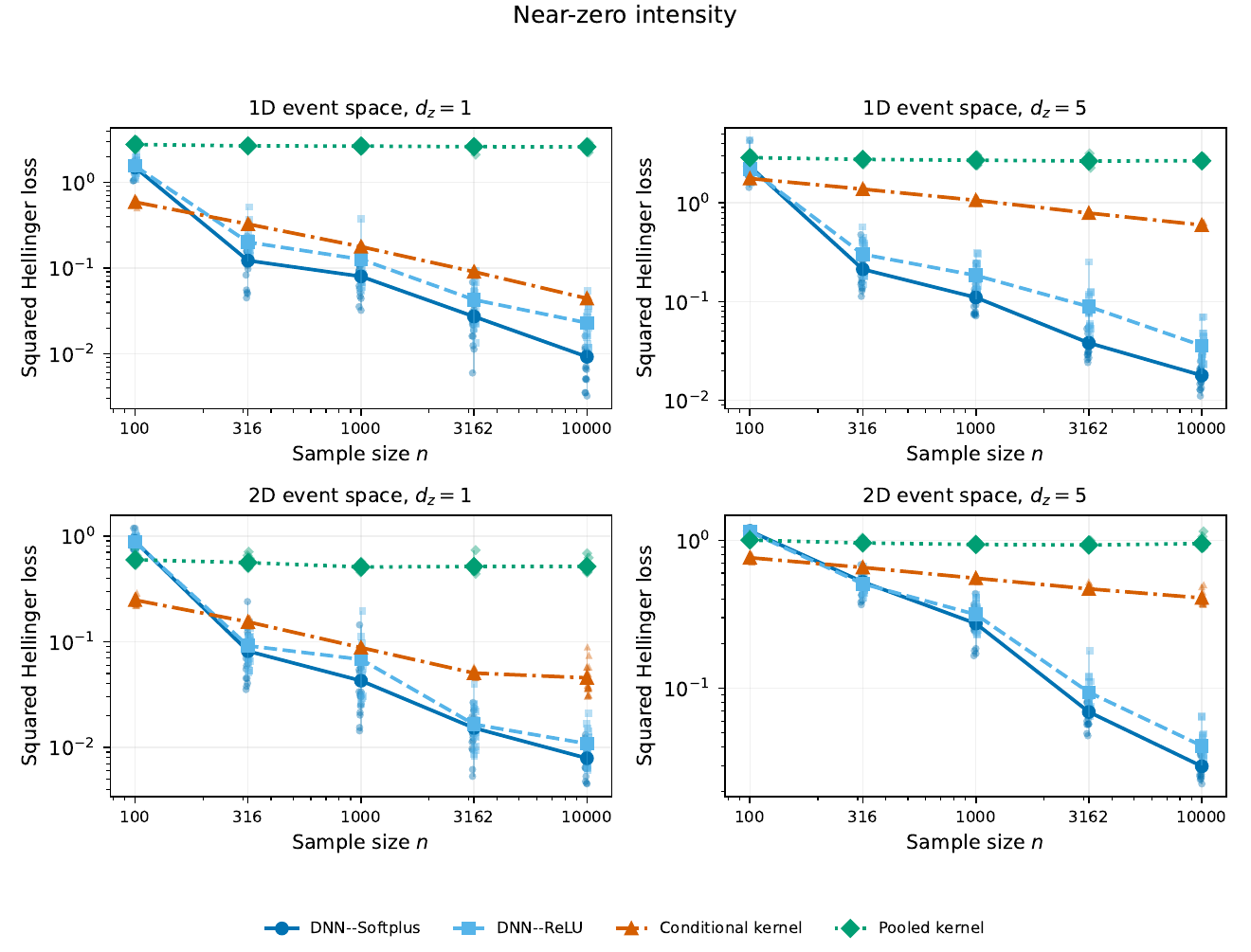}
 \caption{Squared Hellinger loss for the intensity with a moving valley.
 The intensity floor is $10^{-4}$ before rescaling. Points, connected
 markers, and vertical segments show individual repetitions, means,
 and ranges, as in Figure~1.}
 \label{fig:numerical-near-zero}
\end{figure}

\begin{figure}[p]
 \centering
 \includegraphics[width=0.88\textwidth]{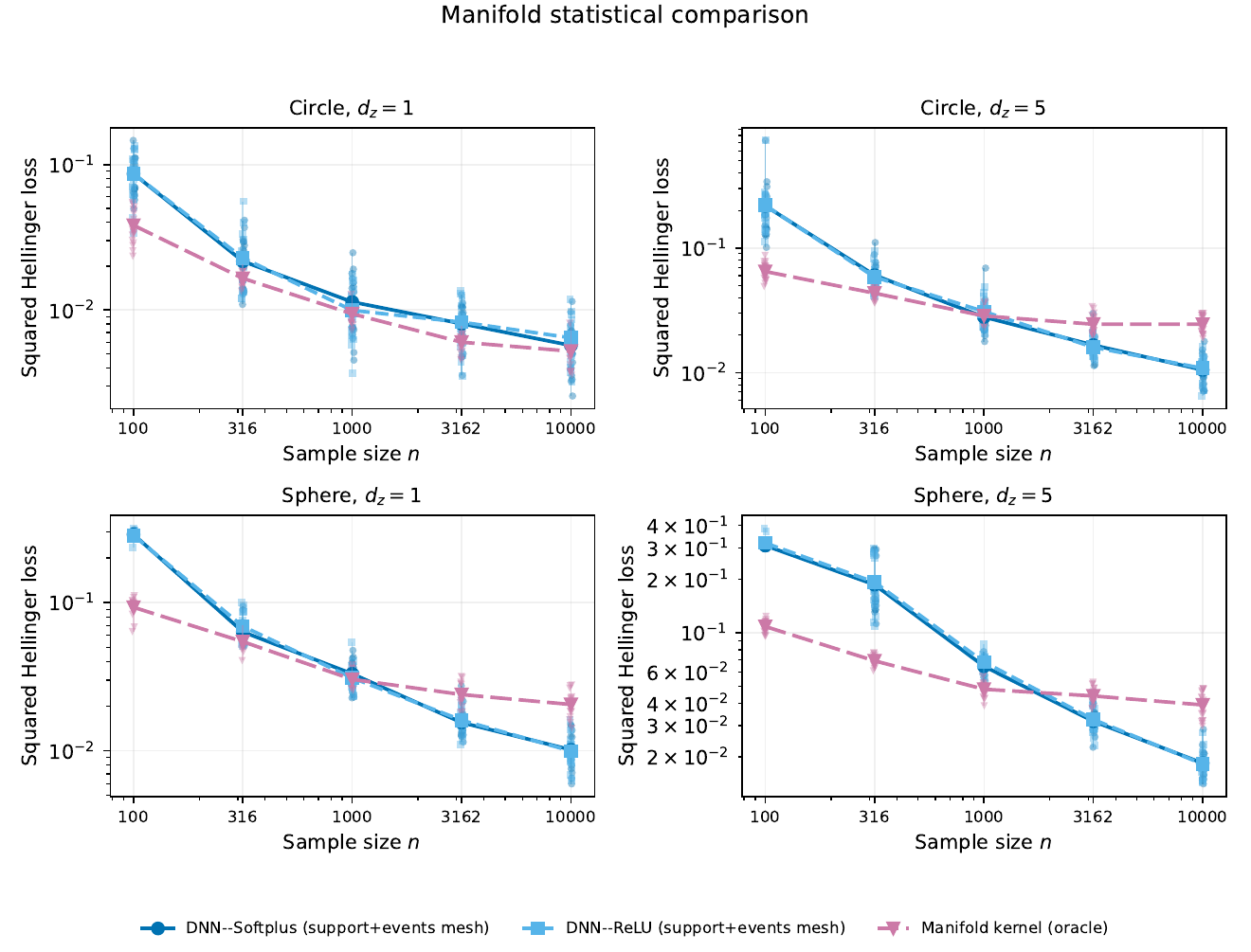}
 \caption{Squared Hellinger loss on the circle and sphere. The DNNs use
 ambient coordinates and reconstructed integration meshes; the kernel
 uses known geodesic distances and geometry. Graphical summaries follow
 Figure~1.}
 \label{fig:numerical-manifold-statistical}
\end{figure}

\clearpage
\subsection{No-covariate cases and manifold integration}
\label{subsec:simulation-additional-results}

\begin{figure}[p]
 \centering
 \includegraphics[width=\textwidth]{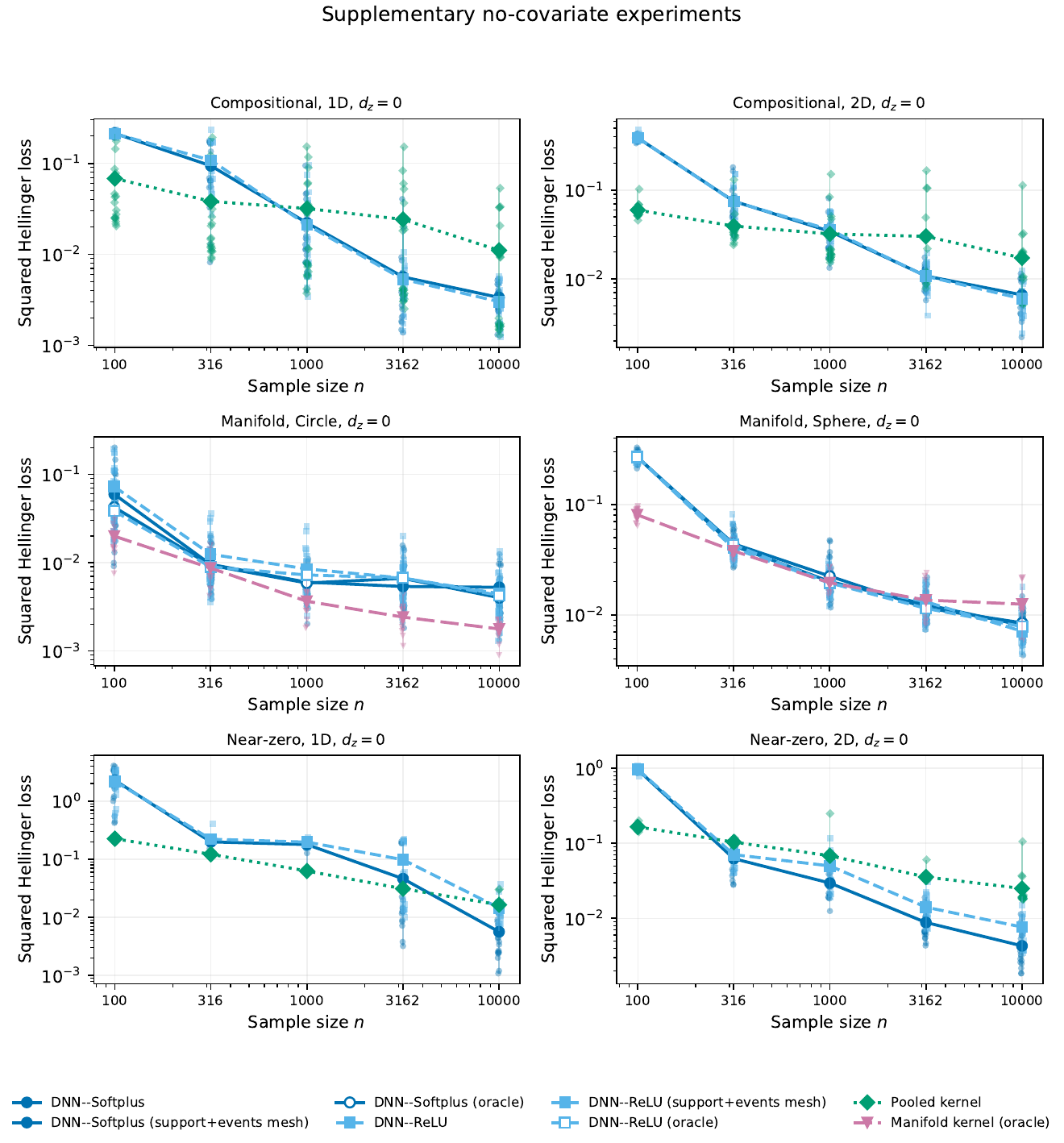}
 \caption{Squared Hellinger loss without covariates. The Euclidean kernel
 estimates the unconditional intensity. Points, connected markers, and
 vertical segments show individual repetitions, means, and ranges.}
 \label{fig:numerical-dz0}
\end{figure}

Figure~\ref{fig:numerical-dz0} shows the results without covariates.
These cases complement the main comparisons by separating estimation
of an ordinary intensity from estimation of its covariate dependence.
For Euclidean supports, the pooled and conditional kernel estimators
coincide, and we fitted only the pooled version.

\begin{figure}[p]
 \centering
 \includegraphics[width=\textwidth]{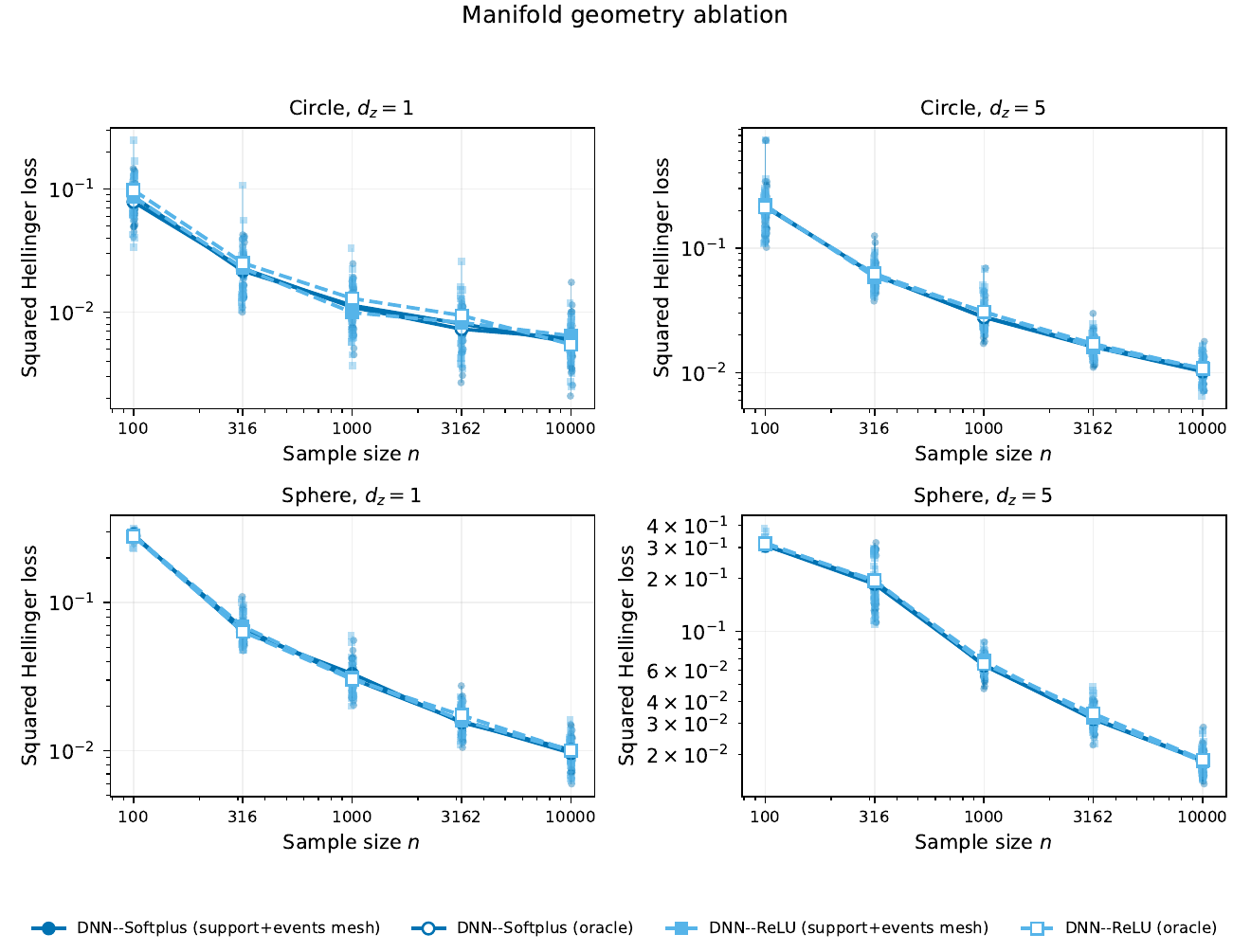}
 \caption{Neural estimators with reconstructed-mesh (M) and known-measure
 (O) integration. Both use ambient inputs, with different dimensions in
 the architecture rule \eqref{eq:simulation-architecture}.}
 \label{fig:numerical-manifold-ablation}
\end{figure}

Figure~\ref{fig:numerical-manifold-ablation} compares the mesh-based
and known-measure neural estimators. For DNN--Softplus at $n=10000$,
the respective mean losses on the circle were $0.00570$ and $0.00615$
with one covariate, and $0.0105$ and $0.0102$ with five covariates.
On the sphere, the corresponding pairs were $0.0102$ and $0.00967$,
and $0.0184$ and $0.0185$. Neither version uniformly dominated the
other. As described above, the architecture rule also differs between
the two versions, so these differences cannot be attributed to the
integration measure alone.


\subsection{Computation times}
\label{subsec:numerical-runtime}

\begin{figure}[p]
 \centering
 \includegraphics[width=0.95\textwidth]{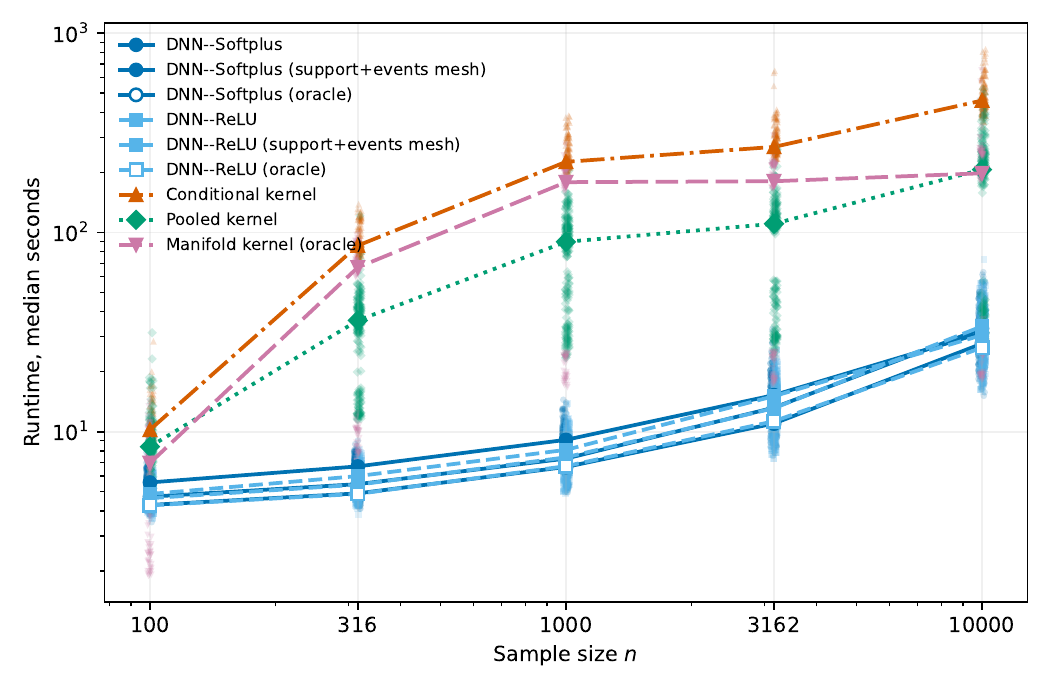}
 \caption{Median total computation time by sample size, pooled over the
 configurations and repetitions available for each method. Totals include
 data generation, fitting, evaluation, and output.}
 \label{fig:numerical-runtime}
\end{figure}

Figure~\ref{fig:numerical-runtime} and Table~\ref{tab:runtime-results}
summarize computation times. At $n=10000$, median total times ranged
from 26.4 to 33.8 seconds across the six DNN variants, compared with
198 seconds for the manifold kernel, 207 seconds for the pooled kernel,
and 459 seconds for the conditional kernel. Total time includes data
generation, fitting, evaluation, and file output. The configurations included and computing environments can differ
across methods, so these times are descriptive rather than a benchmark
under common hardware and workloads. Neural fitting used CUDA when
available, with CPU fallback otherwise.

\begin{table}[htbp]
\centering
\caption{Median computation times in seconds at $n=10000$, pooled over the configurations available for each method. M and O denote reconstructed-mesh and known-measure integration. Total time includes data generation, fitting, evaluation, and file output; component medians need not sum to the median total.}
\label{tab:runtime-results}
\small
\begin{tabular}{lrrrrr}
\toprule
Method & Runs & Total & Fit & Simulation & Evaluation \\
\midrule
DNN--Softplus & 240 & 31.3 & 26.4 & 2.03 & 0.714 \\
DNN--ReLU & 240 & 30.5 & 26.3 & 2 & 0.0455 \\
Conditional kernel & 160 & 459 & 198 & 0.751 & 252 \\
Pooled kernel & 240 & 207 & 68.9 & 0.762 & 135 \\
Manifold kernel O & 120 & 198 & 173 & 0.885 & 21.9 \\
DNN--Softplus M & 120 & 33 & 28.9 & 2.18 & 0.0408 \\
DNN--ReLU M & 120 & 33.8 & 29.3 & 2.19 & 0.0397 \\
DNN--Softplus O & 120 & 27.7 & 23.4 & 2.39 & 0.0409 \\
DNN--ReLU O & 120 & 26.4 & 22.1 & 2.34 & 0.0402 \\
\bottomrule
\end{tabular}
\end{table}

\section{Supplementary Neural Spike-Train Analysis}
\label{sec:allen-supplement}

\subsection{Data selection and preprocessing}

\begin{figure}[p]\centering
\includegraphics[width=\linewidth]{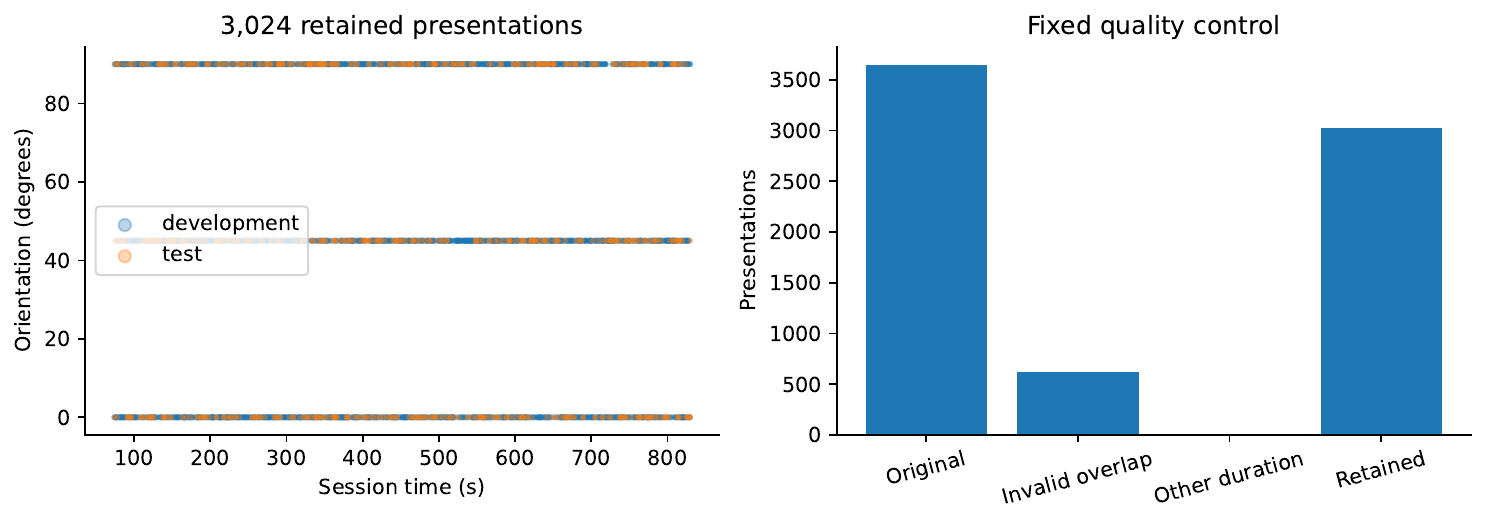}
\caption{Gabor presentation timing and orientation in the development and
test samples (left), and presentation counts before and after exclusions
(right). The split was stratified by position--orientation combination.}
\label{fig:allen-design-supp}
\end{figure}

We required an available NWB recording file and at least 20 VISp units
with $\texttt{isi\_violations}<0.5$,
$\texttt{amplitude\_cutoff}<0.1$, and
$\texttt{presence\_ratio}>0.9$, together with finite spike times.
Session 715093703 was the smallest eligible identifier, and 62 units
met these criteria. Here a unit is a spike-sorted recording unit, which
need not correspond to a perfectly isolated neuron
\citep{allensdkneuropixels}. We removed presentations with missing
stimulus covariates, overlap with an invalid recording interval, or a
duration different from the modal duration rounded to one millisecond.
The exclusions comprised 620 interval overlaps and one duration mismatch.
Figure~\ref{fig:allen-design-supp} summarizes presentation timing, orientations, and exclusions.

We assigned presentations to development and test samples using
within-condition deterministic hashes with seed 20260902. The test
sample contained the rounded 20\% of presentations in each condition,
subject to retaining observations in both samples when possible.
The 12-unit selection used development spike counts, with unit identifier
breaking ties. Unit 950931315, whose development count was nearest the
panel median, was used to illustrate fitted intensities and for the
single-unit supplementary analyses. Its earlier pilot test results had
already been viewed, as disclosed in the main text. That pilot used
normalization based on a reduced sample; all models reported here were
refitted with normalization based on the full development sample before
their test scores were calculated.

Horizontal and vertical positions were mapped to $[0,1]$ using the
development-sample extrema of $-40$ and $40$ degrees. We retained the
original $[-1,1]$ scale of the doubled-angle sine and cosine coordinates.
Both development and test samples contained all 243 combinations of
position and orientation.

\subsection{Fitting procedures}

The DNNs used AdamW with learning rate $10^{-3}$, weight decay
$10^{-4}$, batches of 64 presentations, and gradient-norm clipping at
5. Training integrals were approximated by stochastic quadrature with
256 points. We allowed at most 200 epochs and stopped after 20 epochs
without improvement in validation NLL. After selecting the epoch count
using the 1,935/482 fit/validation split, we refitted from initialization
on all 2,417 development presentations for that number of epochs.
Optimization and quadrature seeds were 20260903 and 20260904,
respectively. As in the simulations, no hard parameter bound or fixed
upper intensity bound was imposed.

The kernel estimators used Gaussian smoothing with boundary
renormalization on normalized time. Temporal bandwidths were selected
from $\{0.02,0.035,0.05,0.075,0.10,0.15,0.225\}$; for the
conditional estimator, covariate bandwidths were selected from
$\{0.10,0.15,0.225,0.35,0.50,0.75,1.00\}$. Five-fold
presentation-level Poisson cross-validation evaluated the full candidate
grid using all development presentations and spikes. The denominator
in the conditional estimator was bounded below by $10^{-6}$.

\subsection{Score aggregation and uncertainty}

Table~\ref{tab:allen-all-pairs} gives all pairwise comparisons, and
Table~\ref{tab:allen-unit-method} reports unit-specific scores and count
prediction errors. In addition to the equal-unit, or macro, information
gain used in the main text, we calculated a micro gain by pooling
log-likelihood improvements and spike counts across units before
converting to bits per spike. The 12-unit test sample contained 38,823
spikes. Relative to the conditional kernel, the micro gains were
$0.00785$ bits per spike for DNN--Softplus and $0.00563$ for DNN--ReLU.
Differences between macro and micro gains reflect the unequal numbers
of spikes across units.

For each of the 2,000 bootstrap draws, we sampled 243 stimulus-condition
labels with replacement, using seed 20260905. A presentation received
the multiplicity of its condition. We recomputed the presentation-weighted
mean score difference within each unit and then the equal-unit mean.
Using common multiplicities across methods and units preserves their
pairing. The percentile intervals hold the selected units and fitted
models fixed and describe uncertainty under condition resampling;
they do not account for refitting, between-session variation, or serial
dependence across conditions. Intervals are pointwise rather than
simultaneous over comparisons.

\begingroup\scriptsize\setlength{\tabcolsep}{3pt}
\begin{longtable}{llll}
\caption{Pairwise predictive comparisons for the 12-unit and 11-unit analyses. Negative mean NLL differences and positive information gains favor the first method. Micro information gain pools score differences and spike counts across units. Brackets give pointwise 95\% paired condition-bootstrap intervals.}\label{tab:allen-all-pairs}\\
\toprule
Units & First vs second method & \shortstack{Macro $\Delta$NLL\\95\% interval} & \shortstack{Micro bits/spike\\95\% interval} \\
\midrule\endfirsthead
\toprule
Units & First vs second method & \shortstack{Macro $\Delta$NLL\\95\% interval} & \shortstack{Micro bits/spike\\95\% interval} \\
\midrule\endhead
12 & DNN--Softplus vs Conditional kernel & -0.029006 [-0.040423, -0.017292] & 0.0078515 [0.0047153, 0.010883] \\
11 & DNN--Softplus vs Conditional kernel & -0.026096 [-0.037913, -0.013776] & 0.0068535 [0.0036527, 0.0099081] \\
12 & DNN--ReLU vs Conditional kernel & -0.020796 [-0.030918, -0.010262] & 0.0056291 [0.0028146, 0.0082415] \\
11 & DNN--ReLU vs Conditional kernel & -0.019545 [-0.029905, -0.0093311] & 0.005133 [0.0024795, 0.0077339] \\
12 & DNN--Softplus vs DNN--ReLU & -0.0082103 [-0.017285, 0.00092019] & 0.0022224 [-0.00025334, 0.0047148] \\
11 & DNN--Softplus vs DNN--ReLU & -0.0065513 [-0.016178, 0.0035349] & 0.0017206 [-0.00092349, 0.0043033] \\
12 & DNN--Softplus vs Pooled kernel & -0.41827 [-0.52087, -0.32659] & 0.11322 [0.090735, 0.13701] \\
11 & DNN--Softplus vs Pooled kernel & -0.38835 [-0.48738, -0.29823] & 0.10199 [0.080291, 0.12445] \\
12 & DNN--ReLU vs Pooled kernel & -0.41006 [-0.51219, -0.31736] & 0.11099 [0.088519, 0.13499] \\
11 & DNN--ReLU vs Pooled kernel & -0.3818 [-0.48114, -0.29125] & 0.10027 [0.078504, 0.12314] \\
12 & Conditional kernel vs Pooled kernel & -0.38926 [-0.4869, -0.29954] & 0.10536 [0.083373, 0.12835] \\
11 & Conditional kernel vs Pooled kernel & -0.36226 [-0.45632, -0.27519] & 0.095138 [0.073976, 0.11694] \\
\bottomrule\end{longtable}\endgroup

\begingroup\scriptsize\setlength{\tabcolsep}{3pt}
\begin{longtable}{llllll}
\caption{Unit-specific prediction scores for the four estimators. Mean NLL is in nats per test presentation; count deviance is summed over the 607 test presentations. Count MAE and RMSE are in spikes per presentation.}\label{tab:allen-unit-method}\\
\toprule
Unit & Method & Mean NLL & Count deviance & Count MAE & Count RMSE \\
\midrule\endfirsthead
\toprule
Unit & Method & Mean NLL & Count deviance & Count MAE & Count RMSE \\
\midrule\endhead
950930407 & DNN--ReLU & -32.091 & 733.63 & 2.8705 & 3.5959 \\*
950930407 & DNN--Softplus & -32.089 & 734.56 & 2.8668 & 3.596 \\*
950930407 & Conditional kernel & -32.072 & 742.38 & 2.8862 & 3.6215 \\*
950930407 & Pooled kernel & -32.009 & 791.48 & 2.9592 & 3.752 \\
950930985 & DNN--ReLU & -27.301 & 931.75 & 3.0347 & 3.7733 \\*
950930985 & DNN--Softplus & -27.302 & 929.87 & 3.0328 & 3.769 \\*
950930985 & Conditional kernel & -27.286 & 929.48 & 3.0395 & 3.7682 \\*
950930985 & Pooled kernel & -27.247 & 967.47 & 3.1041 & 3.8542 \\
950931181 & DNN--ReLU & -5.4646 & 2292.4 & 2.4919 & 3.5381 \\*
950931181 & DNN--Softplus & -5.5336 & 2236.2 & 2.246 & 3.4521 \\*
950931181 & Conditional kernel & -5.4841 & 2307.2 & 2.4948 & 3.5727 \\*
950931181 & Pooled kernel & -4.5327 & 3291.6 & 3.3046 & 4.4735 \\
950931254 & DNN--ReLU & -36.991 & 1055.6 & 3.6299 & 4.503 \\*
950931254 & DNN--Softplus & -36.954 & 1077.8 & 3.614 & 4.5745 \\*
950931254 & Conditional kernel & -36.937 & 1063.6 & 3.6721 & 4.5482 \\*
950931254 & Pooled kernel & -36.594 & 1243.1 & 3.9453 & 4.9961 \\
950931315 & DNN--ReLU & -6.707 & 1272.1 & 1.9907 & 2.6299 \\*
950931315 & DNN--Softplus & -6.7334 & 1254.6 & 2.0037 & 2.6105 \\*
950931315 & Conditional kernel & -6.6724 & 1303.2 & 2.0731 & 2.6873 \\*
950931315 & Pooled kernel & -5.9861 & 1952.2 & 2.6699 & 3.4474 \\
950931458 & DNN--ReLU & -5.6187 & 571.92 & 1.3685 & 1.763 \\*
950931458 & DNN--Softplus & -5.6294 & 562.17 & 1.3511 & 1.7467 \\*
950931458 & Conditional kernel & -5.6495 & 559.75 & 1.3422 & 1.757 \\*
950931458 & Pooled kernel & -5.441 & 734.37 & 1.5685 & 2.0679 \\
950931565 & DNN--ReLU & -6.6409 & 800.71 & 1.7104 & 2.1462 \\*
950931565 & DNN--Softplus & -6.6454 & 794.45 & 1.7068 & 2.1375 \\*
950931565 & Conditional kernel & -6.6396 & 791.37 & 1.7167 & 2.138 \\*
950931565 & Pooled kernel & -6.5978 & 838.68 & 1.7835 & 2.2117 \\
950931656 & DNN--ReLU & -3.1643 & 485.06 & 1.0063 & 1.2609 \\*
950931656 & DNN--Softplus & -3.1617 & 495.25 & 1.0088 & 1.2732 \\*
950931656 & Conditional kernel & -3.1637 & 493.63 & 1.0084 & 1.2855 \\*
950931656 & Pooled kernel & -2.9916 & 633.3 & 1.1747 & 1.5097 \\
950931727 & DNN--ReLU & -9.8035 & 471.29 & 1.4317 & 1.9249 \\*
950931727 & DNN--Softplus & -9.805 & 475.97 & 1.4351 & 1.9465 \\*
950931727 & Conditional kernel & -9.7698 & 477.8 & 1.4584 & 1.9598 \\*
950931727 & Pooled kernel & -9.3176 & 883.66 & 2.0054 & 2.7506 \\
950931770 & DNN--ReLU & -4.5335 & 552.61 & 1.2958 & 1.6014 \\*
950931770 & DNN--Softplus & -4.5413 & 532.31 & 1.231 & 1.5694 \\*
950931770 & Conditional kernel & -4.5402 & 540.2 & 1.2628 & 1.5851 \\*
950931770 & Pooled kernel & -4.4079 & 651.68 & 1.3496 & 1.7743 \\
950931899 & DNN--ReLU & -3.8556 & 556.82 & 1.1241 & 1.4201 \\*
950931899 & DNN--Softplus & -3.8596 & 549.86 & 1.1253 & 1.4085 \\*
950931899 & Conditional kernel & -3.8422 & 553.41 & 1.1316 & 1.4259 \\*
950931899 & Pooled kernel & -3.6331 & 776.51 & 1.416 & 1.7515 \\
950932563 & DNN--ReLU & -7.5589 & 1803.2 & 2.3126 & 2.9565 \\*
950932563 & DNN--Softplus & -7.5743 & 1792.8 & 2.3058 & 2.9641 \\*
950932563 & Conditional kernel & -7.4241 & 1953.7 & 2.5046 & 3.2094 \\*
950932563 & Pooled kernel & -6.0513 & 3390.6 & 3.5601 & 4.7588 \\
\bottomrule\end{longtable}\endgroup

\subsection{Numerical and modeling sensitivity}

\begin{figure}[p]\centering
\includegraphics[width=\linewidth]{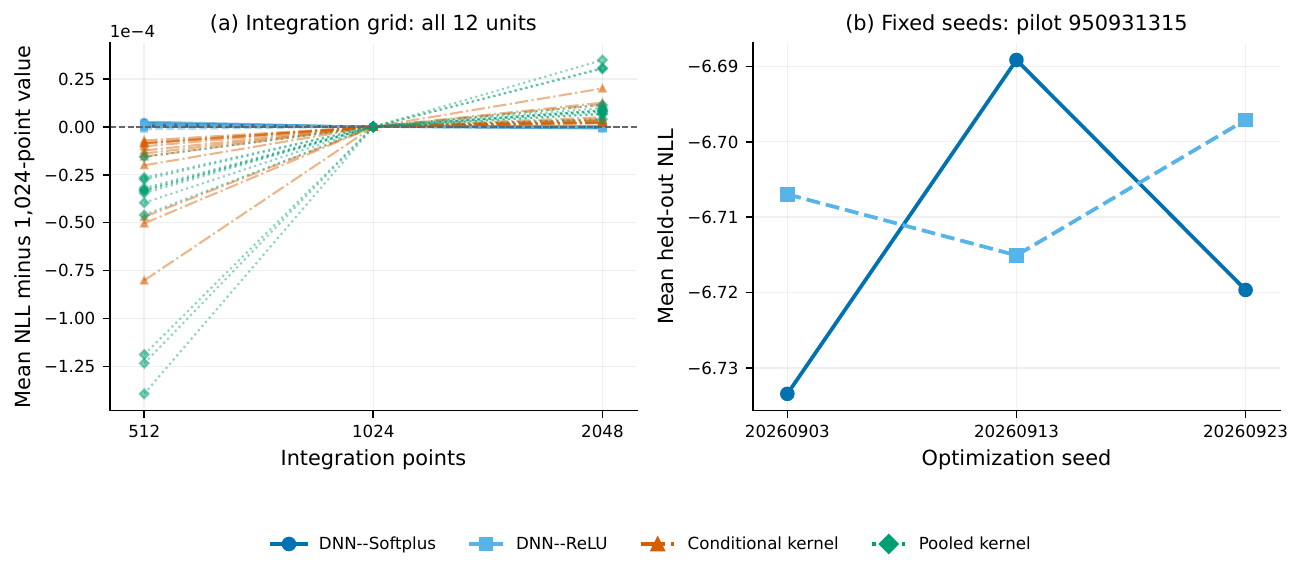}
\caption{Sensitivity to integration resolution and network initialization.
Left: changes in mean NLL from the 1,024-point score for each of the 48
unit--method fits. Right: mean test NLL at three optimization seeds
for unit 950931315.}
\label{fig:allen-sensitivities}
\end{figure}

We first checked the numerical integration used to score all 48
unit--method fits. Relative to the 1,024-point grid, grids of 512 and
2,048 points changed the mean NLL by at most $1.393\times10^{-4}$
and $3.481\times10^{-5}$, respectively
(Figure~\ref{fig:allen-sensitivities}). These changes were small compared
with the panel-average DNN--kernel differences.

The remaining analyses in this subsection used unit 950931315, the
previously examined pilot unit, and were exploratory. With optimization
seeds 20260903, 20260913, and 20260923, DNN--Softplus gave mean test
NLLs of $-6.73342$, $-6.68916$, and $-6.71965$; DNN--ReLU gave
$-6.70696$, $-6.71506$, and $-6.69711$. The ordering of the output
functions therefore depended on initialization. For Softplus networks
with widths 32, 64, and 128, the scores were $-6.70779$, $-6.73342$,
and $-6.71669$. The main results use width 64 and seed 20260903
throughout, without selection among these supplementary fits.

We also included mean signed running speed during the presentation as
an additional covariate. The resulting NLLs were $-6.72147$ for
DNN--Softplus, $-6.71394$ for DNN--ReLU, and $-6.51896$ for the
conditional kernel. Relative to the stimulus-only fits, the DNN changes
were small and differed in direction, whereas the conditional-kernel
score increased. The pooled estimator was unchanged because it does not
use covariates. Running speed was available for all retained
presentations; pupil area was unavailable and was omitted.

To examine prediction at new stimulus conditions, we partitioned the
243 position--orientation combinations into five folds and withheld all
presentations of each fold in turn. Normalization and epoch selection
were repeated using the other four folds. Both DNNs had lower NLL than
the conditional kernel in all five held-out folds
(Figure~\ref{fig:allen-condition-time}). This exploratory analysis repartitions all presentations, including
the primary test sample; it is not an independent replication.
Table~\ref{tab:allen-secondary} summarizes these comparisons.

\begingroup\scriptsize\setlength{\tabcolsep}{3pt}
\begin{longtable}{lllll}
\caption{Sensitivity analyses for unit 950931315. Mean NLL is in nats per test presentation, and count deviance is summed over the test presentations. The seed and width analyses include the corresponding primary fit for reference. Condition-fold scores use different held-out stimuli.}\label{tab:allen-secondary}\\
\toprule
Analysis & Method & Presentations & Mean NLL & Count deviance \\
\midrule\endfirsthead
\toprule
Analysis & Method & Presentations & Mean NLL & Count deviance \\
\midrule\endhead
Seed 20260903 & DNN--Softplus & 607 & -6.7334 & 1254.6 \\
Seed 20260903 & DNN--ReLU & 607 & -6.707 & 1272.1 \\
Seed 20260913 & DNN--Softplus & 607 & -6.6892 & 1318.5 \\
Seed 20260913 & DNN--ReLU & 607 & -6.7151 & 1269.9 \\
Seed 20260923 & DNN--Softplus & 607 & -6.7197 & 1262.4 \\
Seed 20260923 & DNN--ReLU & 607 & -6.6971 & 1275.9 \\
Width 32 & DNN--Softplus & 607 & -6.7078 & 1286.7 \\
Width 64 & DNN--Softplus & 607 & -6.7334 & 1254.6 \\
Width 128 & DNN--Softplus & 607 & -6.7167 & 1271.7 \\
Running speed & DNN--Softplus & 607 & -6.7215 & 1254.1 \\
Running speed & DNN--ReLU & 607 & -6.7139 & 1253.6 \\
Running speed & Pooled kernel & 607 & -5.9861 & 1952.2 \\
Running speed & Conditional kernel & 607 & -6.519 & 1412.7 \\
Condition fold 1 & DNN--Softplus & 620 & -6.8944 & 1507.3 \\
Condition fold 1 & DNN--ReLU & 620 & -6.8842 & 1509 \\
Condition fold 1 & Pooled kernel & 620 & -6.3286 & 1956.2 \\
Condition fold 1 & Conditional kernel & 620 & -6.8342 & 1517.2 \\
Condition fold 2 & DNN--Softplus & 610 & -5.4715 & 1477.9 \\
Condition fold 2 & DNN--ReLU & 610 & -5.4734 & 1474.6 \\
Condition fold 2 & Pooled kernel & 610 & -5.1138 & 1801.1 \\
Condition fold 2 & Conditional kernel & 610 & -5.4171 & 1497.4 \\
Condition fold 3 & DNN--Softplus & 595 & -6.6679 & 1382.4 \\
Condition fold 3 & DNN--ReLU & 595 & -6.6394 & 1408.3 \\
Condition fold 3 & Pooled kernel & 595 & -5.9596 & 2016.4 \\
Condition fold 3 & Conditional kernel & 595 & -6.6126 & 1443.8 \\
Condition fold 4 & DNN--Softplus & 597 & -8.3756 & 1327.9 \\
Condition fold 4 & DNN--ReLU & 597 & -8.39 & 1309.1 \\
Condition fold 4 & Pooled kernel & 597 & -7.4918 & 2140.2 \\
Condition fold 4 & Conditional kernel & 597 & -8.2677 & 1413.7 \\
Condition fold 5 & DNN--Softplus & 602 & -6.3423 & 1403.3 \\
Condition fold 5 & DNN--ReLU & 602 & -6.3328 & 1407.7 \\
Condition fold 5 & Pooled kernel & 602 & -5.644 & 1983.6 \\
Condition fold 5 & Conditional kernel & 602 & -6.2842 & 1455.4 \\
\bottomrule\end{longtable}\endgroup

\subsection{Fitted response profiles and model diagnostics}

\begin{figure}[tbp]
\centering
\includegraphics[width=\linewidth]{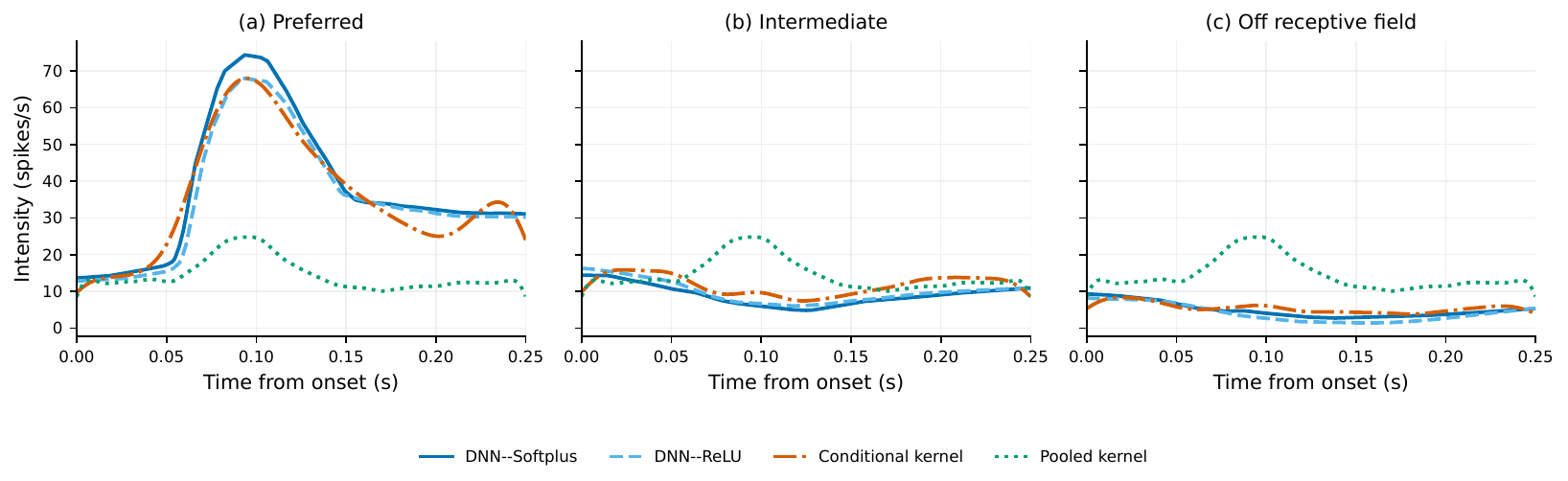}
\caption{Estimated temporal firing rates for unit 950931315 at preferred,
intermediate, and low-response stimulus conditions selected using
development counts. Time is measured from stimulus onset and intensity
is in spikes per second. Curves are fitted intensities.}
\label{fig:allen-intensities}
\end{figure}

Figure~\ref{fig:allen-intensities} displays fitted temporal intensities
for unit 950931315 at three stimulus conditions. We selected these
conditions using development counts to represent preferred, intermediate,
and low-response stimuli. The conditional estimates vary in both response
magnitude and temporal shape, whereas the pooled temporal kernel returns
one curve for all conditions.

\begin{figure}[p]\centering
\includegraphics[width=\linewidth]{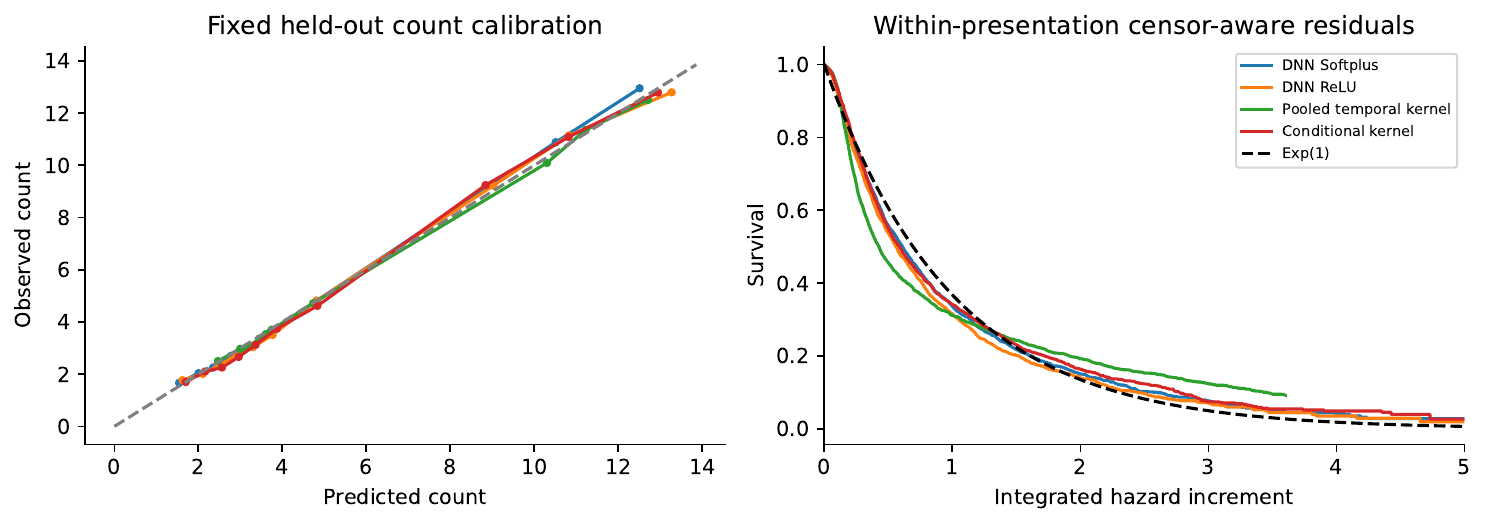}
\caption{Count and timing diagnostics. Left: observed and predicted counts
in method-specific deciles of predicted count, pooled over test
unit--presentation pairs. Right: Kaplan--Meier survival curves of
integrated-hazard increments for unit 950931315, with right censoring at
presentation end, and the unit-exponential reference.}
\label{fig:allen-diagnostics}
\end{figure}

For count calibration, we pooled test unit--presentation pairs and formed
deciles of predicted count separately for each method. The left panel of
Figure~\ref{fig:allen-diagnostics} compares average observed and
predicted counts in these groups. This display assesses aggregate count
prediction rather than calibration separately for every unit and
stimulus condition.

We assessed within-presentation timing for unit 950931315 by integrating
the fitted intensity on a 2,048-point grid. Transformed intervals ending
at a spike were treated as completed observations, and the remaining
interval to the end of each presentation was treated as right-censored;
for a zero-spike presentation, the whole window was censored. The right
panel of Figure~\ref{fig:allen-diagnostics} compares Kaplan--Meier
survival curves of these integrated-hazard increments with the
unit-exponential survival function. The visible deviations indicate
features of the spike trains not fully captured by the conditional
Poisson model, despite improvements in predictive NLL.

\begin{figure}[p]\centering
\includegraphics[width=\linewidth]{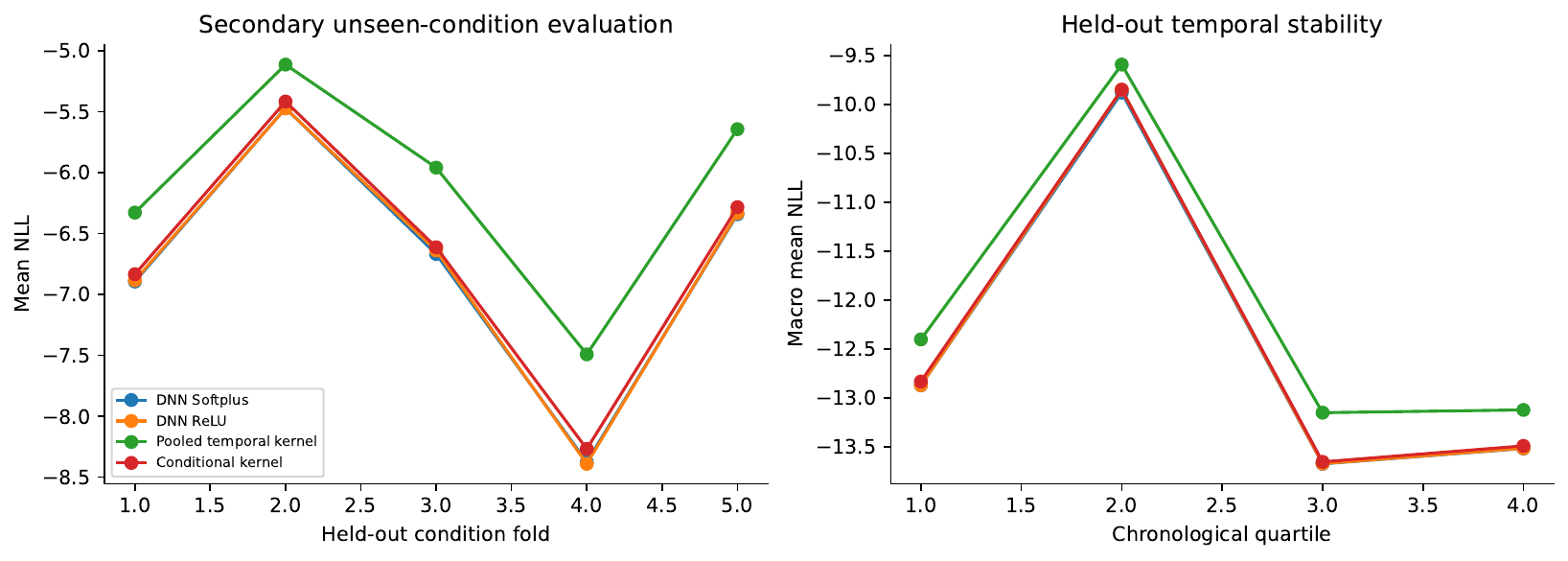}
\caption{Prediction across stimulus conditions and recording time. Left:
mean NLL for unit 950931315 in five condition-held-out folds. Right:
equal-unit mean NLL by chronological quartile of the primary test
presentations, using the original fits.}
\label{fig:allen-condition-time}
\end{figure}

Finally, we examined temporal variation using the original fitted models.
Figure~\ref{fig:allen-condition-time} reports scores by chronological
quartile of the test presentations. Differences in absolute NLL across
quartiles can arise from spike counts and stimulus composition as well
as predictive fit. We also calculated unit-level paired intervals using
noncircular moving blocks of 20 or 50 ordered test presentations, with
2,000 resamples and concatenated blocks truncated to the original sample
size. Because these blocks follow the ordered test subset rather than
contiguous presentations of the full recording, this analysis is a
sensitivity check for dependence within the test sample, not an
assessment of prediction beyond the training period.

\begin{figure}[p]\centering
\includegraphics[width=\linewidth]{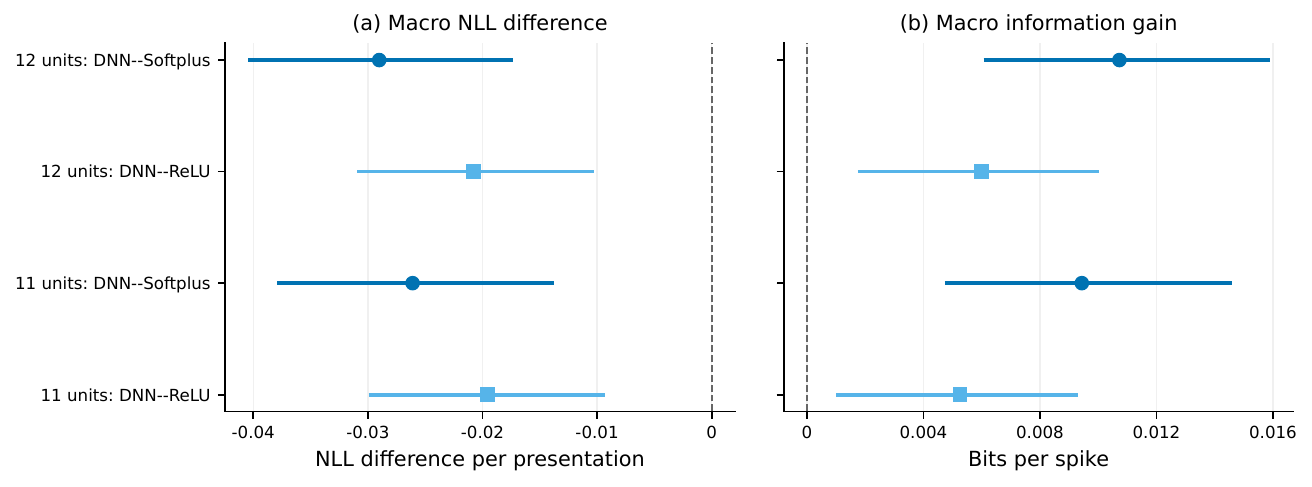}
\caption{Panel-average DNN performance relative to the conditional kernel
for all 12 selected units and the 11 units excluding the pilot. Negative
NLL differences and positive information gains favor the DNN. Segments
are pointwise 95\% paired condition-bootstrap intervals.}
\label{fig:allen-panels-supp}
\end{figure}

\clearpage
\bibliographystyle{plainnat}  
\bibliography{bibliography}  

\end{document}